\documentclass{times}

\usepackage{amsmath,amsfonts,amsthm,amscd}
\usepackage[all]{xy}
\usepackage{graphicx}
\usepackage{pinlabel}
\usepackage{cite}

\newtheorem{thm}{Theorem}[section]
\newtheorem{lem}[thm]{Lemma}
\newtheorem{prop}[thm]{Proposition}
\newtheorem{cor}[thm]{Corollary}

\newtheorem{conj}[thm]{Conjecture}
\newtheorem{add}[thm]{Addendum}

\newtheorem{hyp}[thm]{Hypothesis}
\newtheorem{MainTh}{Theorem}
\newtheorem{Not}[thm]{Notation}

\theoremstyle{definition}
\newtheorem{defn}[thm]{Definition}
\newtheorem{ex}[thm]{Example}

\theoremstyle{remark}
\newtheorem*{rmk}{Remark}
\newtheorem*{rks}{Remarks}

\newenvironment{pf}{ \begin{proof} }{ \end{proof} }

\DeclareMathOperator*{\U}{\mathrm{U}}
\DeclareMathOperator*{\SO}{\mathrm{SO}}
\DeclareMathOperator*{\SU}{\mathrm{SU}}
\DeclareMathOperator{\im}{im}
\DeclareMathOperator{\End}{\mathrm{End}}
\DeclareMathOperator{\Hom}{\mathrm{Hom}}
\DeclareMathOperator{\sym}{Sym}
\DeclareMathOperator{\dbar}{\bar{\partial}}
\DeclareMathOperator{\imag}{Im}
\DeclareMathOperator{\real}{Re}
\DeclareMathOperator{\aut}{Aut}
\DeclareMathOperator{\diff}{Diff}
\DeclareMathOperator{\ham}{Ham}

\DeclareMathOperator{\supp}{supp}
\DeclareMathOperator{\cone}{\mathsf{cone}}

\DeclareMathOperator{\coker}{coker}
\DeclareMathOperator{\gr}{gr}

\DeclareMathOperator{\ind}{ind}
\DeclareMathOperator{\fix}{Fix}

\DeclareMathOperator{\torus}{T}
\DeclareMathOperator{\interior}{int}
\DeclareMathOperator{\crit}{crit}
\DeclareMathOperator{\hor}{\EuScript{H}}

\DeclareMathOperator{\HF}{HF}
\DeclareMathOperator{\CF}{CF}

\DeclareMathOperator{\Fix}{Fix}
\DeclareMathOperator{\Ann}{Ann}
\DeclareMathOperator{\action}{\mathsf{A}}

\newcommand{\ii}{\mathrm{i}}
\newcommand{\id}{\mathrm{id}}
\newcommand{\Tv}{T^{\mathrm{v}}}
\newcommand{\Th}{T^{\mathrm{h}}}
\newcommand{\wh}[1]{\widehat{#1}}
\newcommand{\pr}{\mathsf{pr}}

\begin{document}
\title{A symplectic Gysin sequence}
\author{Timothy Perutz}
\date{July 2008}
\begin{abstract}We use the theory of pseudo-holomorphic quilts to establish a counterpart, in symplectic Floer homology, to the Gysin sequence for the homology of a sphere-bundle.  In a motivating class of examples, this `symplectic Gysin sequence' is precisely analogous to an exact sequence describing the behaviour of Seiberg--Witten monopole Floer homology for 3-manifolds under connected sum. 
\end{abstract}
\maketitle
\section{Introduction}
\subsection{Spherically-fibred coisotropic submanifolds}
This paper is about the following situation. We are given
\begin{itemize}
\item
compact symplectic manifolds $(M, \omega_M)$ of dimension $2n$ and $(N,\omega_N)$ of dimension $2(n-k) $, where $0\leq k \leq n$;
\item
a smooth fibre bundle $\rho\colon V \to N$ whose typical fibre $F$ is a $k$-dimensional manifold admitting a Morse function with precisely two critical points (so $F$ is homeomorphic to $S^k$);
\item
a smooth codimension-$k$ embedding $i\colon V \hookrightarrow  M$, such that
\begin{equation}\label{coiso}
i^*\omega_M=\rho^*\omega_N.
\end{equation}
\end{itemize}
We permit $M$ to have non-empty boundary provided that $i(V)\subset \interior(M)$ and that the boundary is of contact type with outward-pointing Liouville field. However, $N$ must be a closed manifold.

More economically, we can suppose that the data given are $(M,\omega_M)$ and its submanifold $V$, which must be a
\emph{spherically fibred coisotropic submanifold} (a notion which we shall define momentarily). The remaining data, $(N,\omega_N)$ and $\rho$, can then be reconstructed.  Recall that a coisotropic submanifold $V \subset M$ is one such that the annihilator $\Ann(TV)$ of $TV$ with respect to $\omega_M$ is contained in $TV$. As each point $x\in V$, $\omega_M$ then pairs $\Ann(T_xV)$ non-degenerately with $T_xM/T_x V$, hence $\Ann(TV)$ has constant rank $k$. Moreover, $V$ is foliated by isotropic leaves \cite[Lemma 5.33]{McDuffSal} with $k$-dimensional tangent spaces $\Ann(T_x V)$, and there is a quotient map $\rho$ from $V$ to the leaf-space $N$.
Say that $V $ is \emph{fibred coisotropic}\footnote{Or `integral coisotropic', in W. Ruan's terminology \cite{Ruan 05}.} if the isotropic foliation is fibrating: that is, each leaf $F$ has an open neighbourhood $U$ for which there is diffeomorphism $U \cong F \times D^k$ taking the isotropic foliation to the product foliation. When this happens, $N$ has induced smooth and symplectic structures, and $\rho$ becomes a smooth submersion satisfying (\ref{coiso}).

Given $(M,\omega_M;N,\omega_N;V,\rho)$ satisfying (\ref{coiso}),
there is a natural symplectomorphism of $(N,\omega_N)$ with the leaf-space (or reduced symplectic manifold) of $V$, intertwining $\rho \colon V \to N$ with the quotient map from $V$ to its leaf-space.

A fibred coisotropic submanifold is \emph{spherically fibred} when the quotient map to its leaf-space is a bundle whose fibres are homeomorphic to the sphere $S^k$.

One reason to regard fibred coisotropics as a natural class of submanifolds is that, unlike general coisotropic submanifolds, they have a straightforward deformation theory \cite{Ruan05}.  Modulo Hamiltonian diffeomorphisms, the infinitesimal deformations are parametrised by $H^0(N; \mathcal{H}^1 )$ where  $\mathcal{H}^1$ is the local system on $N$  associated with $(V, \rho)$ with fibres  $H^1 (\rho^{-1}(x);\R)$. 

A simple property of fibred coisotropic submanifolds is the following:
\begin{prop}
The graph $\wh{V} := \gr(\rho)$ is Lagrangian in $(M,-\omega_M) \times (N,\omega_N)$.
\end{prop}
\begin{pf}
If $(u_1,v_1)$ and $(u_2,v_2)$ lie in $T_{(x,\rho(x))}V \subset T_x M \times T_{\rho(x)}N$ then $v_1=\rho_*u_1$ and $v_2=\rho_*u_2$, hence 
$\omega_M(u_1,u_2)=(\rho^*\omega_N) (u_1,u_2) = \omega_N(v_1,v_2)$. Thus $\wh{V}$ is isotropic, and by a dimension-count it is Lagrangian. Alternatively, even in infinite dimensions, we can argue that if $(u,\rho_*u')\in \Ann(T_{(x,\rho(x))} \wh{V})$ then for all $w \in T_x V$ one has $\omega_M(u,w)=\omega_N(\rho_*u', \rho_* w)=(\rho^*\omega_N)(u',w)=\omega_M(u',w)$; hence $u-u'\in \Ann(T_x V)=\ker(D\rho)$, so $\rho_*u'=\rho_*u$ and $\wh{V}$ is coisotropic.
\end{pf}
\begin{Not}
We usually indicate that we wish to consider the symplectic form $-\omega_M$ on $M$, rather than $\omega_M$, by writing $M_-$. We distinguish fibred coisotropic submanifolds from their associated Lagrangian correspondences by adding a hat in the latter case, as in the previous proposition.
\end{Not}

The subject of this paper is the Lagrangian Floer homology
\begin{equation}\label{Floer group}
	\HF(\wh{V},(\id_{M_-} \times \mu)\wh{V}),
\end{equation}
taken inside $M_-\times N$, where $\mu\in \aut(N,\omega_N)$. Our aim is to establish a relationship between this and the fixed-point Floer homology $\HF(\mu)$. This takes the form of a long exact sequence, the `symplectic Gysin sequence' of the title. 

\subsection{Statement of results}\label{Statement}
To ensure that our Floer homology groups and the maps between them are well-defined we impose monotonicity hypotheses on the symplectic manifolds $M$ and $N$ and the Lagrangian correspondence $\wh{V}$.   A connected symplectic manifold $(M,\omega)$ is \emph{monotone} if the the symplectic  area homomorphism $\pi_2(M)\to \R$, defined by $f\mapsto \int_{S^2} {f^*\omega}$, is a non-negative multiple of the the Chern class homomorphism $f\mapsto\int_{S^2}{f^*c_1(TM,\omega)}$. A compact connected Lagrangian submanifold $L\subset M$ is monotone if the symplectic area $\pi_2(M,L)\to \R$ is a non-negative multiple of the Maslov index $m_L$. A symplectic manifold containing a monotone Lagrangian is necessarily monotone, since the image of $c_1(TM,
\omega)$ in $\Hom(\pi_2(M,L),\R)$ is $2m_L$. The \emph{minimal Maslov index} $m_L^{\min}$ of a Lagrangian is the non-negative generator of the additive subgroup $m_L(\pi_2(M,L))\subset \Z$.

\begin{hyp}[{\sc mon}]
$ \wh{V}$ is monotone in $M_- \times N$ with $m_{\wh{V}}^{\min} \geq 2$.
\end{hyp}

The hypothesis {\sc mon} is sufficient to guarantee that $\HF(\wh{V}, (\id_M \times\mu) \wh{V})$ is well-defined, with mod 2 Novikov coefficients (when $m_{\wh{V}}^{\min} =2$ this follows from the argument of \cite[Addendum]{Oh93}) but not strong enough to validate the symplectic Gysin sequence. We supplement it with the following:
\begin{hyp}[{\sc mas}]
$m_{\wh{V}}^{\min}\geq k + 2 $.
\end{hyp}

\begin{Not}\label{Novikov}
The universal Novikov ring $\Lambda_R$ of the commutative ring $R$ is the ring of formal series $\sum_{r\in \R}{a(r) t^r}$, with $a(r)\in R$, such that $\supp(a)\cap (-\infty,s]$ is finite for each $s\in \R$.
\end{Not}

\begin{MainTh}[{\sc Gysin sequence}]\label{Gysin}
Assume that {\sc mon} and {\sc mas} hold. For each $\mu \in \aut(N,\omega_N)$ there are then canonical maps forming an exact triangle
\[ \xymatrix{ %
	\HF(\mu) \ar[rr]^{e(V)\,  \frown \, \cdot} && \HF(\mu)
	\ar[dl]^{\rho^!} \\
	& \HF(\wh{V}, (\id_{M_-}\times \mu) \wh{V}) \ar[ul]^{\rho_*}
 }\]
Here $e(V) \in H^{k+1}(N;R)$ is the Euler class of the $S^k$-bundle
$\rho \colon V \to N,$ and the horizontal arrow is its action on fixed-point Floer homology by quantum cap product.  All three terms in the triangle are taken with coefficients in the universal Novikov ring $\Lambda_R$ of a commutative ring $R$ of characteristic 2.
\end{MainTh}
\begin{rks}
(i) All three maps in the sequence have a definite degree, but since Floer-theoretic gradings are relative, only the horizontal arrow has a straightforward \emph{numerical} degree. This degree is $-k-1$. If one declares that $\rho_*$ has degree $0$ (modulo $m^{\min}_{\wh{V}}$) then $\rho^!$ has degree $k$. Note that $(-k-1)+ 0 + k =-1$.

(ii) The hypotheses are sharp in that there is no odd $k$ for which {\sc mas} could be weakened to $m_{\wh{V}}^{\min}\geq k + 1$ (see Example \ref{counterex}). In Theorem \ref{borderline}, we analyse the borderline case, $m^{\min}_{\wh{V}}=k+1$. We show that the sequence is valid when the map $e(V)\,  \frown$ is replaced by $x\mapsto e(V) \frown x + \nu\,  t^a x$, where $a$ is the area of the discs of Maslov index $k+1$, $t$ the formal Novikov variable, and $\nu$ is the mod 2 count of index $k+1$ pseudo-holomorphic discs attached to $\wh{V}$, with two marked boundary points mapping respectively to a point and to a \emph{global angular chain}.
\end{rks}
\begin{add}\label{add1}
The symplectic Gysin sequence has the following additional properties.
\begin{itemize}
\item
The maps $\rho^!$ and $\rho_*$ intertwine the natural actions of the quantum cohomology $QH^*(N;\Lambda_{\Z/2})$.
\item
One can lift the characteristic two assumption and take, for example, $R=\Z$, providing that $(V,\rho)$ is an $R$-orientable bundle and that either $k=1$, or  there is a vector bundle $E\to M$ such that $w_1(E)=0$ and $w_2(E)|_V=w_2(\ker D\rho)$.
\end{itemize}
\end{add} 

If we are in a situation where $\HF_*(\wh{V},\wh{V})\cong H_*(\wh{V};\Lambda_R)$, we obtain from the case $\mu=\id_N$ an exact triangle which resembles the classical Gysin sequence:
\begin{equation}\label{classical} 
\xymatrix{ %
	H_*(N;\Lambda_R) \ar[rr]^{e(V)\frown_q  \cdot} && H_{*-k-1}(N;\Lambda_R) \ar[dl]  \\
	& H_{*-1}(V;\Lambda_R). \ar[ul]^{[-1]}
 }  \end{equation}
The three maps in the classical Gysin sequence are $e(V) \frown   \cdot$, $\rho^!$ and $\rho_*$ (where, of course, these have their classical meanings). We have stated that the horizontal map in our sequence is the quantum cap product  $e(V)\frown_q \cdot $, a deformation of the classical cap product $e(V)\frown \cdot $.  
\begin{add}\label{add2}
The map $\rho^!$ is a deformation of the arrow $\swarrow$ from diagram \ref{classical}. That is,  $\swarrow$ is the lowest order part of the map $\rho^!$; the remaining part has strictly positive Novikov weight.
\end{add}
The proofs of the addenda are routine. However, a further refinement remains unproven. One has $\HF_*(\mu) =\HF^*(\mu^{-1})$ and $\HF_*(\wh{V}, (\id\times \mu)\wh{V})=\HF^*(\wh{V},\id\times \mu^{-1}\wh{V})$ (`Poincar\'e duality').  Conjecture: \emph{under Poincar\'e duality, $\rho^!$ (for $\mu$) corresponds to the adjoint of $\rho_*$ (for $\mu^{-1})$.} We will see that $\rho^!$ is defined geometrically, and is associated a picture (Figure \ref{gysincyl}), while $\rho_*$ is defined algebraically. The conjecture is that $\rho_*$ should be associated with a similar picture, with the orientation reversed.

There is another hypothesis under which one has an exact Gysin sequence: that of \emph{strong negativity}.
\begin{defn}\label{strongly negative}
A symplectic manifold $(X^{2n},\eta)$ is \emph{strongly negative} if every class $\beta\in \pi_2(X)$ with $\langle \beta, [\eta]\rangle >0 $ satisfies $\langle \beta, c_1(TX)\rangle \leq  -n-2$. A compact Lagrangian $L$ is strongly negative if every $\delta\in \pi_2(X,L)$ with $\langle \delta, [\eta] \rangle >0$ satisfies $m_L(\delta) \leq - n - 2$.
\end{defn}

\begin{thm}
The terms in the Gysin sequence are well-defined when $M$ and $N$ are strongly negative and $\wh{V}\subset M_-\times N$ is a strongly negative Lagrangian. Moreover, Theorem \ref{Gysin} and its two addenda continue to hold in this situation.\end{thm}

\subsection{Examples and corollaries}
The basic example of a spherically-fibred coisotropic submanifold is the following.
\begin{ex}\label{sphere ex}
Let $S=S^{2n-1}(r)$ be a sphere of radius $r$ in (standard symplectic) $\C^n$. It is coisotropic, and the isotropic leaves are the orbits under the scaling $S^1$-action; thus the reduced manifold is $\mathbb{CP}^{n-1}$ and the reduction map $S \to \mathbb{CP}^{n-1}$ is the Hopf fibration. The symplectic form on $\mathbb{CP}^{n-1}$ is $r \omega_{\mathrm{FS}}$. Thus there are Lagrangian spheres in $M_- \times \mathbb{CP}^{n-1}$, for any symplectic manifold $M^{2n}$ (we have to take $r$ small enough to find a Darboux chart of radius $r$). 
\end{ex}
Two mechanisms for producing fibred coisotropic
submanifolds are \emph{reduction} and \emph{degeneration}.

\emph{Symplectic reduction.} If the compact group $G$ acts by Hamiltonian diffeomorphisms of $M$ via a moment
map $\mu$, and $\mu^{-1}(s)$ is a regular level set on which $G$ acts freely, then $\mu^{-1}(s)$ is a fibred coisotropic with
$G$ as fibre. In particular, when $G$ is $\U(1)$ or $\SU(2)$, $\mu^{-1}(s)$ is spherically fibred. 

\emph{Degeneration.} Consider a symplectic manifold $(E,\Omega)$ and a proper map $\pi \colon E \to S$ to a Riemann surface $S$. Assume that the critical set $N:=\crit(\pi)$ is smooth and contained in some fibre $E_s$; that there is an almost complex structure $J$ in a neighbourhood of $N$ making $\pi$ holomorphic (so that $N$ is an almost complex submanifold) such that 
the Hessian form $D^2 \pi$ on its normal bundle,  $\nu=\nu_{N/E}$, is non-degenerate as a complex quadratic form on each fibre $\nu_x$. Let $\gamma \colon [0,1] \to S$ be a smooth embedding with $\gamma(1)=s$, and put $M=E_{\gamma(0)}$.
Parallel transport over $\gamma|_{[0,t]}$ defines a symplectomorphism $\rho_t \colon M= E_{\gamma(0)}\to E_{\gamma(t)}$,
for any $t\in [0,1)$. The \emph{vanishing cycle}
\[V_\gamma:= \{ x \in M: \lim_{t \to 1}{\rho_t(x)} \in N\}\]
is then a spherically-fibred coisotropic submanifold. The fibres are copies of $S^k$, where $2k+2$ is the real codimension of $N$ in $E$,
and the reduction map is the limiting parallel transport $\rho=\lim_{t\to 1}{\rho_t}\colon V_\gamma \to N$. The bundle $\rho$
comes with a natural reduction of structure to $\SO(k+1)$. Conversely, any  spherically-fibred coisotropic submanifold with structure group reduced to $\SO(k+1)$ arises in this way. For further details, see \cite[Section 2]{PerutzI}.

\begin{ex}
Let $(M,\omega)$ be a compact symplectic manifold, and fix $p\in M$. Let $\widetilde{M\times \C}$ be a symplectic blow-up of $M\times \C$ near $(p,0)$. Then the blow-down map  $\widetilde{M\times \C}\to M\times \C$ followed by second projection defines a degeneration $\widetilde{M\times \C}\to \C$ of the type just described (with $k=1$). The smooth fibres are symplectomorphic to $M$; the zero-fibre is a normal crossing variety $\widetilde{M}\cup_{\mathbb{CP}^{n-1}}{\mathbb{CP}^{n}}$, so its singular locus $N$ is a copy of $\mathbb{CP}^{n-1}$. The normal bundle $\nu$ is $\mathcal{O}(1)\oplus \mathcal{O}(-1)$, and the Hessian form is the natural pairing of $\mathcal{O}(1)$ with $\mathcal{O}(-1)$. The vanishing cycle (for the path $[-1,0]$) is just the sphere $S^{2n-1}\subset M$ of  Example \ref{sphere ex}. Its radius in a Darboux chart is the weight of the blow-up, i.e., the area of a projective line $\mathbb{CP}^1$ in the exceptional divisor.
\end{ex}

\begin{rmk}
(i) The `complex Morse--Bott'-type degenerations described above occur routinely in algebraic geometry. The case $k=1$ is
especially prevalent, because this is the situation which occurs when the singular fibre $E_s$ has a smooth normal crossing
singularity. At the end of the paper we shall take up the examples studied in  \cite{PerutzI}, in which $M$ is a symmetric product of a Riemann surface.

(ii) The symmetric product examples just mentioned arise in a context related to 3- and 4-dimensional topology. The Floer homology groups we shall consider, as well as the maps between them, appear as part of a conjectural symplectic model for Seiberg--Witten gauge theory (see \cite{PerutzI}). We shall explain that, for this class of examples, our exact sequence is precisely analogous to a connected sum formula in Seiberg--Witten monopole Floer homology.
\end{rmk}
Taking $\mu=\id_{N}$, the symplectic Gysin sequence gives a method for computing Lagrangian Floer \emph{self}-homology. The standard `PSS' isomorphism  $\HF(\id_N)\cong H_*(N;\Lambda_R)$ respects the actions on the two sides by the (small) quantum cohomology ring of $N$ \cite[Section 12.1]{McDuffSal04}. One can compare Theorem \ref{Gysin} with a result of Oh about Lagrangian Floer self-homology (which in the monotone case may be defined over any commutative ring $A$ of characteristic $2$ without the need for Novikov rings):
\begin{thm}[Oh \cite{Oh96}]
If $L$ is a monotone Lagrangian with $m_{L}^{\min}\geq \dim(L) + 2$, one has $\HF(L,L;A) \cong H_*(L; A)$.
\end{thm}

If $N$ is a point, so that $\dim(V)=n=k$, hypothesis {\sc mas} coincides with Oh's $m^{\min}_{\wh{V}}\geq \dim(V)+2$. On the other hand, if $N$ is positive-dimensional, {\sc mas} is a weaker assumption than Oh's.  Our result implies that if the quantum cap product on $N$ is undeformed, or if $e(V)=0$, one has $\HF_*(\wh{V},\wh{V};\Lambda_R) \cong H_*(V;\Lambda_R)$.

\begin{ex}\label{counterex}
Consider the spherical Lagrangian $\wh{S}=S^{2n-1}\subset (\C^n)_-\times \mathbb{CP}^{n-1}$ from Example \ref{sphere ex}. One has $m_{\wh{S}}^{\min}=2n = \dim(S)+1$. But $\wh{S}$ is Hamiltonian displaceable, by $(\mathrm{translation})\times \id_{\mathbb{CP}^{n-1}}$, so $\HF_*(\wh{S},\wh{S})=0$. If we take $M= D^{2n}(0;R)_-\times \mathbb{CP}^{n-1}$ with $R\gg 0$, and $N=\{\mathrm{pt.}\}$, {\sc mas} is, reassuringly, not fulfilled. One sees from this example that hypothesis {\sc mas} is sharp when $k=2n-1=\frac{1}{2}\dim(M)$. 

If in the same example, with $n\geq 2$, we take $M$ to be the $2n$-disc $D^{2n}(0;R)$ and $N=\mathbb{CP}^{n-1}$, then $k+2=3 < 2n$, so the Gysin sequence \emph{does} apply. The Euler class $e$ is the hyperplane class $h$, and quantum cap product by $h$ is invertible; indeed, the quantum cap operator $H=h\frown_q\cdot $ on $H_*(\mathbb{CP}^{n-1};\Lambda_R)$ has as inverse $t^{-1}H^{n-1}$, where $t\in \Lambda_R$ is the formal Novikov parameter (see, e.g., \cite{McDuffSal}). Thus Theorem \ref{Gysin} correctly predicts that $\HF_*(\wh{S},\wh{S})=0$ (or, reversing one's viewpoint, it shows that in $QH^*(\C P^{n-1})$ one has a relation $h^{*n}= f(t)$ for some $f\neq 0$).
\end{ex}
The last example points to a corollary of Theorem \ref{Gysin}.
\begin{cor}
When {\sc mon} and {\sc mas} hold, and $(V,\rho)$ is $R$-orientable, $\wh{V}$ is Hamiltonian displaceable in $M_-\times N$ only if $e(V)$ acts invertibly on the quantum homology $QH_*(N;\Lambda_R)$.
\end{cor}
Indeed, this is the criterion for the vanishing of $\HF_*(\wh{V},\wh{V})$. The unusual aspect of this corollary is that it invokes only \emph{closed} holomorphic curves.

It has recently been shown by Biran and Cornea \cite{BirCor} that, for a large class of monotone Lagrangians $L$, there is a dichotomy: either the Floer self-homology (over $\Lambda_{\Z/2}$) vanishes, or else it is isomorphic to $H_*(L;\Lambda_{\Z/2})$. Their hypotheses are that $m_L^{\min} \geq 2$ and that the cohomology ring $H^*(L;\Z/2)$ is generated in degrees $\leq m_L^{\min} - 1$. This, together with our vanishing criterion, has the following corollary.

\begin{thm}
Suppose that the hypotheses of Theorem \ref{Gysin} hold, and that $H^*(\wh{V};\Z/2)$ is generated in degrees $\leq m_{\wh{V}} - 1$. Then either $e(V)$ acts invertibly on $QH_*(N;\Lambda_{\Z/2})$, or $\HF_*(\wh{V},\wh{V};\Lambda_{\Z/2})\cong H_*(\wh{V};\Lambda_{\Z/2})$.
\end{thm}

\subsection{Remarks on the proof}
One approach to Theorem \ref{Gysin}, along the lines of a standard proof of its classical counterpart,  would be to study the spectral sequence associated with a natural filtration on $\CF_*(\wh{V}, (\id_{M_-}\times \mu) \wh{V})$. A generator for this complex projects to a fixed point of $\mu$, and the filtration is by the  Floer-theoretic (relative) degree of this fixed point, considered as a generator of $\CF_*(\mu)$. The task is to identify the $E^2$-page and the higher differentials.

We shall proceed by a different (more geometric?) route, exhibiting the sequence as the exact triangle associated with
a mapping cone. Under hypothesis {\sc mon}, we shall construct a chain map
\begin{equation}
  C\rho \colon \CF_*(\mu) \to \CF_*(\wh{V}, (\id_{M-}\times \mu) \wh{V}). 
  \end{equation}  
When {\sc mas} also holds the composite, quantum cap $e(V)\frown \cdot $ followed by $C\rho$, is nullhomotopic; by defining a nullhomotopy we shall construct a chain map 
\begin{equation}\label{iso}
 C\Phi  \colon  \cone \left( \CF_*(\mu)\stackrel{e(V)\frown \cdot }{\longrightarrow}
 	\CF_{*-k-1}(\mu)\right) \to
	\CF_*^{M_- \times N}(\wh{V}, (\id_{M-}\times \mu) \wh{V}).
\end{equation}
\begin{MainTh}\label{Quasiiso}
The map $C \Phi$ induces an isomorphism $\Phi$ on homology.
\end{MainTh}
Theorem \ref{Gysin} is an immediate consequence, since the homology of the mapping cone sits inside a canonical exact triangle. 
 
The classical Gysin sequence can be proved by considering the chain map
\begin{align} 
 C\Psi\colon \cone \left( \mathsf{S}_*(N) \stackrel{\tilde{e}\frown \cdot }{\longrightarrow}  \mathsf{S}_{*}(N)[-k-1] \right) & \to \mathsf{S}_{*}(V)[-k-1],\\
 \notag
 (\sigma, \sigma' ) &\mapsto   \psi\frown \rho^! \sigma + \rho^! \sigma'. \end{align}
Here $\mathsf{S}_*(N)$ and $\mathsf{S}_*(V)$ are singular chain complexes,  $\tilde{e}$ a cochain representative for $e(V)$, and $\psi\in \mathsf{S}^k (V)$ a cochain satisfying $d\psi = \rho^*\tilde{e}$ and restricting to each fibre as an orientation cocycle. (In de Rham theory, $\psi$ is what Bott and Tu \cite{BottTu} call a `global angular form'). One proves that $C\Psi$ is a quasi-isomorphism.

Our map $C\Phi$ can be viewed as a Floer-theoretic analogue of $C\Psi$. Its natural habitat is the  theory of `Lagrangian matching conditions' in the sense of \cite{PerutzI}, alias the  `quilted' Floer theory of Wehrheim--Woodward \cite{WehrWood}. To prove that $\Phi$ is an isomorphism we show that its low-action part defines a linear isomorphism and then fall back on the algebra of $\R$-graded complexes. We introduce an algebraic criterion (Lemma \ref{double}) which is applicable even though there is no monotonicity hypothesis on $\mu$.

To close the discussion, we remark that there are also Gysin-type sequences entirely within Lagrangian Floer homology, and it may well be that a variant of our result can be sharpened to a statement about of Fukaya categories. Assume that we have a rigorous formulation of the monotone Fukaya $A_\infty$-category $\EuScript{F}_N$, along the lines of that of the exact Fukaya category in \cite{Seidel06}. Fix a monotone Lagrangian $L\in \mathrm{Ob}(\EuScript{F}_N)$. There will then be an $A_\infty$-module over $\EuScript{F}_N$ \cite[Chapter 1]{Seidel06}, $\EuScript{C}\mathsf{one} (e(V)\smile \cdot)$, sending $L'\in \mathrm{Ob}(\EuScript{F}_N) $ to the chain complex $\cone(e(V)\smile \cdot \colon \CF^*(L',L)\to \CF^{*+k+1}(L',L))$; and a second $A_\infty$-module sending $L'$ to the Floer complex 
$\CF^*(L'\times \wh{V}, \wh{ V}^{\circ}\times L)$ taken in $N \times M_- \times N$,  where $\wh{V}^\circ$ denotes the `opposite' correspondence to $\wh{V}$. (For this we need $L$, $L'$ and $\wh{V}^\circ$ to be normalised monotone; see Section \ref{Lag match}.)
\begin{conj}
Under the hypotheses of Theorem \ref{Gysin}, these two $A_\infty$-modules are quasi-isomorphic.
\end{conj}
Our proof readily adapts to show that $\HF^*(L'\times \wh{V}, \wh{ V}^{\circ}\times L)\cong H\cone (e(V)\smile\cdot)$, so the substance of the conjecture lies in the compatibility with the $A_\infty$-structure.

By the main theorem of \cite{WehrWood}, $\HF^*(L'\times \wh{V}, \wh{ V}^{\circ}\times L) \cong HF^*_{M_-}(\rho^{-1} L' ,\rho^{-1}L)$. As the present paper was nearing completion, the author learned that the relation of $\HF^*(\rho^{-1}L,\rho^{-1} L)$ to $\HF^*(L,L)$ has been analysed directly, in the case $k=1$,  by Biran--Khanevsky \cite{BirKhan}. They focus on circle-bundles appearing as the boundary of a tubular neighbourhood of symplectic hypersurfaces, and proceed by a direct comparison of the `pearl complexes' computing Floer self-homology. The results appear to be similar to ours where both apply, but the motivations, methods and applications are quite different.

\subsection{Navigation}
Section \ref{review} reviews the basic Floer-theoretic constructions, as well as some basic facts about symplectic fibrations. Section \ref{Lag match},  which is largely to blame for the extravagant length of the paper, explains the framework of `Lagrangian matching conditions'. There is more generality here than is really needed, but this framework undoubtedly has further applications. Readers who find such formalism disagreeable are urged to press on to Sections \ref{geom} and \ref{alg}, wherein we prove Theorem \ref{Gysin}. In Section \ref{refinements} we study the `borderline' case, and briefly consider orientations. Section \ref{3manifolds} explains the connection with Seiberg--Witten gauge theory and the monopole Floer homology of connected sums which provided the initial impetus for this paper.  

\subsection{Acknowledgements}
My thanks to Robert Lipshitz, Ciprian Manolescu, Dusa McDuff,  Ivan Smith, Katrin Wehrheim and Chris Woodward for stimulating discussions, and to MIT for generous hospitality. I acknowledge financial support during the initial stages of this research from EPSRC research grant EP/C535995/1.

\section{Review: fibre bundles, Floer homology}\label{review}

\subsection{Locally Hamiltonian fibre bundles}
\begin{defn}\label{LHF}
A \emph{locally Hamiltonian fibration} (LHF) is a triple $(E,p,\Omega)$ consisting of a smooth, proper fibre bundle $p\colon E \to B$ equipped with
a closed two-form $\Omega \in \Omega^2(E)$ which restricts to each fibre $E_b$ as a non-degenerate (hence symplectic) form.
\end{defn}
An LHF carries a canonical connection $A=A(\Omega)$. The horizontal distribution $\Th E$ is given by $\Th_x E  = \Ann^{\Omega}(\Tv_x E)$: that is, a tangent vector $u$ is horizontal iff $\Omega(u,v)=0$ for all vertical vectors $v\in 
\Tv E:=\ker Dp$.
This connection is symplectic in the sense that $(\mathcal{L}_X\Omega )|_{E_b} =0 $ when $X$ is a horizontal vector field defined in a
neighbourhood of $E_b$. The monodromies of loops based at $b$ are therefore
symplectic; so, writing $(M,\omega_M)=(E_b,\Omega|_{E_b})$, and assuming $B$ to be connected, the structure group is reduced to
$\aut(M,\omega_M)$. The monodromies of null-homotopic loops are in fact Hamiltonian automorphisms \cite[Chapter 6]{McDuffSal}.
\begin{ex}\label{mapping torus}
A fundamental example is the \emph{mapping torus} $(\torus_\mu,p_\mu,\omega_\mu)$ of an automorphism $\mu \in \aut(M,\omega)$.
Here $\torus_\mu$ is the quotient of $M \times \R$ by the relation $(\mu(x),t)\sim (x,t+1)$.
The projection map $p_\mu \colon \torus_\mu \to S^1$ is an $M$-bundle. The closed, fibrewise-symplectic two-form
$\omega_\mu$ is the form which pulls back to $M \times \R$ as $\pr_1^*\omega$. Any locally Hamiltonian fibre bundle over $S^1$ is canonically isomorphic to the mapping torus of its monodromy.
\end{ex}
Once we progress to higher-dimensional bases, the connection $A$ is usually not flat. Its curvature $F_\Omega$ is a two-form on $B$ with values in the bundle $\EuScript{V}$ of symplectic vector fields on the fibres, defined via horizontal lifting $X\mapsto X^h$ of vector fields on $B$
in the usual way:
\begin{equation} \label{curvature}
F_\Omega \in \Omega^2(B; \EuScript{V}); \quad  F_\Omega (X,Y) = [X^h,Y^h]-[X,Y]^h.
\end{equation}
The function $\kappa_{X,Y} = \Omega(X^h,Y^h)$ on $E$ has the property that
\begin{equation}
d\kappa_{X,Y} = \mathcal{L}_{X^h}\iota(Y^h) \Omega - \mathcal{L}_{Y^h}\iota(X^h) \Omega =
\iota([X^h,Y^h])\Omega;
\end{equation}
hence $\kappa_{X,Y}|_{E_b}$ is the Hamiltonian generator for the vertical vector field
$F_\Omega(X,Y)_b$. We can therefore understand $F_\Omega$ through the map $(X,Y)\mapsto \kappa_{X,Y}$, which is a two-form on $B$
with values in functions on the fibres.  One may in particular speak of \emph{flat} LHFs, meaning those for which $F_\Omega=0$, or equivalently, those for which the function-valued two-form $(X,Y)\mapsto \kappa_{X,Y}$ on $B$ reduces to an ordinary two-form. The techniques of \cite[Chapter 6]{McDuffSal} readily yield the following.
\begin{lem}\label{flat}
A flat LHF $(E,p,\Omega)$ over the cylinder $S^1\times I$, with fibre $(M,\omega)=(E_{(z,0)},\Omega|_{E_{(z,0)}})$, is equivalent to $( T_\mu \times I , p_\mu \times \id,  \pr_1^*\omega_\mu+ \beta)$, where $\mu$ is the monodromy around $S^1\times \{0\}$ and $\beta$ is the pullback of a form on the base.
\end{lem}
If $u\colon B\to E$ is a horizontal section, the pullback bundle $u^*\Tv E \to B$ inherits a symplectic connection $\nabla^u$, the result of linearizing $A_\Omega$: as noted in \cite[Section 2]{Seidel03}, one has $\nabla^u_X Z = u^* [X^h, \tilde{Z}]$, where $X$ is a vector field on $B$ and $\tilde{Z}$ any section of $\Tv E$ such that $u^*\tilde{Z}=Z$. The curvature $F(\nabla^u) \in \Omega^2_B(\mathfrak{sp}(u^* \Tv E))$  is given by
\begin{equation}\label{curvature2} 
F(\nabla^u) (X,Y)\colon Z\mapsto  u^* [ F_\Omega(X,Y), \tilde{Z} ] \in C^\infty_B(u^*\Tv E).   
\end{equation} 
\subsection{Floer homology for fixed points and Lagrangian intersections}
There are two variants of symplectic Floer homology:
\begin{itemize}
\item Floer homology for fixed points of symplectic automorphisms (not necessarily isotopic to the identity); and
\item Floer homology for intersections of Lagrangian submanifolds.\footnote{The first can be seen as a case of the second, since the Floer homology of an automorphism $\phi$ is canonically isomorphic to that of the diagonal $\Delta \subset M_-\times M$ with the graph of $\phi$, by an isomorphism that essentially goes back to Floer. However, the author does not know a source in the literature that explains this carefully. Note that the graph may not be monotone.}
\end{itemize}
Both are modules over the universal Novikov ring 
$\Lambda_R$ of a commutative ring $R$, as defined at (\ref{Novikov}). Until further notice, we assume that $R$ has characteristic 2.

\subsubsection{Fixed-point Floer homology}
Suppose that $(T,p,\Omega)$ is a LHF (\ref{LHF}) over a closed oriented 1-manifold $Z$. We suppose the typical fibre $N$ to be compact and either monotone or strongly negative (see (\ref{strongly negative})).

To $(T,p,\Omega)$ we attach its Floer homology $\HF(T,\Omega)$.
When $Z$ is connected, and $(N,\omega_N) = (E_z,\Omega|_{E_z})$
for a basepoint $z$, the monodromy along $Z$ is an automorphism $\mu\in \aut(N,\omega_N)$, and there is a
two-form-preserving bundle isomorphism of $T$ with $T_\mu$. We often write $\HF(\mu)$ instead of $\HF(T,\Omega)$.

$\HF(T,\Omega)$ is initially defined when $Z$ is connected and the monodromy $\mu$ has non-degenerate fixed points. It is then the homology of a complex $\CF(T,\Omega;J^v)$, where the underlying  $\Lambda_{R}$-module $\CF(T,\Omega)$ is the $\Lambda_R$-module freely generated by the set $\hor(T,\Omega)$ of \emph{horizontal sections} of $T \to Z$ (a set in bijection with $\Fix(\mu)$). The differential $\partial_{J^v}$ depends on a choice of $J^v \in \EuJ(\Tv T)= \EuJ(\Tv T,\Omega)$, the set of $\Omega$-compatible almost complex structure in the vertical tangent bundle $\Tv T$.  The complex also depends on an orientation-preserving diffeomorphism $Z\cong S^1$. 

Let $\pr_2 \colon \R \times S^1\to S^1$ be the second projection, and form the LHF $(\pr_2^* T,  p \circ \pr_2, \pr_2^*\Omega)$ over $\R\times S^1$ (we use the parametrisation $Z\cong S^1$ here). Any $J^v\in \EuJ(\Tv T)$ extends to a unique almost complex structure $\tilde{J^v}$ in $p^*T$ such that (i) $\tilde{J^v}$ is translation-invariant; (ii) $p \circ \pr_2$ is $\tilde{J^v}$-holomorphic (where the cylinder is identified with $\C/(z\sim z+ i)$ and inherits the latter's complex structure); and (iii) $\tilde{J^v}$ preserves the horizontal distribution. One then considers the moduli space $\EuM(T ; J^v)$ of $\tilde{J^v}$-holomorphic sections $u\colon \R\times S^1 \to \pr_2^*T$ of finite \emph{action} 
\begin{equation}\label{action} 
\action(u)=\int_{\R\times S^1}{u^* (\pr_2^*\Omega)}.
\end{equation}
We give it the $L^p_{1,loc}$ topology, for some $p>2$. When $T$ is non-degenerate, any $u\in \EuM(T ; J^v)$, thought of as a path of sections $u(s,\cdot)$ of $T$, converges uniformly as $s\to \pm\infty$ to horizontal sections $\nu_\pm \in \hor(T,\Omega)$ \cite{Salamon94}. Hence $\EuM(T;  J^v)$ is partitioned into components $\EuM(T; \nu_-, \nu_+; J^v)$. The \emph{Conley--Zehnder index} $\ind \colon \EuM(T; \nu_-, \nu_+; J^v)\to \Z$, which coincides with the Fredholm index of the linearisation of the Cauchy--Riemann equation $\tilde{J^v} \circ Du = Du \circ i$, further partitions $ \EuM(T; \nu_-, \nu_+; J^v)$ into subsets $\EuM(T; \nu_-, \nu_+; J^v)^q:=\ind^{-1}(q)$. 

\begin{defn}\label{regular}
We define a Baire subset $\EuJ^{\mathrm{reg}}(\Tv T)\subset \EuJ(\Tv T)$ of \emph{regular} complex structures as follows. When $N$ is monotone, $J^v\in \EuJ^{\mathrm{reg}}(\Tv T)$ if the moduli spaces  $\EuM(T; \nu_-, \nu_+; J^v)^q$ are all transversely cut out (in particular, empty when $q\leq 0$). In the strongly negative case, we require  additionally that   there should be no non-constant $J^v$-holomorphic spheres in the fibres of $T$. 
\end{defn}
\begin{rmk} Baire subsets of Fr\'echet spaces need not be dense. To get around this, we can make an initial choice of $J\in \EuJ(\Tv T)$, and then define regularity among those $J'\in\EuJ(\Tv T)$ which are close to $J$ in Floer's $C^\infty_\epsilon$ norm \cite{HS}; the norm depends on a choice of rapidly decreasing sequence $\epsilon=(\epsilon_m)_{m=1,2,\dots}$ and on a metric on $\Tv T$. This detail will be glossed over from now on.
\end{rmk} 
For a proof that  transversality of the moduli spaces is generic see, e.g., \cite{HS}. To understand why non-existence of $J^v$-spheres is  generic in the strongly negative setting,  look at the moduli space $\EuScript{N}^s(T)$ of non-constant, simple, parametrised $J^v$-spheres $\beta\colon S^2\to T$, mapping to some (variable) fibre of $p$. The real, fibrewise index $l= 2 \int_{S^2}{\beta^*c_1(TN)}  + \dim(N)$ breaks this moduli space into components  $\EuScript{N}^s(T)^l$. Impose the (generic) condition on $J^v$ that each $\EuScript{N}^s(T)^l$ is transversely cut out. The quotient of $\EuScript{N}^s(T)^l$ by $\aut(S^2)=\mathrm{PSL}(2,\C)$ is then a smooth manifold of dimension $l +1 - \dim \aut(S^2)=l-5$ \cite{McDuffSal04}. But when $c_1^{\min}(N)\geq n - 2$, this dimension is $\leq -1$ for all $\beta$ with $\langle \beta,[\omega]\rangle \geq 0$, and so $\EuScript{N}^s(T)^l$ is empty for all $l$.  Since every sphere factors through a simple one, there are then no holomorphic spheres in the fibres of $T$.

We shall need to consider the subset $\EuM(T; J^v)_q^r \subset \EuM(T; J^v)^q$ of sections $u$ whose action $\action(u)$ equals $r$. We shall also need to divide out by the $\R$-action, by translation, on the moduli spaces. For any $J^v\in \EuJ^{\mathrm{reg}}(\Tv T)$, any $\nu_-,\nu_+\in \hor(T,\Omega)$, and any $r\geq 0$, the space $\EuM(T; \nu_-,\nu_+; J^v)_1^r/\R$ is compact and $0$-dimensional. Say it consists of $n(\nu_-,\nu_+)_r $ points, modulo 2. Compactness of $\EuM(T; \nu_-,\nu_+; J^v)_1^r/\R$ follows from Gromov--Floer compactness and the regularity assumptions, which forbid bubbling in the index 1 case. The differential $\partial_{J^v}$ is defined by 
\begin{equation}\partial_{J^v}\langle \nu_- \rangle = 
\sum_{\nu_+} 
{\left (\sum_{r\in \R} {n(\nu_-,\nu_+)_r  t^r } \right)  \langle \nu_+ \rangle}. 
\end{equation}
One has $\partial_{J^v}^2=0$, so that this is a well-formed complex $\CF(T,\Omega;J^v)$. One then uses `continuation maps' to prove the invariance of $\HF(T,\Omega)$ under change of almost complex structure, of parametrisation of $Z$, and under changes $\Omega \leadsto \Omega+ d\alpha$ where $d\alpha$ vanishes on the fibres. Since the monodromy can always be made non-degenerate by such perturbations, the definition can be relaxed to allow degenerate $(T,\Omega)$. See \cite{HS} for proofs (trivial changes must be made to generalise from the case of Hamiltonian symplectomorphisms). 

Finally, when $Z$ is disconnected, one takes the tensor product of the chain complexes attached to the components of $Z$. 

\subsubsection{Lagrangian intersection Floer homology}
Suppose that $L_0$ and $L_1$ are monotone Lagrangian submanifolds in an $(M,\omega_M)$ which is closed or has convex contact-type boundary. Assume that $\min( m^{\min}_{L_0}, m^{\min}_{L_1})\geq 3$. One then has a Lagrangian Floer homology module $\HF_*(L_0,L_1)= H(\CF_*(L_0,L_1; J))$ over $\Lambda_R$ \cite{Oh93}. 

The chain complex $\CF(L_0,L_1;J)  $ is defined on the free $\Lambda_R$-module $\CF(L_0,L_1)$ generated by $L_0 \cap L_1$ (here one first moves $L_1$ by a Hamiltonian diffeomorphism to ensure the intersection is transverse). The definition of the differential $\delta_J$ invokes an element $J=\{J_t\}_{t\in [0,1]}$ of the space $\EuScript{P}\EuJ(M)$ of smooth paths $[0,1]\to \EuJ(M)$ in the space of $\omega_M$-compatible almost complex structures. One forms a moduli space $\EuM(L_0,L_1;J)$ of $L^p_{1,loc}$ maps $u\colon \R\times [0,1]\to M$ with $u(\R\times\{i\})\subset L_i$ ($i=1$, $2$), of finite action $\action(u)= \int_{[0,1] \times\R} {u^*\omega_M} $, satisfying the Cauchy--Riemann equation $J_t \circ D_{(s,t)} u = D_{(s,t)} u\circ i $. The subset $ \EuScript{P}\EuJ^{\mathrm{reg}}(M)\subset \EuScript{P}\EuJ(M)$ is defined by the requirement that the components of the quotient space $\EuM(L_0, L_1;J)/\R$ of formal dimension $\leq 1$ are cut out transversely.

For $J\in  \EuScript{P}\EuJ^{\mathrm{reg}}(M)$, monotonicity and the Maslov index assumption allow one to rule out the bubbles and boundary bubbles which potentially lead to the failure of compactness for moduli spaces of broken trajectories. This is done by an index argument (applied twice: once in proving that $\delta_J\circ \delta_J=0$, and once in defining the continuation maps which prove invariance). Thus $\HF(L_0,L_1)$ is well-defined and independent of $J$ and of Hamiltonian perturbations of $L_0$ and $L_1$.  This is also true when $m^{\min}_{L_0} = 2$ and $L_1=\phi(L_0)$ for some $\phi\in\aut(M,\omega)$ \cite[Addendum]{Oh93} (discs bubble off in cancelling pairs).

Floer homology is well-defined and invariant also when $(M,\omega)$ and the two Lagrangians $L_0$ are $L_1$ are all strongly negative (compare \cite{Laz}). The definition of $ \EuScript{P}\EuJ^{\mathrm{reg}}(M)$ is slightly more elaborate than in the monotone case: the additional requirements are that the moduli spaces of simple $J_t$-spheres, somewhere-injective $J_0$-discs attached to $L_0$, and somewhere-injective $J_1$-discs attached to $L_1$ are all cut out transversely. This ensures that these moduli spaces are negative-dimensional manifolds, hence empty. The possibility of \emph{nowhere}-injective discs is ruled out by a lemma of Kwon--Oh and Lazzarini \cite{KwonOh, Laz} to the effect that the image of any pseudo-holomorphic disc attached to a Lagrangian contains that of a somewhere-injective disc.

\subsection{Open manifolds}
In Theorem \ref{Gysin} we allow $(M,\omega)$ to have a convex contact-type boundary: thus there is a  non-vanishing, outward-pointing vector field $X$ defined near $\partial M$,  satisfying the Liouville condition $L_X\omega=\omega$. For the purposes of pseudo-holomorphic theory on such manifolds, one imposes a convexity condition on almost complex structures (see \cite[Section 2.1]{Seidel03}) which is understood as part of the notion of  a `compatible' almost complex structure. The adaptations described in \cite{Seidel03} do not interact at all with the proof of the theorem in Sections \ref{geom} and \ref{alg}; therefore, to avoid irrelevant complication, in our exposition we shall  only consider closed manifolds. 

\section{Floer homology and Lagrangian matching conditions}\label{Lag match}
\subsection{A cobordism category for Lagrangian correspondences}
Floer homology groups for Lagrangian correspondences fit into a sort of TQFT, or sigma-model, of a non-standard kind, which we now outline. It could be more accurately described as a $(1+1)$-dimensional topological field theory coupled to locally Hamiltonian fibrations and subject to  `Lagrangian matching conditions'. It is under  systematic development  by Wehrheim--Woodward  \cite{WehrWood}, though its existence was observed independently by the author. 
I have decided, after some hesitation, to outline the general theory, since it is an appropriate setting for the symplectic Gysin sequence. The theory as we shall describe it differs from that of \cite{WehrWood} only in some rather minor aspects.

We first describe the objects in a cobordism category, $\EuScript{C}$ (the morphisms will be described later).

\begin{defn}
A \emph{cycle of Lagrangians}, of \emph{length} $l\geq 1$, is a sequence of  symplectic manifolds $\mathsf{M}= (M_1,\dots, M_{l})$ together with a cycle 
\[\mathsf{L} = (L_{1,2}, L_{2,3},\dots, L_{l-1 , l},L_{l,1})\] 
of Lagrangian correspondences\footnote{A Lagrangian correspondence between $M_1$ and $M_2$ is, by definition, a Lagrangian submanifold of $(M_1)_-\times M_2$.}
\[  L_{i, i+1 } \subset  (M_i)_- \times M_{i+1} , \quad i \in \Z/l. \] 
\end{defn}

\begin{defn}\label{object}
An object $\mathsf{O}$ of $\EuScript{C}$ comprises a collection $\{ (\mathsf{M}_\alpha,\mathsf{L}_\alpha)\}_{\alpha \in A}$ of cycles of Lagrangians $(\mathsf{M}_\alpha,\mathsf{L}_\alpha)$, indexed by a finite set $A$, and a collection $\{(M_\beta, \omega_\beta, \mu_\beta)\}_{\beta\in B}$ of compact symplectic manifolds $(M_\beta,\omega_\beta)$ with symplectic automorphisms $\mu_\beta\in \aut(M_\beta,\omega_\beta)$, indexed by another finite set $B$.
\end{defn}

To give a symplectic automorphism $\mu\in\aut(M,\omega)$ is (essentially) equivalent to giving an oriented 1-manifold $Z \cong S^1$ with a marked point $z$, and an LHF $(E,\pi,\Omega)$ over $Z$: given $\mu$ one forms its mapping torus; given $(E,\pi,\Omega)$ and $z\in Z$ one takes the monodromy on $E_z$. 

Cycles of Lagrangians can also be described in fibre-bundle language. Define a \emph{cycle of LHFs} in terms of the following data: an oriented 1-manifold $Z \cong S^1$ with a finite set of marked points enumerated in cyclic order as $(z_1,\dots,z_l)$ dividing $Z$ into intervals $[z_i,z_{i+1}]$; and, for each $i\in \Z/l$, an LHF $(E_i,\pi_i ,\Omega_i)$ over  $[z_{i-1},z_i]$, together with a Lagrangian 
\[\Lambda_{i,i+1} \subset \pi_i^{-1}(z_{i})_- \times \pi_{i+1}^{-1}(z_i).\] 

Write $M_i$ for $\pi_i^{-1}(z_{i-1})$. Parallel transport $m_i\colon M_i \to \pi_i^{-1}(z_{i})$ trivialises $(E_i,\pi_i,\Omega_i)$. One has Lagrangian correspondences $L_{i,i+1} = (\id\times m_i^{-1}) \Lambda_{i,i+1} \subset (M_i)_- \times \pi_i^{-1}(z_i)$. Thus one obtains a cycle of Lagrangians $(M_i; L_{i,i+1})_{i\in \Z/l}$. Conversely, given a cycle of Lagrangians, one can build a cycle of (trivial) LHFs. These constructions are mutually inverse (though only up to isomorphism when one starts with a cycle of LHFs).

\begin{figure}

\centering
\labellist
\pinlabel {$Z\cong S^1$}  at  -40 110 
\pinlabel {$E \cong T_\mu$}  at  95 255 
\pinlabel $z$ at 160 105

\pinlabel $z_1$ at 485 105
\pinlabel $z_3$ at 370 105
\pinlabel $z_2$ at 425 160
\pinlabel $z_4$ at 425 45

\pinlabel $L_{1,2}$ at 525 105
\pinlabel $L_{3,4}$ at 325 105
\pinlabel $L_{2,3}$ at 425 195
\pinlabel $L_{4,1}$ at 425 10

\pinlabel $E_1$ at 550 -10
\pinlabel $E_2$ at 550 220
\pinlabel $E_3$ at 300 220
\pinlabel $E_4$ at 300 -10

\endlabellist
\includegraphics[scale=0.5]{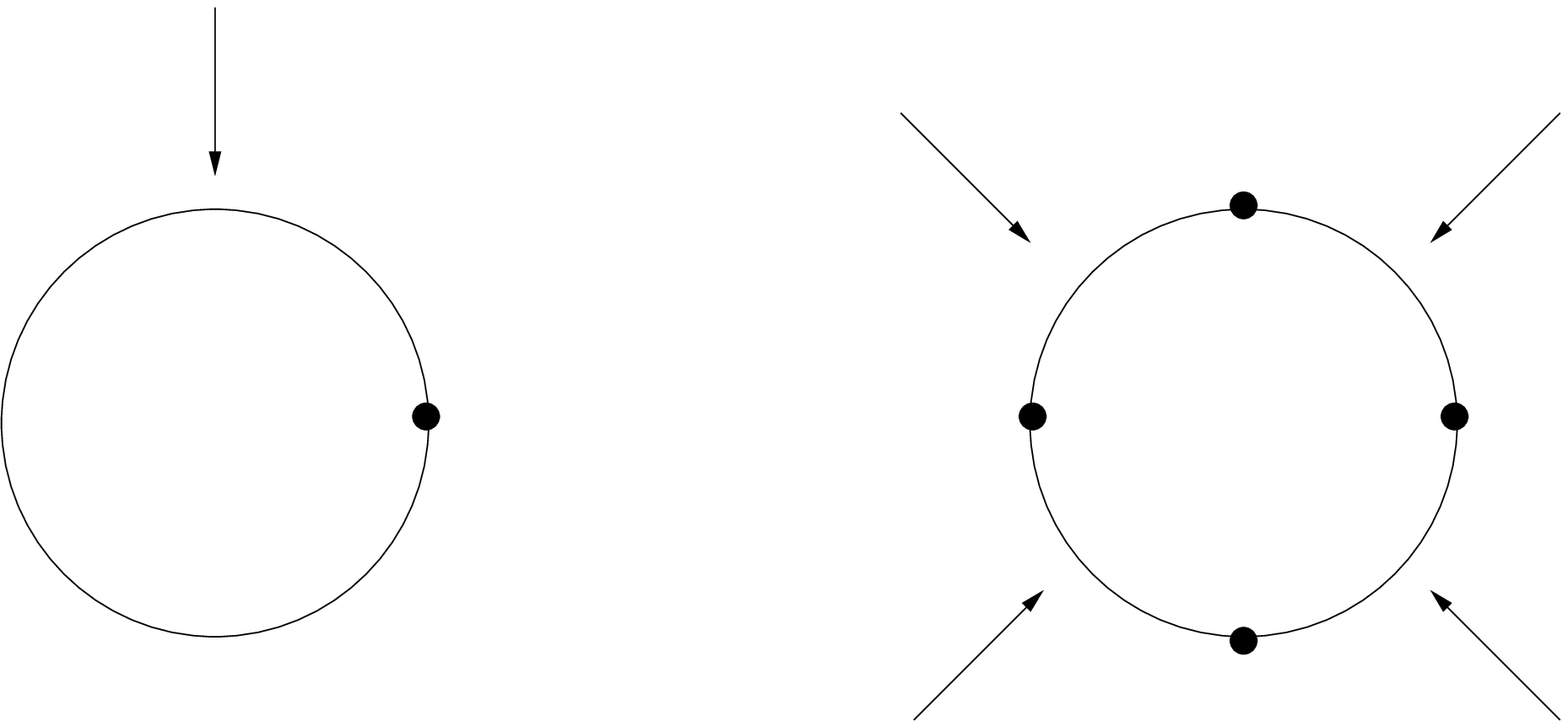}
\caption{Schematic of a mapping torus (left) and of a cycle of Lagrangians of length  4 (right). }\label{cycles}
\end{figure}

For schematic pictures of a mapping torus and a cycle of Lagrangians, see Figure \ref{cycles}.

\subsubsection{Non-degeneracy}
A cycle of Lagrangians $L_{i,i+1}\subset (M_i)_-\times M_{i+1}$ ($i\in \Z/l$) is called \emph{non-degenerate} when the product $\mathbb{L}:=\prod_{i\in \Z/l }{L_{i,i+1} }$ (a Lagrangian in $\prod_{i\in \Z/l}{(M_i)_-\times M_{i+1}}$) meets the multi-diagonal $\Delta_l = \{ (x_1,x_2, x_2, \dots, x_l,x_l, x_1)\}$  (another  Lagrangian) transversely. 

\begin{lem}
$\mathbb{L}=\prod_{i\in \Z/l }{L_{i,i+1}}$ may be moved by a Hamiltonian isotopy to a new Lagrangian, close to the old one in Floer's $C^\infty_\epsilon$-norm, but meeting $\Delta_l$ transversely.   Moreover, the Hamiltonian can be taken to be a product of Hamiltonian isotopies $\pr_2^*(H_{i,t})$ in the factors $(M_i)_- \times M_{i+1}$.
\end{lem}
\begin{pf}
Consider the map 
\[ q\colon  \prod_i{ C^\infty([0,1]\times M_{i+1}) } \times \mathbb{L} \to \mathbb{M}\] 
given by
\[ \left( H_1,\dots, H_l;  (x_1,y_1),\dots, (x_l,y_l) \right) \mapsto  \left((x_1,\phi_{H_1}(y_1)),\dots, (x_l, \phi_{H_l}(y_l)) \right).\]
Observe that $\im(Dq)$ is the subspace $\bigoplus_i {  T_{\phi_{H_i}(y_i)} M_{i+1} }\subset \bigoplus_i {T_{x_i}M_i \oplus T_{\phi_{H_i}(y_i)} M_{i+1}}$. Thus $q$ is transverse to $\Delta_l$.  Form the fibre product
\[Z= \left( \prod_i{ C^\infty([0,1]\times M_{i+1}) } \times \mathbb{L} \right)\times_q \Delta_l,\] 
and note that transversality of $ (\prod {\id_{M_i}\times\phi_{H_i}} ) (\mathbb{L})$  with $\Delta_l$ is equivalent to $(H_1,\dots, H_l)$ being a regular value of the projection $Z\to  \prod_i{ C^\infty([0,1]\times M_{i+1}) }$. The projection is Fredholm, so the Sard--Smale theorem gives the density of the regular values. 
\end{pf}

There is a corresponding notion of non-degeneracy for cycles of LHFs. A \emph{cycle of Hamiltonian chords} $(\sigma_1,\dots,\sigma_l)$ is a sequence of $\Omega_i$-horizontal sections $\sigma_i\colon [z_{i-1},z_i] \to E_i$ where $(\sigma_{i-1}(z_i),\sigma_i(z_i))\in L_{i-1,i}$ for each $i\in \Z/l$. Cycles of Hamiltonian chords are in natural bijection with intersections $\mathbb{L}\cap \Delta_l$ for the corresponding cycle of Lagrangians, and one says that the cycle of LHFs is non-degenerate when the intersection is transverse. It follows from the genericity property just discussed that non-degeneracy can be achieved by making a small adjustment to the forms $\Omega_i$, replacing each by $\Omega_i+d(K_t dt)$ where $t$ is the circle coordinate. 

\subsection{Floer homology from cycles of correspondences}
Before turning to the morphisms in the category $\EuScript{C}$, we explain how to attach homology groups to objects.

\begin{defn}
A compact symplectic manifold $(M,\omega)$ is \emph{normalised monotone} if $\int_{S^2} {f^*\omega} = 2 \langle f^* c_1(TM,\omega) , [S^2]\rangle$ for all smooth maps $f\colon S^2\to M$.  A compact Lagrangian $L\subset M$ is normalised monotone if $\int_{D}{f^*\omega}=m_L(f)$ for all $f\colon (D^2,\partial D^2)\to(M,L)$.
\end{defn}
Notice that the product of normalised monotone symplectic manifolds is normalised monotone. From a cycle of normalised monotone Lagrangian correspondences $L_{i,i+1}\subset (M_i)_- \times M_{i+1}$ between normalised monotone manifolds $\{M_i\}_{i\in \Z/l}$, we form the product
\begin{equation}
  \mathbb{M} = \prod_{i\in \Z/l}  { (M_i)_- \times  M_{i+1} }  = \prod_{i\in \Z/ l}  { M_i \times  (M_i)_-}       
  \end{equation} 
(also normalised monotone) and two normalised monotone Lagrangians inside it: the product $\mathbb{L}= \prod_{i\in \Z/l }{L_{i,i+1}}$, and the multi-diagonal $\Delta_l$, consisting of cycles whose components in $M_i$ and in $(M_i)_-$ coincide for each $i\in \Z/l$.

For each $\mathbb{J}\in \EuJ^{\mathrm{reg}}(\mathbb{M};\mathbb{L},\Delta_l)$ one has a Floer complex 
\begin{equation}  \label{complex}
\CF(\mathbb{L},\Delta_l; \mathbb{J} ), 
\end{equation}
whose chain-homotopy type is independent of $\mathbb{J}$. We wish to show that regularity can be attained when $\mathbb{J}$ is drawn from a smaller pool of almost complex structures.

Recall that $\EuJ(M_i)$ is the space of compatible almost complex structures on $M_i$. Put $\EuJ(M_1, \dots, M_l)=\prod_{i\in \Z/l} { \EuJ(M_i)}$. 
The \emph{pleating map} $\mathrm{pl}$,
\begin{equation}\label{pl} 
\mathrm{pl} \colon \EuScript{P}\EuJ(M_1,\dots, M_l) \to   \EuScript{P}\EuJ(\mathbb{M})\end{equation} 
sends the path $t\mapsto (J_1(t), \dots, J_l(t))$ to the path
\[   t\mapsto \left(J^1(t), - J^1(1-t); \dots ;  J^l(t), -J^l(1-t)\right).  \]
(Here we think of $\mathbb{M}$ as $\prod_{i\in \Z/ l}  { M_i \times  (M_i)_-} $).
Let $ \EuScript{P}\EuJ^{\mathrm{reg}}$ denote those products of $t$-dependent almost complex structures which induce regular almost complex structures for the pair $(\mathbb{L},\Delta_l)$:
\begin{equation}\label{PJ reg}
 \EuScript{P}\EuJ^{\mathrm{reg}}(M_1,\dots, M_l;\mathbb{L},\Delta_l ) = \mathrm{pl}^{-1}( \EuScript{P}\EuJ^{\mathrm{reg}}(\mathbb{M}; \mathbb{L},\Delta_l)).
 \end{equation}

\begin{lem}
$ \EuScript{P}\EuJ^{\mathrm{reg}}(M_1,\dots, M_l;\mathbb{L},\Delta_l )$ is a Baire subset of $ \EuScript{P}\EuJ(M_1,\dots, M_l)$. 
\end{lem}
\begin{pf}
This is a consequence of the Sard--Smale theorem, and only the last stage of the proof  diverges from the standard  arguments found in, e.g., \cite[Lemma 2.4]{Seidel03}. The last stage is the following: let $u\colon [0,1]\times \R \to \mathbb{M}$ be a non-constant, finite-action $(j,\mathbb{J})$-holomorphic map satisfying the boundary conditions $u(\{0 \}\times \R)\subset \mathbb{L}$ and $u(\{1 \}\times \R)\subset \Delta_l$. Let $F\to [0,1]\times \R$ be the bundle dual to $\Lambda^{0,1}(u^* T\mathbb{M})$, and suppose that $\eta$ is a section of $F$, of class $L^{p/(p-1)}$, such that $D^*_{u,\mathbb{J}}\eta =0$ on $ (0,1)\times\R$ and such that
\[   \int_{[0,1]\times \R}{\langle \eta, \mathbb{Y}\circ Du \circ  j \rangle} =0   \] 
for all $\mathbb{Y}\in T_{\mathbb{J}} \EuJ(M_1,\dots, M_l)$. 
Concretely, 
\[ \mathbb{Y}(t)=(Y_1(t),-Y_1(1-t) ; \dots ; Y_l(t),-Y_l(1-t)),\] 
where $Y_i(t) \in \End(TM_i)$ satisfies $J_i(t) \circ Y_i(t) + Y_i (t)\circ j =0$. We must prove that $\eta=0$. 

By elliptic bootstrapping, $\eta$ is actually continuous. The set of points at which $Du =0$ is discrete \cite{FHS}, so if $\eta$ is not everywhere zero, there is a point $(t,s)  \in [0,1]\times \R$ where $\eta(t,s)\neq 0$ and  $D u(x,y) \neq 0$. We may assume $t\neq 1-t$. Let $\alpha=(Du\circ j) (t,s) \in \Hom(\R^2, T_{u(t,s)}\mathbb{M})$; it has components $(\beta_1,\gamma_1; \dots ; \beta_l,\gamma_l)$. Say $(\beta_i,\gamma_i)\neq 0$; then there exist 
$Y_i(t) \in \End(TM_i)$ satisfying $J_i(t) \circ Y_i(t) + Y_i (t)\circ j$ and $Y_i(1-t)\in \End(TM_i)$ satisfying $J_i(1-t) \circ Y_i(1-t) + Y_i (1-t)\circ j$, such that also 
\[ Y_i(t)(\beta_i) \neq Y_i(1-t)(\gamma_i) .\] 
Define $\mathbb{Y}(t)$ by putting $(Y_i(t),-Y_i(1-t))$ in the $i$th slot, and $(Y_j(t),Y_j(1-t))=(0,0)$ for $j\neq i$. Then $\langle \eta(s,t), \mathbb{Y}(t) \circ \alpha  \rangle \neq 0$. Extend $\mathbb{Y}(t)$ arbitrarily to a $\tilde{Y}\in T_J \EuJ(M_1,\dots, M_l)$, then replace $\tilde{\mathbb{Y}}$ by $\mathbb{Y}=\chi \tilde{\mathbb{Y}}$ where $\chi$ is a bump function supported near $(s,t)$. When the support of $\chi$ is small enough, one will then have $ \int_{[0,1]\times \R}{\langle \eta, \mathbb{Y}\circ Du \circ  j \rangle} \neq 0$, contrary to assumption. Hence $\eta \equiv 0$.

Strictly, what this argument proves is that (for a Baire subset of $\mathbb{J}$'s), all non-constant Floer-theoretic trajectories are cut out transversely. In the strongly negative case, there are additional transversality conditions on bubbles. However, these can be handled by almost identical arguments to those for trajectories.
\end{pf}
\begin{rmk}
In forming the complex (\ref{complex}), we shall generally choose $\mathbb{J}$
from the space $ \EuScript{P}\EuJ^{\mathrm{reg}}(M_1,\dots, M_l;\mathbb{L},\Delta_l )$. Thus we assemble $\mathbb{J}$ from $t$-dependent almost complex structures on $l$ (rather than $ 2l$) symplectic manifolds.
\end{rmk} 

\begin{prop}
When $l$ is even, and $\mathbb{J}=\mathrm{pl}(J_1,\dots,J_l)\in \EuJ^{\mathrm{reg}}(M_1,\dots, M_l;\mathbb{L},\Delta_l )$, there is a canonical isomorphism of chain complexes
\[  \CF(\mathbb{L},\Delta_l; \mathbb{J})\cong \CF( L_{1,2}\times L_{3,4}  \times \dots \times L_{l-1,l}, L_{2,3}\times L_{4,5} \times \dots, L_{l,1} ; J_1,\dots,J_l) , \] 
where the latter complex is taken in $M_1\times (M_2)_-\times \dots \times M_{l-1}\times (M_l)_-$. In particular, when $l=2$, we have
\[ \CF(L_{1,2}\times L_{2,1}, \Delta_2 ) \cong \CF(L_{1,2}, L_{2,1}).  \]
\end{prop}

\begin{pf}
A generator for $ \CF(\mathbb{L},\Delta_l; \mathbb{J})$ is an intersection point $\left((x_1,x_2), (x_2,x_3),\dots, (x_l,x_1)\right)$ with $(x_i,x_{i+1})\in L_{i,i+1}$; we map this to the generator  $(x_1,x_2,\dots, x_l)$ for the complex on the right and so obtain a linear isomorphism between the two complexes. Checking that this is a chain map is left to the reader.
\end{pf}

An object $\mathsf{O}$ of $\EuScript{C}$ is, according to Definition \ref{object},  a finite collection (indexed by some set $A$) each of whose elements is a sequence of symplectic manifolds $(M_1^\alpha,\dots, M_{l^\alpha}^\alpha)$ and a cycle of Lagrangian correspondences $L_{i,i+1}^\alpha$, together with a finite collection (indexed by $B$) of symplectic manifolds $N_\beta$ equipped with automorphisms $\mu_\beta$. 

Given a monotone (or strongly negative) $N$ and a non-degenerate automorphism $\mu\in \aut(N)$, we have a family of chain complexes $\CF(\mu; J^v)$ indexed by $J^v \in \EuJ^{\mathrm{reg}}(\Tv \torus(\mu))$. 

We associate with $\mathsf{O}$ the homology
\begin{equation} 
\HF(\mathsf{O}) = H\left( \bigotimes_{\alpha\in A}{\CF_*(\mathbb{L}_\alpha,\Delta_{l_\alpha})} \otimes \bigotimes_{\beta\in B}{\CF(\mu_\beta) } \right ),  
\end{equation} 
which is independent, up to canonical and coherent isomorphisms, 
of the additional choices of $\mathbb{J}_\alpha\in \EuJ^{\mathrm{reg}}(\mathbb{L}_\alpha,\Delta_{l_\alpha})$ and $J^v_\beta \in \EuJ(\Tv \torus(\mu_\beta)$).

\subsubsection{Morphisms}
A morphism in $\EuScript{C}$ is an isomorphism class of `cobordisms', where a cobordism is a \emph{matched collection of fibrations}. The latter is defined to be a formidable assembly of data 
\begin{equation}\label{matched collection} \EuScript{F}=(S,\Gamma; E,\pi,\Omega; Q; \{\zeta_\alpha\}_{\alpha \in \pi_0(\partial S)}),\end{equation} 
consisting of the following components:
\begin{enumerate}
\item[(i)] A quilted surface $(S,\Gamma)$;
\item[(ii)] an LHF $(E, \pi ,\Omega)$ over the surface $ S_\Gamma$ obtained by cutting along $\Gamma$;
\item[(iii)] a Lagrangian matching condition $Q \to \Gamma$ for $E$;
\item[(iv)] distinguished points $\zeta^\alpha$ in the boundary.
\end{enumerate}
Here are the definitions.

\emph{Component (i):}
\begin{defn}
A \emph{quilted surface} (Wehrheim--Woodward's term) $(S,\Gamma)$ consists of a compact oriented surface $S$, possibly with boundary, and a smoothly embedded 1-manifold with boundary $\Gamma$, intersecting $\partial S$ transversely, with $\partial \Gamma \subset \partial S$.
\end{defn}
We denote by  $S_\Gamma$ the compact surface with boundary and corners obtained by cutting $S$ open along $\Gamma$, i.e. by adjoining boundary arcs and circles to
compactify $S\setminus \Gamma$.  The boundary of $S_\Gamma$
consists of two copies of $\Gamma$ which we denote $\Gamma^+$ and $\Gamma^-$ (they may be taken in either order), and a union of closed arcs and circles  $T$ coming from $\partial S$ and meeting $\Gamma^+\cup \Gamma^-$ at the corners. Thus
\[ \partial S_\Gamma = \Gamma^+ \cup \Gamma^- \cup T. \]
\begin{figure}
\begin{center}
\labellist
\pinlabel {$S$}  at  250 410 
\pinlabel {$\Gamma$}  at  480 370 
\pinlabel {$\Gamma^+$}  at  530 185 
\pinlabel {$\Gamma^-$}  at  715 300 
\pinlabel {$T$}  at  780 270 
\pinlabel {$T$}  at  590 270 
\pinlabel {$T$}  at  190 290 
\endlabellist
\includegraphics[scale=0.4]{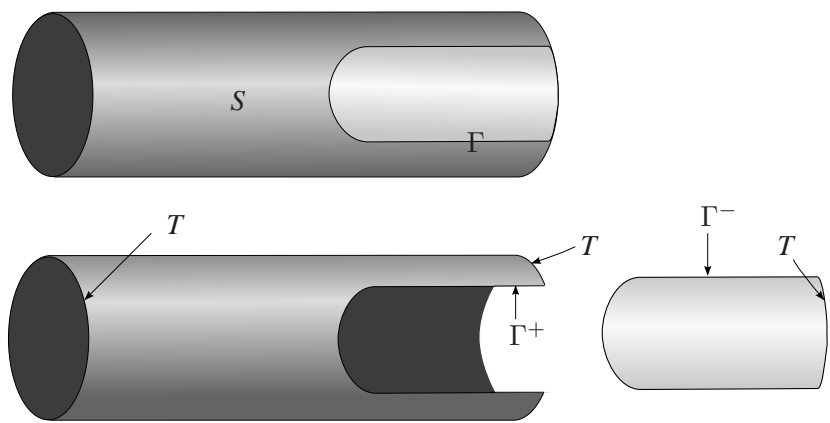}
\caption{Above: a quilted cylinder $(S,\Gamma)$. Below: the corresponding surface $S_\Gamma$ with the three parts of its boundary labelled.}\label{quilt}
\end{center}

\end{figure}
\emph{Component (ii):} this is an LHF $(\pi\colon E\to S_\Gamma,\Omega)$. It is allowed to have differing dimensions over different components of $S_\Gamma$.  However, \emph{we insist that the fibration is flat over a neighbourhood of $T$.} We shall be particularly interested in its restrictions $E|_{\Gamma^\pm}\to \Gamma^\pm$, and in the closed two-forms $\Omega^\pm$ on  $E|_{\Gamma^\pm}$ obtained by restricting $\Omega$.

It is important to note that $\Gamma^+ $ is naturally identified with $\Gamma$ (and that the same goes for $\Gamma^-$). Thus \emph{we may regard $E|_{\Gamma^\pm}$ as LHFs over $\Gamma$.} Form the
fibre product
\begin{equation}\label{E Gamma}
  E_\Gamma : =   (E|_{\Gamma^+}) \times_{\Gamma} (E|_{\Gamma^-}) .
  \end{equation}
It becomes an LHF over $\Gamma$ when we endow it with the two-form $\Omega_\Gamma$ obtained as the restriction of the form $(-\Omega^+) \oplus \Omega^-$ on $  (E|_{\Gamma^+}) \times (E|_{\Gamma^-}) $ to its subspace $E_\Gamma$. The minus sign is crucial (though to which factor it is attached is not). 

\emph{Component (iii):}
\begin{defn}
A \emph{Lagrangian matching condition} $Q \to \Gamma$ for $E$ is a sub-fibre bundle $Q\subset E_\Gamma$,  proper over $\Gamma$, such that (i) $\Omega_\Gamma |_Q=0$, and (ii) the fibres of $Q$ have half the dimension of those of $E_\Gamma$.  
\end{defn}
We can equally regard $Q$ as a sub-fibre-bundle of $(E|_{\Gamma^-}) \times_{\Gamma} (E|_{\Gamma^+})$; the ordering of $\Gamma^+$ and $\Gamma^-$  is immaterial.

\emph{Component (iv):}  a collection of reference points $\{ \zeta^\alpha\}_{\alpha \in \pi_0(\partial S)}$, one in each component $(\partial S)^\alpha$ of $\partial S$. We require that $\zeta^\alpha \in \partial \Gamma$ when $\partial \Gamma \cap (\partial S)^\alpha \neq \emptyset$. 

\emph{This concludes the definition of matched collections of fibrations.}

Now consider what all this gives us over the original boundary $\partial S$. For any component $Z=(\partial S)^\alpha$ of $\partial S$, the matched collection of fibrations gives by restriction a cycle of LHFs. The marked points are those of the finite set $\partial \Gamma \cap Z$, enumerated in cyclic order, consistently with the boundary orientation, as $(z_1^\alpha,\dots, z_{l_\alpha}^\alpha)$, starting with $z_1^\alpha=\zeta^\alpha$. When $\partial \Gamma \cap Z=\emptyset$, we just take the one-element sequence given by $z^\alpha$. As described above, a cycle of LHFs can be trivialised so as to give a cycle of Lagrangians. 

\begin{defn}
A \emph{cobordism} from $(\mathsf{M}_-,\mathsf{L}_-)$ to $(\mathsf{M}_+,\mathsf{L}_+)$ is defined to be a matched collection of fibrations, together with a symplectic identification of the cycles of Lagrangians obtained from its boundary with $\overline{\mathsf{M}}_-\amalg \mathsf{M}_+$. Here $\overline{\mathsf{M}}_-$ is the object obtained from $\mathsf{M}_-$ by reversing the direction of each of its constituent sequences; thus $(M_1,\dots, M_l)$ becomes $(M_l,\dots M_1)$, and the Lagrangian correspondences are obtained in the obvious way (when $l=0$, $\mu$ changes to $\mu^{-1}$).  A morphism in the category $\EuScript{C}$ is defined as an isomorphism class of cobordisms. 
\end{defn}
The definition of `isomorphism' here is the obvious one, involving diffeomorphisms of the base surfaces and two-form-preserving bundle-maps; it is left to the reader to elaborate. 

Morphisms are composed by concatenation of cobordisms. The picture to have in mind is that one can \emph{cut} the quilted surface $(S,\Gamma)$ (and the fibrations lying above it) along a collection of circles transverse to $\Gamma$; sewing is just the inverse operation. For the sewing of fibrations, it is convenient to be able to assume that the LHFs are flat, hence (by Lemma \ref{flat}) locally trivial near the boundaries being glued (compare \cite[Section 1]{Seidel03}). The formalities are again left to the reader.  

\subsection{Cobordism-maps}
We now turn to the construction of maps on Floer homology from matched collections of fibrations.
\subsubsection{Cylindrical ends; almost complex structures}
Let $\EuScript{F}$ be a matched collection of fibrations. We elongate the base surface $S$ to a surface $\wh{S}$ with cylindrical ends. To do this, choose for each boundary component $(\partial S)_\alpha$  a boundary collar $e_\alpha \colon (-1,0]\times S^1 \hookrightarrow S$ such that $e_\alpha^{-1}(\Gamma)$ is a union of line segments $(-1,0]\times \{t\} $. Then define
\begin{equation}\label{elongation}
 \wh{S}=  S \cup \bigcup_\alpha {(-1,\infty)\times S^1}  
 \end{equation}
where the $\alpha$th copy of $(-1,\infty)\times S^1$ is glued on along the region $(-1, 0]\times S^1$ via the collar $e_\alpha$.  Extend $\Gamma$ linearly to a 1-manifold $\wh{\Gamma}\subset \wh{S}$: that is, if $e_\alpha^* \Gamma$ contains a segment $(-1,0] \times \{t\}$, one extends it to $(-1,\infty)\times \{t\}$. Let $\wh{S}_\Gamma$ be the surface with boundary obtained by cutting $\wh{S}$ along $\wh{\Gamma}$. 

We may choose $e$ so that the LHF $(E\to S_\Gamma,\Omega)$ is flat near $T$ (recall the flatness assumption in the definition). It then extends uniquely to an LHF $(\wh{E}\to \wh{S}_\Gamma,\wh{\Omega})$ which is flat over the ends. The Lagrangian matching condition $Q$ extends to a Lagrangian matching condition $\wh{Q}$ over $\wh{\Gamma}$ in a unique way; indeed, the requirement that it be isotropic implies that  the extension must be obtained by symplectic parallel transport along the ends of $\wh{\Gamma}$.

\begin{defn}\label{a c conditions}
(a) Given an LHF $(E,\pi,\Omega)$ over a Riemann surface $(S,j)$, an almost complex structure $J$ on $E$ is \emph{partly compatible} if $J\circ D\pi = D\pi \circ j$, and the restriction of $J$ to each fibre $E_x$ is compatible with the symplectic form $\Omega|_{E_x}$.

(b) A partly compatible $J$ is \emph{fully compatible} if it preserves the horizontal distribution.

(c) Suppose $S$ has a cylindrical end, holomorphically modelled on $(-\infty, -T] \times S^1$ or $[T,\infty)\times S^1$. We say $J$ is \emph{adapted} if it is partly compatible everywhere, fully compatible over the end, and translation-invariant (after shortening the end by increasing $T$). 

(d)  When the cylindrical ends of $\wh{S}$ have seams $\wh{\Gamma} =(-\infty, -T) \times \{ z_1,\dots, z_l \}\subset (-\infty, -T) \times  S^1$, we choose $j$ to be \emph{standard} so that the end is modelled holomorphically on $(-\infty, -T) \times  S^1$ with seams $(-\infty, -T) \times \{ e^{2\pi i/l},e^{4 \pi i/l}\dots , e^{2\pi il/l} \}$, i.e., the $z_k$ are equally spaced. Given an LHF over the resulting Riemann surface $\wh{S}_\Gamma$, we consider almost complex structures $J$ which are adapted in the same sense as (c). 
\end{defn}
\subsubsection{Holomorphic sections}
The ends of the surface $\wh{S}$ are cylindrical, with seams $\wh{\Gamma}$. Fix a complex structure $j$ on $\wh{S}$ (inducing $\hat{j}$ on $\wh{S}_\Gamma$) which is standard over the ends, as in Definition \ref{a c conditions}, and consider the space $\EuJ(\wh{E}_\Gamma)=\EuJ(\wh{E}_\Gamma,j)$ of $j$-adapted almost complex structures in $\wh{E}_\Gamma$. Over each end $e$, $J$ is eventually translation invariant and preserves the horizontal distribution, and is there determined by its restriction $J^v(e)$ to the vertical tangent bundle. We write  $\EuJ(\wh{E}_\Gamma; \{ J^v_e\}) $ for the subspace of $\EuJ(\wh{E}_\Gamma)$ of almost complex structures limiting to the vertical complex structures  $J^v_e$ over the ends. The $J^v_e$ are the complex structures used in the definition of the Floer complexes, and so must be chosen \emph{regular} (cf. Definition \ref{regular}).

Given a section $u$ of $\wh{E} \to \wh{S}_\Gamma$ one obtains, by restriction to the boundary, a section $u_\Gamma$ of $\wh{E}_\Gamma$ (the completed version of $E_\Gamma$). At a point $x\in \Gamma$, the value of $u_\Gamma(x)$ is the pair $(u(x_+),u(x_-))$, where $x_+$ and $x_-$ are the points of $\Gamma^+$ and $\Gamma^-$ that correspond to $x$ under the identifications $\Gamma^+\cong \Gamma\cong \Gamma^-$.

For a fixed $p>2$, the configuration space $\EuScript{B}^p$ associated with our non-degenerate matched collection of fibrations is the disjoint union of subspaces $\EuScript{B}^p(\nu)$ indexed by cycles of Hamiltonian chords $\nu$. The space $\EuScript{B}^p(\nu)$ consists of the sections $u$ of  $\wh{E} \to \wh{S}_\Gamma$ of class $L^p_{1,\mathrm{loc}}$ with the properties (i) that  $u_\Gamma$ is a section of $Q\subset \wh{E}_\Gamma$, and (ii) the function $e_\alpha^* u(s,t) - \nu_\alpha(t)$ is an $L_1^p$-function on the cylinder for each cylindrical end $e_\alpha$. The $L^p_1$ topology makes this space into a Banach manifold. The \emph{action} of $u\in \EuScript{B}^p(\nu)$ is 
\begin{equation}
\action(u)= \int_{\wh{S}_\Gamma}{u^*\Omega},
\end{equation} 
when the integral exists.  Given an adapted almost complex structure $J$, the moduli space of $J$-holomorphic sections is the subspace
\begin{equation} 
 \EuScript{Z}(\nu;J) := \{ u\in  \EuScript{B}(\nu):   J\circ Du  - Du\circ j  = 0 ,\; \action(u)<\infty \} \subset \EuScript{B}^p(\nu) . 
 \end{equation}
Thus  $\EuScript{Z}(\nu;J) = \dbar_J^{-1}(0)$ where $\dbar_J$ is the section $\dbar_J = J \circ D - D\circ j$ of the natural Banach vector bundle $\EuScript{V}^p\to \EuScript{B}^p (\nu)$ whose fibre at $u$ is the subspace 
\begin{equation} (\EuScript{V}^p)_u \subset L^p(\Lambda^{0,1}_{\wh{S}_\Gamma}(u^*T\wh{E})) \end{equation}
consisting of those $\C$-antilinear homomorphisms $\alpha\colon T\wh{S}_\Gamma \to u^*T\wh{E}$ whose restrictions to $\Gamma$ satisfy the linearised matching  condition 
\begin{equation}
\left(\alpha_{x_+}(\xi),\alpha_{x_-}(\xi)  \right)\in T_{\left((u(x_+),u(x_-)\right)} Q \;  \text { for all }  \; \xi\in T_x\wh{\Gamma} .
\end{equation}
The linearisation of $\dbar_J$ at $u$ is a map 
\begin{equation}\label{DuJ}
D_{u,J}=D_u\dbar_J\colon \EuScript{T}_u^p\to \EuScript{V}_u^p, 
\end{equation} 
where $\EuScript{T}^p_u =T_u\EuScript{B}^p(\nu)$ is the space of globally $L^p_1$ sections of $u^* \Tv \wh{E}$ tangent to $Q$ at the boundary.
\begin{lem}
When $\nu$ is non-degenerate, the map $D_{u,J}$ from (\ref{DuJ}) is Fredholm, with index independent of $p$. 
\end{lem}
\begin{pf}[Sketch proof]
The Fredholm property in the quilted setting turns out to be no more difficult than that of the construction of Floer homology \cite{Floer} and its continuation maps.
We first discuss the situation when $p=2$. 

We can find a finite cover $\{U_\alpha\}_\alpha$ for $\wh{S}_\Gamma$ of the following type: each $U_\alpha$ is the image in $\wh{S}_\Gamma$ of a corresponding subset $U_\alpha'\subset \wh{S}$. Moreover,  $(U_\alpha', U_\alpha' \cap\wh{\Gamma}')$ is the image of an embedding of pairs whose domain is one of the following: 
\begin{enumerate}
\item[(i)]
$\left(D^2,\emptyset \right)$; 
\item[(ii)] 
$\left(S^1\times [0,\infty),\emptyset \right)$ or $\left(S^1\times (-\infty,0], \emptyset\right)$; 
\item[(iii)]
$\left(D^2, (-1,1) \times \{0\}\right)$; 
\item[(iv)]
$\left(S^1 \times (-\infty,0], \{ e^{2\pi i k/m}: k= 1,\dots, m\} \times  (-\infty,0]\right)$  or \\
$\left( S^1\times [0, \infty), \{ e^{2\pi i k/m}: k= 1,\dots, m\} \times [0,\infty) \right).$
\end{enumerate}
Let $i_\alpha \colon U_\alpha \to \wh{S}_\Gamma$ be the inclusion. 
In case (iii), $i_\alpha^*\xi$ is a pair $(\xi_+,\xi_-)$ of sections of bundles $E_\pm$ over $D^2_\pm = \{ z \in D^2: \pm \im(z)\geq 0\}$. But $ \xi_+ +\xi_-\circ \tau$ then defines a section of $E_+\oplus\tau^*E_-\to D^2_+$, where $\tau\colon D^2_+\to D^2_-$ is complex conjugation. This section is subject to a Lagrangian boundary condition, the linearisation of $Q$. Standard elliptic theory then shows that for any precompact open set $U\subset \wh{S}_\Gamma$ and any compact $K\subset U$, one has 
\begin{equation}\label{elliptic}  
\| i_\alpha^*\xi\|_{L^2_k(K)} \leq \mathrm{const.}\big(\| i_\alpha^*\xi \|_{L^2(K)} + \| i_\alpha^*D_{u,J} \xi\|_{L^2_{k-1}(U)}\big).  
\end{equation}
Now consider first an end of type (iv). One sees that in this case $i_\alpha^*\xi $ is a section of a bundle over the strip $\{ e^{2\pi i t}: t\in [0,1/m] \} \times (-\infty,0] $, subject to non-degenerate linear Lagrangian boundary conditions. (The reduction is constructed by `pleating', as in our treatment of almost complex structures.)   We may take the bundle to be a trivial complex bundle, of rank $r$, say, and the boundary conditions to be two transverse totally real subspaces $(\Lambda_0,\Lambda_1)$. The operator $D_{u,J}$ reduces to one of shape $\frac{d}{ds} + L$, where $L=i \frac{d}{dt}\colon L^2_1([0,1/m], \C^r; \Lambda_0,\Lambda_1)\to L^2([0,1/m],\C^r)$ is formally self-adjoint and invertible. The operator $\frac{d}{ds} + L\colon L^2_1(\R\times [0,1/m], \C^r;\Lambda_0,\Lambda_1) \to L^2(\R\times [0,1/m], \C^r)$ is an isomorphism \cite[Section 3.1]{Don}. The story for ends of type (ii) is similar. For both types, the $L^2$-norm of $\xi\in \ker D_{u,J}$, over the infinite strip, is controlled by its $L^2$-norm over a compact band. It follows that the unit ball in $\ker D_{u,J}\cap L^2$ is compact, hence that the $L^2$ kernel is finite-dimensional.  

As to finite-dimensionality of the cokernel, the point is to construct a parametrix for $D_{u,J}$, i.e. a bounded linear operator  $P\colon\EuScript{V}_u^2\to \EuScript{T}^2_u$ such that $D_{u,J}P $ differs from the identity map by a compact operator. The construction is done by patching parametrices $P_\alpha$ for sections supported in $U_\alpha$. For the patching procedure (essentially that of \cite[p. 54]{APS}) see \cite[Section 3.2]{Don}.   The existence of a parametrix in cases (i) and (iii) follows from elliptic theory. In cases (ii) and (iv) we take the inverse $P_\alpha$ to the invertible operator $d/ds+ L$.

The $L^p$ theory introduces some additional complications, but the hard work is in dealing with the strip itself \cite{Floer, Don}. Elements of $\ker D_{u,J }\cap L^p$ decay exponentially over the ends, and so lie in any $L^q$. Similarly for the kernel of the formal adjoint operator. 
\end{pf}

We define an enlargement $\widetilde{\EuScript{Z}}(\nu,J)$ of $\EuScript{Z}(\nu,J)$ as the space pairs $(u, [\gamma])$, where (i) $u\in \widetilde{\EuScript{Z}}(\nu,J)$, (ii) $\gamma$ is a finite disjoint union of trees of simple parametrised $J$-holomorphic spheres in fibres of $\pi$, rooted at points of $\im(u)$. Thus $\widetilde{\EuScript{Z}}(\nu,J)$ contains  $\EuScript{Z}(\nu,J)$ as the subspace where $\gamma$ is the empty union. The topology is that of Gromov convergence. 

The pleating trick allows us to establish Gromov--Floer compactness, which says that says that any bounded-action sequence in $\widetilde{\EuScript{Z}}(\nu)$ converges, in the Gromov--Floer topology, to a \emph{broken section} plus finitely many trees of spherical bubbles in fibres plus finitely many trees of discs in boundary fibres. The point here is that the process of bubble formation is local (in the domain) and therefore the usual proofs need only cosmetic changes.

\subsubsection{Definition of the cobordism-maps}
Let $\mathsf{O}_-$ and $\mathsf{O}_+$ be objects in $\EuScript{C}$, and $\EuScript{F}$ a matched collection of fibrations defining a cobordism between them. Let $J_-^v$ and $J^v_+$ be regular complex structures defining the Floer complexes $\CF_*(\mathsf{O}_\pm ;J^v_\pm)$. Also pick $J \in \EuJ(\wh{E}_\Gamma; \{ J^v_\pm\})$.  Given cycles of Hamiltonian chords $\nu_-$ and $\nu_+$ for the respective objects, there is a locally constant action function $\mathsf{A}$ on  $\widetilde{\EuScript{Z}}(\nu_-,\nu_+;J)$ with sub-level sets $\widetilde{\EuScript{Z}}( \nu_-,\nu_+;J)_{\leq c}:= \mathsf{A}^{-1}((-\infty,c])$.  Assuming we are in the normalised monotone (or strongly negative) case, the zero-dimensional part $ \widetilde{\EuScript{Z}}(\nu_-,\nu_+;J )_{\leq c}^0$ is compact, because its only possible ends arise from breaking trajectories, and that can only occur when the index is $\geq 1$.  Hence the zero-dimensional part of the level set of $c$, $\widetilde{\EuScript{Z}}(\nu_-,\nu_+;J)^{0}_{c}=\mathsf{A}^{-1}(c)\cap   \widetilde{\EuScript{Z}}(\nu_-,\nu_+;J )^0$, is also compact. Put
\begin{equation}\label{counting}
  n_J(\nu_-|\EuScript{F} |\nu_+) =   \sum_{r\in \R} {\#\widetilde{\EuScript{Z}}(\nu_-,\nu_+;J)_{r}^0 t^r} \in \Lambda _R.    
  \end{equation}
These coefficients define a chain map $\CF(\EuScript{F})\colon \CF(\mathsf{O}_- ;J_-)\to \CF(\mathsf{O}_+;J_+)$: 
\begin{equation}\label{cob map}
\CF(\EuScript{F})(\nu_-) = \sum_{\nu_+}{n(\nu_-|\EuScript{F}|\nu_+) \langle \nu_+ \rangle}.\end{equation} 
These maps compose compatibly with sewing of fibrations. This is a consequence of the standard gluing principles for Cauchy--Riemann operators. We shall only invoke sewing along cylindrical ends, where the relevant gluing theorem is (in all analytical essentials) the one that underpins Hamiltonian Floer homology (cf. \cite{Salamon94}, for instance).

\subsubsection{Quantum cap product}
Assuming $N$ monotone, the quantum cap product
\[ H^p(N;\Lambda_R) \otimes_{\Lambda_R} \HF_*(\mu)\to \HF_{*-p}(\mu), \quad c\otimes x \mapsto c \frown x,   \]
could in principle be defined by means of cocycles for any ordinary cohomology theory in which one can make sense of transversality. We opt for smooth singular cohomology. 

A smooth singular $q$-chain $\sigma$ in $N$ is a formal linear combination $\sum_k{a_k \sigma_k}$, where $a_k\in \Lambda_R$ and $\sigma_k \colon\Delta^q \to N$ is a smooth $q$-simplex for each $k$. `Smooth' means that $\sigma_k$ extends to a $C^\infty$ map  defined on an open neighbourhood of $\Delta^q$ in $\R^q$, so in particular, the restriction of $\sigma_k$ to a face of $\Delta^q$ is also smooth. These chains form a complex $\mathsf{S}_*^{\mathrm{sm}}(N)$; its dual is the smooth singular cochain complex $\mathsf{S}^*_{\mathrm{sm}}(N)$. 

Take a partly $\omega_\mu$-compatible almost complex structure $J$ on $\R\times \torus(\mu)$, eventually translation invariant and fully compatible over the two ends of the cylinder. For $\nu_\pm \in \hor(\torus(\mu))$, one has the `continuation map' moduli spaces $\widetilde{\EuScript{Z}}(\nu_-,\nu_+;J)$, and evaluation maps 
\[ \ev \colon \widetilde{\EuScript{Z}}(\nu_-,\nu_+;J) \to N,\quad u\mapsto u(0,[0]).\]

Given any countable collection of smooth simplices $\sigma^k\colon \Delta^{d(k)}\to N$, we can find perturbations $J'$ of $J$, supported over a neighbourhood of $(0,[0])\in \R\times S^1$, so that the $q$-dimensional moduli spaces $ \widetilde{\EuScript{Z}}(\nu_-,\nu_+;J')^q$ are transverse to the $\sigma^k$ for all $k$ and all $q\leq q_0$ (see \cite[Lemma 2.5]{Seidel03}, for instance). One then has fibre products
\begin{equation}\label{fibre products}
  \widetilde{\EuScript{Z}}( \nu_-,\nu_+;J' )\times_{\mathrm{ev} } \sigma^k. 
  \end{equation} 

Choose a codimension $p$ smooth singular chain $\tau$ such that $\tau$ and $\partial \tau$ are transverse to the $ \widetilde{\EuScript{Z}}(\nu_-,\nu_+;J')^q$. By counting the points of fixed energies in the fibre products (\ref{fibre products}) one defines the matrix entries in a linear map
\begin{equation}  \tau \frown \cdot \colon \CF_*(\mu)\to \CF_{*-p}(\mu).  
\end{equation} 
The usual trajectory-breaking argument demonstrates that 
\begin{equation} \partial (\tau\frown x)= (\partial \tau)\frown x + \tau \frown \partial x.\end{equation} 
This argument requires transversality for the 1-dimensional fibre products, and the fact that these 1-dimensional spaces compactify to (topological) 1-manifolds with boundary. The latter fact is part of the gluing theory for holomorphic spheres \cite{McDuffSal}. To define the quantum cap product action of a smooth singular \emph{co}chain $c$ is defined to be that of its Poincar\'e dual $c\cap [N]$.
\begin{rks}
(i) In this paper we will only need the quantum cap-product action of a few specific chains. We have not defined a chain-level action of the \emph{whole} of $\mathsf{S}^{\mathrm{sm}}_*(N)$.

(ii) There was nothing special about the base $\R\times S^1$ used in this construction; there is a similar story for any matched collection of fibrations with a marked point in the base (not on a seam), and it is compatible with sewing of these fibrations. We will not say any more about the formalities.

(iii) The quantum cap product also exists in the strongly negative case. The construction is easy, because if the almost complex structures are chosen generically one does not encounter bubbles at all (bubbles appear  in generic 3-parameter families of almost complex structures). 
\end{rks}\

\section{Proof of the Gysin sequence: geometric aspects}\label{geom}
\subsection{Setting up the isomorphism}
The maps $\rho^!$ and $e(V)\frown \cdot $ in the Gysin sequence (\ref{Gysin}) fit into the field theory described in the previous section.
Let $\wh{V}'= (\id_{M_-}\times \mu) \wh{V} \subset M_- \times N$. As promised in the introduction, we will construct an isomorphism
\[ \Phi  \colon H \cone(e(V) \frown \cdot) \to
\HF_*(\wh{V}, \wh{V}')  \]
as the map on homology induced by a quasi-isomorphism
\begin{equation}\label{H Cg}
(H,C\rho) \colon  \cone(e(V) \frown \cdot) \to
\CF_*(\wh{V}, \wh{V}'). \end{equation}
So
\[ C\rho\colon \CF_*(\mu)[-k-1] \to \CF_*( \wh{V},  \wh{V}') \]
is a chain map, and $H$ a nullhomotopy of $C\rho \circ Ce$. Here $Ce$ is the map defined by cap product with a cochain representing $e(V)$.

\subsubsection{Global angular cochains}\label{global angular}
The Euler class of an oriented sphere bundle vanishes when pulled back to the total space:
\begin{equation} \rho^*e(V)=0. 
\end{equation}
When $V$ is an $\SO(k+1)$-bundle, this can be seen very simply: $\rho^*V$ has a tautological section. A more precise statement, valid for any smooth, $R$-oriented fibre bundle whose fibres are $R$-homology $k$-spheres, is that that the orientation class $o$ on the fibre $F$ \emph{transgresses} to the Euler class  $e(V)$, where we work with coefficients in the ring $R$. Recall that the transgression $\tau$ is the homomorphism
\[ d_{k+1}^{0,k}\colon E^{0,k}_{k+1}\to E^{k+1,0}_{k+1}  \]
in the Leray--Serre spectral sequence $(E^{p,q}_r, d_r^{p,q} )$ of a fibre bundle $\pi\colon E\to B$. One thinks of it as a map from a subgroup of $H^k(F;R)$ to a quotient of $H^{k+1}(B;R)$. The cochain interpretation is that $\tau( [x] )=[y] $ when $x$ extends to a cochain $w\in C^k(E;R)$ such that $dw=\pi^*y$.

In the case of an $R$-oriented homology-$S^k$-bundle $\rho\colon V\to N$, there is a transgression homomorphism $\tau \colon H^k(F;R)\to H^{k+1}(N;R)$. One has $\tau(o)= e(V)$ (indeed, $e(V)$ is sometimes defined this way). Hence there exists $\psi\in C^k(V;R)$ such that  $d\psi=\rho^* e(V)$ and $[\psi|_{S^k}]=o$. We shall call such a cochain $\psi$ a `global angular cochain'.\footnote{Bott and Tu call the `de Rham' analogue of $\psi$ a \emph{global angular form} \cite[p. 121]{BottTu}.} Dually, if $Z$ is a (codimension $k+1$) chain Poincar\'e dual to $e(V)$, one can find a `global angular chain': a codimension $k$ chain $Y$ on $V$ such that $\partial Y = \rho^! Z$ and $[Y\cap \rho^{-1}(x)]= [ \mathrm{pt.} ]$.

References for this discussion are \cite{Serre, BottTu}.  A small technical concern is that Serre \cite{Serre} sets up the spectral sequence via cubical singular (co)chains. However, from this one can readily deduce its analogue for the smooth singular (co)chains complexes used here.

\subsubsection{Defining $C\rho$}
\begin{figure}
\centering
\labellist
\pinlabel {$C$} at 200 390
\pinlabel {$B$} at 450 390
\pinlabel {$\Gamma$} at  450 440
\pinlabel {$p$} at 325 390
\pinlabel {$S$} at 300 490
\pinlabel {$z_1$} at 565 350
\pinlabel {$z_2$} at 565 440
\endlabellist
\includegraphics[scale=0.6]{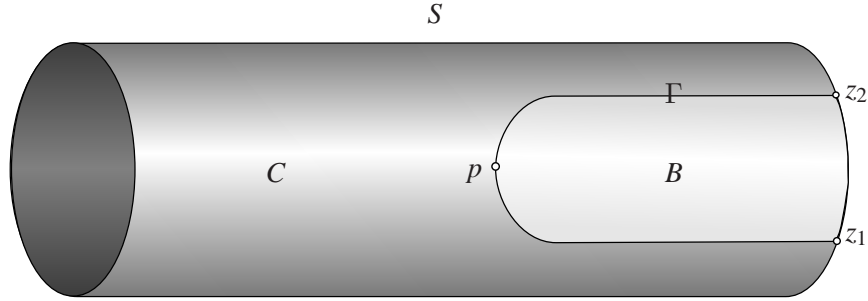}
\caption{The base of the matched pair of fibrations underlying $C\rho$.}\label{gysincyl}
\end{figure}

The map $C\rho$ is associated with a matched collection of fibrations over a quilted surface depicted  schematically in Figure \ref{gysincyl}. 
The base surface $S$ is a finite cylinder; its elongation $\wh{S}$ (see (\ref{elongation})) an infinite cylinder. To be precise, we set $\wh{S}=\C/ 4i\Z
$, and
\begin{equation} 
S=\{[z]\in\wh{S}: |\real(z)| \leq 2 \}, \end{equation}
then define $q\colon \C \to \wh{S}$ as the quotient map. The arc $\wh{\Gamma}\subset S$ is the union of two line segments and a circular arc\footnote{Pedantically, we ought to smooth the joins between the line segments and the arc, but this detail will be suppressed.}:
\begin{equation}
\wh{\Gamma} = q(\{ z: \imag(z)=\pm 1,\real(z)\geq 1 \}) \cup q(\{z: |z-1|=1, \real(z) \leq 1\}) .\end{equation} 
Figure \ref{gysincyl} shows the finite part $\Gamma:=\wh{\Gamma}\cap S$. The point $p$ shown in the diagram (which is nothing more than a convenient reference point) is $q(0)$.  The closures of the two components of $\wh{S}\setminus \wh{\Gamma}$ are denoted by $B$ (the convex part) and $C$, as in the figure.

To construct an LHF over $C$, let $R=(N \times \C)/\Z$ where $n\in \Z$ acts by $(x,z)\mapsto (\mu^n(x) ,z+4in)$. There is a projection map $\pi_R \colon R \to \wh{S}$, $[x,z]\mapsto [z]$, and a natural closed two-form $\Omega_R$ on $R$ whose pullback to $N\times \C$ is $\pr_1^*\omega$. Now let $R=\pi_R^{-1}(C)$. Let 
\begin{equation}
\pi_Q\colon Q\to C
\end{equation} be the restriction of $\pi_R$ to $Q$, and $\Omega_Q$ the restriction of $\Omega_R$ to $Q$. Then $(Q,\pi_Q,\Omega_Q)$ is a flat LHF over $C$. Let 
\begin{equation} 
P = B\times M,
\end{equation} 
and let $\pi_P\colon P\to B$ be the trivial projection. Give $P$ the closed 2-form 
\begin{equation}\label{Omega P}
\Omega_P : = \omega_M + d  (K dt)
\end{equation} 
for some function $K\in C^\infty(M)$. Here $z=s+it$ is the complex coordinate on $\C$ (or on $\wh{S}$). For later reference, we record the following easy observation.
\begin{lem} \label{P is flat} 
$(P,\pi_P,\Omega_P)$ is a flat LHF over $B$.
\end{lem}

The fibre $Q_p= \pi_Q^{-1}(p)$ is identified with $N$ by thinking of it as the image of $N\times \{0\} \subset N\times \C$ under the $\Z$-action. The fibre $P_p=\pi_P^{-1}(p)$ is identified with $M$ because $\pi_P$ is by definition a trivial bundle. Our Lagrangian submanifold $\wh{V}\subset M_-\times N$ appears here:
\begin{equation}
\wh{V}_p :=  \wh{V} \subset (P_p)_- \times Q_p = M_- \times N. 
 \end{equation}
We denote this submanifold by $\wh{V}_p$ to emphasise that it lies in the product of the fibres over $p \in \wh{S}$.  Symplectic parallel transport along the arc $\wh{\Gamma}$, starting at $p$, carries $\wh{V}_p$ to a Lagrangian $\wh{V}_x\subset (P_x)_-\times Q_x$ for each $x\in \wh{\Gamma}$. Write
\begin{equation} 
\wh{\mathsf{V}} =\bigcup_{x\in \wh{\Gamma}} {\wh{V}_x}.
\end{equation}
\begin{defn}
The map $C\rho$ is the cobordism-map (\ref{cob map}) associated with the matched pair of fibrations $(P,Q, \wh{\mathsf{V}})$. Thus it counts finite-action pairs $(u,v)$, where $u$ is a section of $\pi_P$, $v$ a section of $\pi_Q$, both pseudo-holomorphic for adapted almost complex structures (\ref{a c conditions}), and $(u(x),v(x))\in \wh{V}_{x}$ for each $x\in \wh{\Gamma}$.
\end{defn}
This definition needs fleshing out. The incoming boundary of $S$ (which is the one on the left) is a circle without marked points, and over it we have an LHF isomorphic to the mapping torus $\torus(\mu)\to S^1$. The outgoing boundary (on the right) is a circle with two marked points, $z_1=[2-i]$ and $z_2=[2+i]$. Now, let us identify the fibre of $P$ over $z_1$ with $M$ \emph{via symplectic parallel transport along $\Gamma$, starting at $p$}. Similarly, identify the fibre of $Q$ over $z_2$ by with $N$ via symplectic parallel transport along $\Gamma$ from $p$. The restriction of $P$ to the directed boundary-arc $[z_1,z_2]$ is trivialised by parallel transport starting at $z_1$. The restriction of $Q$ to the directed boundary-arc $[z_2,z_1]$ is trivialised by parallel transport starting at $z_2$. By using these trivialisations we see that the target for $C\rho$ is $\CF_*( (\id \times \mu^{-1})( \wh{V} ),\wh{V} )=
\CF_*( \wh{V} ,(\id \times \mu)( \wh{V} ))$. Thus this matched pair of fibrations defines a map
\begin{equation} \label{chain g}
C\rho\colon \CF_*(\mu)\to \CF_*(\wh{V} ,(\id \times \mu)( \wh{V} ) ).
\end{equation}

\subsubsection{Defining $Ce$}
Choose a smooth singular chain $Z$ representing the Poincar\'e dual to 
the Euler class $e(V)$. Define the quantum cap product map 
\[ Ce =  Z \frown \cdot \colon \CF_*(\mu) \to  \CF_{*-k-1}( \mu). \] 
Let 
\[ e \colon \HF_*(\mu) \to  \HF_{*-k-1}( \mu) \]
be the induced map on homology.
\subsubsection{The composite} 
\begin{figure}[t]
\centering
\labellist
\pinlabel {$p$} at 325 390
\pinlabel {$p$} at 180 80
\pinlabel {$p'$} at -145 85
\endlabellist
\includegraphics[scale=0.3]{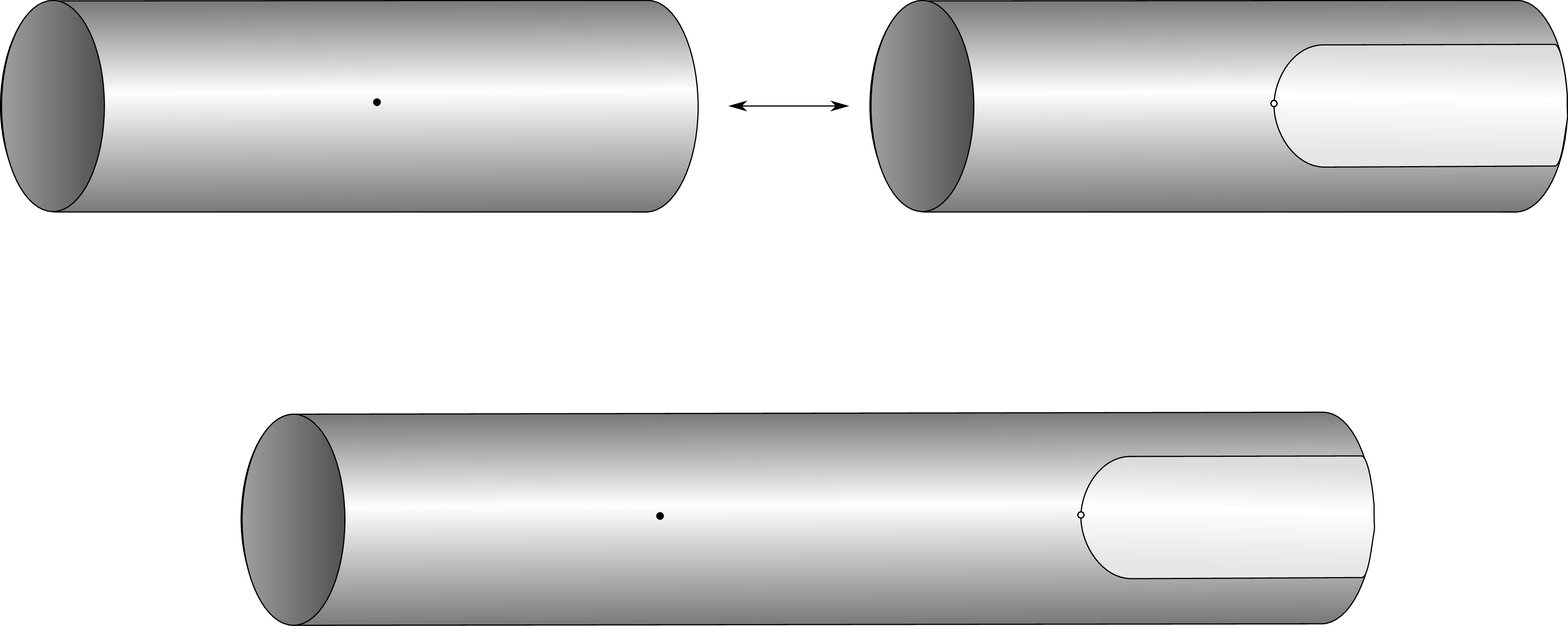}
\caption{}\label{gysingluing}
\end{figure}

The gluing (or sewing) theorem for fixed-point Floer homology expresses the composite $C\rho \circ Ce$ in terms of a sewed fibration, as in Figure \ref{gysingluing}. Namely, $Ce$ is defined using the LHF $(R,\pi_R,\Omega_R) $ over $ \widehat{S}$ defined above; $C\rho$ using a matched pair of LHFs over $\widehat{S}$. We can glue these together, as in the figure, to obtain a new matched pair of LHFs over $\widehat{S}$; in fact, this new pair are naturally identified with $P$ and $Q$. The only new feature is that we have an additional marked point in the interior of $C$. For definiteness, let us say that this marked point is $p':=q(-1)$. In the diagram, we have indicated the marked points where incidence conditions are specified with solid black dots; marked points which are needed just for reference are shown as hollow dots.

The composite  $C\rho\circ Ce$ counts pairs of pseudo-holomorphic sections  $(u,v)$ as above, but now $v$ sends the marked point $p'$ to a point of $Z\subset M = \pi_Q^{-1}(p')$.

\subsubsection{The nullhomotopy}
\begin{figure}[t]
\centering
\labellist
\pinlabel {$p'$} at -135 65
\pinlabel {$p$} at 180 60
\endlabellist
\includegraphics[scale=0.4]{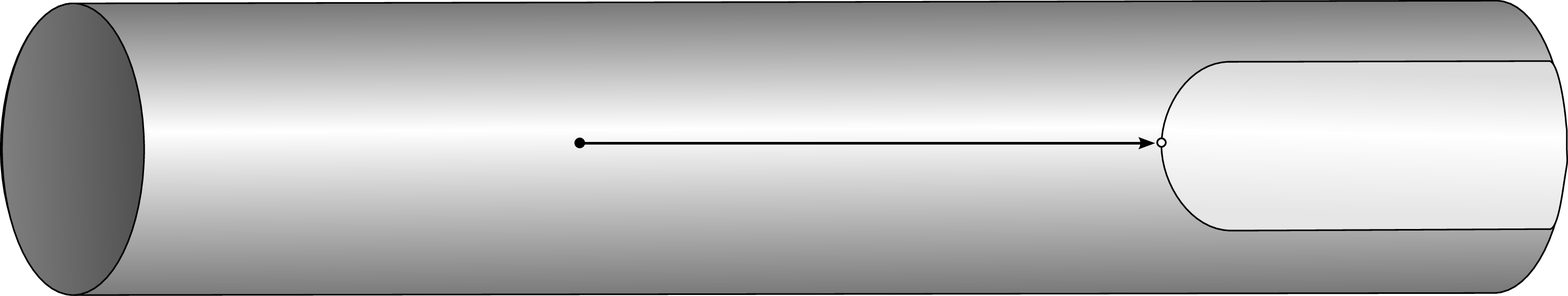}
\caption{}\label{gysinhomotopy}
\end{figure}

The nullhomotopy $H$ promised in (\ref{H Cg}) will be constructed as the concatenation of two homotopies, $H_1$ and $H_2$.  The first homotopy, $H_1$, arises by moving the marked point from $p'$ to $p$, as in Figure \ref{gysinhomotopy}.

The matrix entries of this homotopy are defined as counts of isolated triples $\{ (t,u,v)\} $, where $(u,v)$ is a pair of pseudo-holomorphic sections as before, $t \in[-1,0]$, 
and $u(q(t) ) \in \phi_{t}(Z)$, where $\phi_t \colon Q_{p}\to Q_t$ is (leftward) symplectic parallel transport along the path $q([t,0])$. Formally,
\begin{equation}  H_1 \langle \nu \rangle =  \sum_{x \in \wh{V}\cap \wh{V}' } {\sum_{r\in \R} {\#  \left ( \widetilde{\EuScript{Z}}( \nu ,x; J  )^k_r \times [ -1,0]) \times_{\ev, \phi}  (Z\times [-1,0] )\right) t^{r}\langle x \rangle}}
\end{equation}
where $\ev(z,t)=u(q(t))$ and $\phi(\cdot ,t)= \phi_t \circ Z$. This definition requires $J$ (which amounts to a pair of complex structures $(J_P,J_Q)$ on the two fibrations $P$ and $Q$) to be regular in the sense that the $\ev$ is transverse to $\phi$. 

We have to take a little care over compactness: whilst spherical bubbles in the fibres are tolerable, disc-bubbles attached to $\wh{V}$ are not. This is where the assumption {\sc mas} comes into play. When $(t;u,v)$ contributes to the moduli space, the pair $(u,v)$ has index $k+1$; consequently, when $\wh{V}$ is monotone, disc-bubbles in the compactified moduli space have Maslov index $\leq k+1$, and consequently do not occur when $m_{\wh{V}}^{\min}\geq k+2$.  
In the strongly negative case, we can choose our almost complex structures so as to have the generic property that there are no disc-bubbles at all (note that this can be done using product complex structures on $M_-\times N$, because these suffice for a transversality of evaluation arguments along the lines of \cite[Lemma 2.10]{Seidel03}).

\begin{cor}
We have $\partial H_1+H_1 \partial =  C\rho \circ C e + f$, where $f$ is the homomorphism defined via pairs $(u,v)$ with $v(1)\in Z$. 
\end{cor}
\begin{pf}
This is almost an instance of a standard principle, but its non-standard aspect it that we initially mark the \emph{interior} point $p'\in C$, but then slide this to the \emph{boundary} point $p$. We have to show that, despite this, the (formally) 1-dimensional component $\EuScript{Z}_1$ of the moduli space of triples $(t; u,v)$, where $v(t)\in \phi_t(Z)$, is a 1-manifold with boundary. This comes down to showing that the evaluation map $\ev \colon [0,1]\times \EuScript{Z}_1 \to N $, $(t,u,v)\mapsto v(t)$, can be made transverse to the cycle $Z$.

First consider an index $k$ moduli space $\EuScript{Z}_k(x,y)$  of pairs $(u,v)$ without point constraints.  We have an evaluation map $ev_1 \colon \EuScript{Z}_k(x,y) \to \wh{V}_1$. We can choose the complex structures on $P$ and $Q$ so as to make $ev_1$ transverse to the codimension $k$ cycle $\rho^! Z$ in $\wh{V}$, by an argument explained by Seidel \cite[Lemma 2.5]{Seidel03}, for example. 

It follows that $\pr_N \circ ev_1$ is transverse to $Z$. This projected evaluation map is the same as the evaluation map $ev^v_1$, where we evaluate only the second component of the pair $(u,v)$. Now, $ev^v_t \colon \EuScript{M}_k(x,y) \to N $ depends smoothly on $t$. Since transversality is an open property, it follows that for small $\epsilon>0$ and $t\in [1-2\epsilon,1]$,  $ev^v_t$ is also transverse to $Z$. Combining this with an unproblematic transversality argument for $[0,1-\epsilon]$ now gives the result.
\end{pf}

The second homotopy, $H_2$, arises via a global angular cochain. As in (\ref{global angular}), let $Y = \psi \cap [V] \in \mathsf{S}^{\mathrm{sm}}_{2n-2k}(V)$ be a global angular chain. Thus $\partial Y = \rho^! Z$ and $[Y\cap \rho^{-1}(x)] = [ \mathrm{pt}. ] \in H_0(V;R)$. We use $Y$ to define a nullhomotopy of $f$. Define 
\begin{equation} 
H_2 \colon \CF(\mu) \to
	\CF(\wh{V},(\id_{M} \times \mu) \wh{V}), \quad
	c \mapsto  \sum_{x}{n_1(\widetilde{c} | \mathrm{gr}(\rho|_Y) )| x) \, \langle x \rangle}.  
\end{equation}
The coefficient $n_1(\widetilde{c}| \mathrm{gr}(\rho|_Y)  |x)$ counts the isolated points in a one-parameter moduli space, weighted by $t^\epsilon$ with $\epsilon$ the action.  We have
\begin{equation} \partial H_2 + H_2 \partial =  f.  \end{equation}
Indeed, the right-hand side corresponds to pairs of sections $(u,v)$ where $u(p)$ lie in the image of $\rho^! Z$ (that is, in $\rho^{-1}( \im Z)$); but this is equivalent to $v(p)$ passing through $Z$.
Hence
\begin{equation}
 \partial (H_1+H_2) + (H_1+H_2) \partial = C\rho\circ Ce, 
 \end{equation} 
i.e.  $H:=H_1+ H_2$ is a nullhomotopy of $C\rho \circ Ce$.

\subsection{Lagrangian intersections}
The chain complexes defined in the previous section are well-defined up to homotopy-equivalence, and the maps between them are well-defined up to homotopy. The resulting homology modules are independent of choices up to  isomorphism (in fact, canonical isomorphism) and the maps between them are canonical. However, the proof of Theorem \ref{Gysin} appeals to specific choices, to be elaborated in this subsection and the next. 

Our first task is to analyse the intersections of $\wh{V}$ with a particular Hamiltonian perturbation of $(\id_{M_-}\times \mu)(\wh{V})$. To begin, $\wh{V} \cap (\id_{M_-} \times \mu)\wh{V} $ is the set of
pairs $(x,\rho x)$ such that $x\in V $ and $\rho x \in \fix(\mu)$:
\begin{equation}  \pr_1 \left( \wh{V} \cap (\id_{M_-} \times \mu)\wh{V}\right)
= \rho^{-1}( \fix(\mu ) ). 
\end{equation}
After a $C^k$-small Hamiltonian perturbation, $\mu$ has non-degenerate fixed points, according to the following well-known lemma which follows from \cite[Theorem 3.1]{HS}.
\begin{lem}
There is a dense set of functions $h\in C^\infty(N)$ such that $\phi_h\circ \mu$ has non-degenerate fixed points, where $\phi_h$ is the time-one Hamiltonian flow. 
\end{lem}
We choose such an $h$ and rename $\phi_h \circ \mu$ as $\mu$ (this does not affect Floer homology). The intersection
$\wh{V} \cap (\id_{M_-} \times \mu)\wh{V} $ in $M_- \times N$ is then clean and diffeomorphic to $\fix(\phi_h \circ \mu)\times F$, where $F \simeq S^k$ is the typical fibre of $\rho$. 

To make the intersection transverse, we shall make a controlled Hamiltonian perturbation to $\id_{M_-}$. By the coisotropic neighbourhood theorem, the symplectic structure in a neighbourhood of $V\subset M$ is determined by that of $N$. For a specific model, choose a connection $\alpha$ in the fibre bundle $\rho\colon V\to N$. It determines an embedding
$i_\alpha\colon  T^*_{v}(V)= (\Tv V)^* \to T^*V$ as a sub-bundle complementary to $\Ann(\Tv V)$. On the space
\begin{equation}\label{def W}
W:= i_\alpha(T^*_v V) \subset T^* V, 
\end{equation}
we consider the symplectic form
\begin{equation}\label{sigma}
   \sigma = i_\alpha^* d\lambda_V + R^* \rho^*  \omega_N ,  
\end{equation}
where $R\colon W\to V$ is the projection to $V$ and $\lambda_V$ the canonical one-form on $T^*V$. Fix a Riemannian metric on $V$, and let $W_\epsilon \subset W$ be the 
subset of vertical cotangent vectors of length $<\epsilon$.

For small $\epsilon$, there is a symplectomorphism $\mathrm{nd}(V) \to W^\epsilon$ extending the zero-section $V\to W_\epsilon$; we may therefore consider $W_\epsilon$ as a neighbourhood of $V$. 

The space $W_\epsilon$, with its projection $p= \rho\circ R \colon W_\epsilon \to N$, is an LHF. Its Hamiltonian curvature (\ref{curvature}) is determined by the curvature of $\alpha$. It is convenient to choose the connection form $\alpha$ to be flat over a neighbourhood
$U= \bigcup_{\bar{x}\in \Fix(\mu)}{U_{\bar{x}}}$ of the finite set $\Fix(\mu)$.
The Hamiltonian connection is then flat over $U$, hence (after shrinking $U$) symplectically trivial.
For each $\bar{x}\in \Fix(\mu)$, choose a symplectically trivialising projection
\begin{equation} \label{triv}
t_{\bar{x}} \colon p^{-1} U_{\bar{x}} \to (T^* F)_{< \epsilon},  
\end{equation}
so that the symplectic form on $p^{-1}U_{\bar{x}}$ is $p^* \omega_N + t_{\bar{x}}^*d\lambda_{F}$.

For each $\bar{x} \in \fix(\mu)$, fix a bump function $\eta_{\bar{x}}\in C^\infty(N)$,
supported in $U_{\bar{x}}$. Choose a Morse function $m$ on $F$, and let $\hat{m} \in C^\infty(T^*F)$ be the function
$\hat{m} =m \circ s$, where $s \colon T^*F \to F $ is the projection.
The Hamiltonian flow $\phi_{\hat{m}}$ with respect to the canonical form $-d\lambda_{can}$ is `vertical', and moves the zero
section to the graph of $dm$.

Now consider functions $H\in C^\infty(M)$ of shape $\chi \cdot H_m$, where $\chi$ is a cutoff function for
the neighbourhood $W_\epsilon$ of $V$, and
\begin{equation}\label{first Hamiltonian} 	
H_m =  \sum_{\bar{x} \in \fix(\mu) } ( R^* \rho^* \eta_{\bar{x}} ) \cdot (\hat{m} \circ t_{\bar{x}}) .
\end{equation}

\begin{lem}\label{int points}
After rescaling $m$, all intersection points \[ (x,\rho x)\in \wh{V} \cap (\phi_{\chi H_m} \times \mu)\wh{V}\] satisfy
$(x,\rho(x))\in \Fix(\phi_{\chi H_m}\times \mu)$.
\end{lem}

\begin{pf}
Since $\bar{x}$ is a non-degenerate fixed point of $\mu$, there are Euclidean
coordinates $v$ on $U_{\bar{x}}$ centred on $\bar{x}$, and a constant $C>0$, such that
\[ \| \mu(v)-v \| \geq C \|v \| \]
when $\|v \|$ is sufficiently small. We may assume that this inequality is valid in $\supp(\mathrm{\eta_{\bar{x}}})$.
On the other hand, by rescaling $m$ by small positive reals, $\phi_{\chi H_m}$ can be made arbitrarily $C^0$-close to the
identity. In particular, we may arrange that $| \rho (y) - \rho\circ \phi_{\chi H_m}(y) | \leq C/4$ for all $y\in W_\epsilon$.

Take a point $(x,\rho (x)) = (\phi_{H_m} (y), \mu \rho (y)) \in \wh{V} \cap (\phi_{\chi H_m} \times \mu)\wh{V}$. Then $\rho(x)\in \supp(\eta_{\bar{x}})$ for some $\bar{x}\in\fix(\mu)$, since if not, $H_m\equiv 0$ near $x$, so $x=y$ and $\rho(x)\in \fix(\mu)$, which is absurd. We have
\[ 0 = \| \mu\rho (y ) - \rho(x) \| =
\| \mu\rho (y ) - \rho \phi_{\chi H_m} (y)  \| \geq \frac{1}{2} C \| \rho(y) - \rho \phi_{H_m}(y) \|,
\]
so $\rho(x) = \rho \phi_{\chi H_m}(y) = \rho (y)$. Thus $\rho(x) \in \Fix(\mu)$. Now, $\phi_{\chi H_m}$ preserves the fibre
$\rho^{-1}(\rho(x))$, and (identifying the fibre with $(T^*F)_{< \epsilon}$ via $t_{\rho(\bar{x})}$) it acts as $\phi_{\hat{m}}$.
So for $\phi_{\chi H_m} (y)$ to be a point of $V$, we must have $x \in \crit(m)$ and $x = \phi_{\chi H_m}(y)= y$. Thus
$(x,\rho(x))\in \Fix(\phi_{\chi  H_m}\times \mu)$, as required.
\end{pf}
The intersection points in the previous lemma are transverse because $m$ was Morse.

The functions $H_m$ are not Morse, and it will be useful to observe that result continues to hold for a convenient class of Morse functions on $M$. Start with a Morse function $l \in C^\infty(N)$ such that each
$\bar{x}\in \Fix(\mu)$ occurs as a local maximum of $l$. On $W_\epsilon$, put
\begin{equation}\label{second Hamiltonian}
 K'_{l,m} = H_m +  \delta \sum_{c \in \crit(l)} ( R^* \rho^* l ) \cdot (\hat{m} \circ t_c)  \end{equation}
where $t_c$ are local trivialising projections (like the projections $t_{\bar{x}}$ from (\ref{triv})), and $\delta \ll 1$. 
Put
\begin{equation}\label{third Hamiltonian}
K_{l,m} = |\cdot|^2 + K'_{l,m}
\end{equation}
where $|\cdot|^2$ denotes the norm-squared of cotangent vectors to $F$.
The functions $K_{l,m}$ are Morse; the critical points lie on $V$, and project to the critical points of $l$ on $N$.

\begin{lem}
Lemma \ref{int points} remains true when $H_m$ is replaced by $K'_{l,m}$ or $K_{l,m}$.
\end{lem}
\begin{pf}
The proof is similar but a little more complicated.

Given $\epsilon>0$, we can find a constant $C'$ so that if $\mathrm{dist}(y,\mu(y))<C'$ then $y$ lies within $\epsilon$ of a fixed point of $\mu$ (an easy consequence of the compactness of $N$). Hence, when $\delta$ is small enough, the following holds: take any point $(x,\rho (x)) = (\phi_{K_{l,m}'} (y), \mu \rho (y)) \in \wh{V} \cap (\phi_{K_{l,m}'} \times \mu)\wh{V}$; then
 $\rho(x)$ lies in $\bigcup_{\bar{x}\in \fix(\mu)} {\supp(\eta_{\bar{x}}) }$.

Making $\delta$ still smaller, we may suppose that $| \rho (y) - \rho\circ \phi_{K_{l,m}'}(y) |\leq C/4$ for all $y\in W_\epsilon$. 
It is still true that 
$\phi_{K_{l,m}}-|\cdot|^2$ preserves the fibre $\rho^{-1}(\rho(x))$ when $\rho(x)\in \fix(\mu)$. By the same argument as before, one sees that for any intersection point $(x,\rho(x))\in \wh{V}\cap  (\phi_{K_{l,m}' } \times \mu)\wh{V}$ one has $\rho(x)\in \fix(\mu)$. 

Moreover, the functions
$K_{l,m}'$ and $|\cdot|^2$ Poisson-commute, which implies that
\[ \phi_{K_{l,m}} = \phi_{ | \cdot |^2}\circ
\phi_{K_{l,m}' }. \]
But $\phi_{|\cdot|^2}$ preserves $V$, and hence for any intersection point $(x,\rho(x))\in \wh{V}\cap  (\phi_{K_{l,m}} \times \mu)\wh{V}$ one has $\rho(x)\in \fix(\mu)$ (because this is true for $K_{l,m}')$. 

Finally, take an intersection point $(x,\rho(x))\in \wh{V}\cap  (\phi_{K_{l,m}' } \times \mu)\wh{V}$ (or in $\wh{V}\cap  (\phi_{K_{l,m}} \times \mu)\wh{V}$) with 
$\rho(x)\in \fix(\mu)$. Then $\rho(x)=\rho(y)$ and $x\in \crit (m)$ as before.
\end{pf}

\begin{prop}\label{bijection}
When $m$ is small in $C^0$ and has precisely two critical points, $\rho$ induces a two-to-one map $\beta\colon \wh{V} \cap (\phi_{K_{l,m}} \times \mu)\wh{V} \to \fix(\mu)$. 
\end{prop}
\begin{pf}
Immediate from the previous lemma.
\end{pf}

\begin{rmk}
$\beta$ gives a bijection between the generators of the complexes
$\cone(e)$ and $\CF_*(\wh{V}, (\phi_{K_{l,m}}  \times \mu)\wh{V})$.
\end{rmk}
\subsubsection{The definition of $K$}\label{def K}
Recall our standing hypothesis that $F$ is has a Morse function with precisely two critical points. Thus, by an observation of Reeb (see \cite{Milnor}), $F$ is homeomorphic to $S^k$.\footnote{Conversely, by results due to Smale---the h-cobordism theorem and the non-existence of exotic 5-spheres \cite{Smale}---any $F$ homeomorphic to $S^k$ with $k\neq 4$ has such a Morse function.}
We are in a position to stipulate that, from now on, 
\begin{center}
\framebox{
\parbox[c]{0.75\textwidth}{
$K$ should be one of the functions $\chi K_{l,m}$,
where $K_{l,m}$ is from (\ref{third Hamiltonian}),
$\chi$ is a cutoff function for the neighbourhood $W_\epsilon$ of $V$,
and $m\in C^\infty(F)$ has precisely two critical points.
}}\end{center}
\subsection{Almost complex structures}
\subsubsection{On $M$}
We continue to work on  $W \cong T^*_v V$ (\ref{def W}). Notice that the connection $\alpha$ induces isomorphisms  
\begin{align}\label{splitting}
T^*W & \cong R^*(T^*V) \oplus R^* (T^*_v V)  \\
\notag & \cong R^*(T^*_v V)\oplus \left( R^* (T_v^* V) \oplus R^* \rho^*(T^*N) \right).
\end{align} 
We shall consider almost complex structures $I$ which have a block decomposition 
\begin{equation} \label{block matrix}
I  =  \left [  \begin{array}{ccc}
0 & - \id & 0 \\   
\id & 0 & 0 \\
0 & 0 & R^*\rho^*J_N
 \end{array}\right]    
 \end{equation}
with respect to the splitting (\ref{splitting}), for some almost complex structure $J_N$ on $N$ compatible with $\omega_N$. One has $\sigma(Iu,Iv)=\sigma(u,v)$ and $\sigma(v,Iv)>0$ for $v\neq 0$, so there is an induced metric $g_I$ making the splitting orthogonal. 
\begin{lem}\label{tangent}
The gradient $\nabla f$ of a function on $W$, with respect to $g_I$, is tangent to the zero-section $V\subset W$ at $x\in V$  if and only if $df (x) $ annihilates the first summand  $R^*(T^*_v V)_x$. In the case of  the functions $H_m$, $K_{l,m}'$ and $K_{l,m}$, this holds for all $x\in V$.
\end{lem}
\begin{pf}
The first assertion is easily checked. That $dH_m$ annihilates the `vertical' summand $R^*(T^*_v V)$ follows from the Leibnitz rule  and the facts that both $d(\hat{m}\circ t_{\bar{x}})$ and $R^*\rho^* d\eta_{\bar{x}}$ kill the relevant summand. The same argument applies to $dK_{l,m}'$. Since $d( | \cdot |^2)$ vanishes along the zero-section $V$, $dK_{l,m}$ too annihilates  $R^*(T^*_v V)$ along $V$.
\end{pf}
\begin{Not}
Write $\EuJ(M; V)$ for the space of complex structures $I\in  \EuJ(M)$, compatible with $\omega_M$, whose restriction to $\mathrm{nd}(V)$ have the block decomposition (\ref{block matrix}) for some identification of $\mathrm{nd}(V)$ with $W_\epsilon$ and some connection 
$\alpha$.
\end{Not}
\subsubsection{On the matched pair of fibrations}
Recall that  $\pi_P\colon P\to B$ is a trivial $M$-bundle, whilst  $\pi_Q\colon Q\to C$ is an $N$-bundle. We endowed $Q$ with the closed 2-form $\Omega_Q$ which makes it a flat LHF with monodromy $\mu$ around the puncture. We gave $P=B\times M$ the closed 2-form $ \Omega_P = \omega_M + d( K dt) $,
where $t$ is the vertical (or imaginary) coordinate on $B$ and $K\in C^\infty(M)$; in Lemma \ref{P is flat} we observed that this LHF is also flat. Remember also that we built the Lagrangian matching condition $\wh{\mathsf{V}}$ using the symplectic parallel transport defined by $\Omega_P$ and $\Omega_Q$.

We now stipulate conditions on our almost complex structures.
\begin{itemize}
\item
$C$ has an incoming cylindrical end $e\colon (-\infty, 0]\times \R/\Z \to C$, $(s,t)\mapsto (s-2 , 4t)$. The pullback $e^*Q\to  (-\infty, 0]\times \R/\Z$ is also cylindrical: it is $\pr_2^* T_\mu $. We choose an adapted, regular almost complex structure $J_{Q,\infty}$ on $e^*Q$, so that the chain complex $\CF_*(\mu; J)$ is defined. 
\item
The quilted surface $\wh{S}$ also has an outgoing cylindrical end $f \colon [0,\infty) \times \R/\Z$, embedded by $f(s,t)=(s+2, 4t)$, with seams $[0,\infty)\times \{1/4\}$ and $[0,\infty)\times \{ -1/4\}$. The matched pair of fibrations, pulled back by $f$, is also translation-invariant. 

We choose a regular pair of adapted almost complex structures   $J_{P,+\infty}$ and $J_{Q,+\infty}$ on the two LHFs. We take $J_{P,+\infty}$ to be of a special form, namely, $J_{P,+\infty}|_{N\times \{(s,t)\} }= \phi_t^* J_M$ for an almost complex structure $J_M\in \EuJ(M;V)$.
Here $\phi_t$ is the time-$t$ Hamiltonian flow corresponding to the function $K$. 
\item
We draw $J_Q$ from the space of fully compatible almost complex stuctures (\ref{a c conditions}) on $Q$. For $J_P$ on $P=B\times M$ we are still more prescriptive. Choose a compatible almost complex structure $J_M$ on $M$. \emph{Define $J_P=J_P(J_M,K)$ to be the unique fully compatible almost complex structure on the bundle $B\times M\to B$ which restricts to the fibre $M\times \{(s,t)\}$ as $\phi_t^*J_M$.}
\item
Let $\EuScript{J}(M,Q)$ denote the set of pairs $(J_M, J_Q)$ where $J_M\in \EuJ(M;V)$ and $J_Q$ is an almost complex structure on $Q$  fully compatible with $\pi_Q$. 
\end{itemize}
\subsubsection{Tweaking the almost complex structures}\label{tweaking}
The definition of $\EuScript{J}(M,Q)$ suffers from a predictable defect: it is
so stringent that $\EuScript{J}(M,Q)$ might not contain any regular almost complex structures. To get around this, we will take some $(J_M, J_Q)\in \EuJ(M,Q)$ and tweak it slightly, without disturbing the features which are useful to us.

\begin{figure}[t]
\centering
\labellist
\pinlabel {$B$} at 460 505
\pinlabel {$C$} at 200 505
\pinlabel {$U_B$} at 360 505 
\pinlabel {$U_C$} at 300 505
\pinlabel {$\Gamma$} at 460 565
\endlabellist
\includegraphics[scale=0.5]{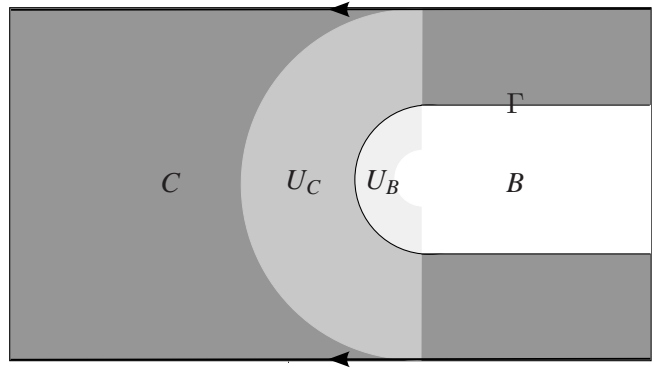}
\caption{Schematic of $S$, pictured as a rectangle with two sides identified. It is divided along $\Gamma$ into its two regions $B$ (white, with its shaded sub-region $U_B$) and $C$ (dark grey, with its lighter-shaded sub-region $U_C$).}
\label{regions}
\end{figure}

Consider the region $U_B = \{ z\in B :  \real(z)\leq 1,\; 2 |z-1| > 1 \}$ (see Figure \ref{regions}). Define $\tau \colon U_B \to \C$ by $\tau(z)= 1+ (\bar{z}-1)^{-1}$. The quotient map $q\colon \C \to \widehat{S}=\C/Z$ is injective on $\tau(U_B)$ and maps it to $C$; let $U_C\subset C$  denote $q\circ \tau(U_B)$. We shall write $\tau$ also for the induced map $U_B\to U_C$, which is a diffeomorphism, mapping $U_B\cap \Gamma$ to $U_C\cap \Gamma$ by the identity map. Notice that the closure $\overline{U_C}$ wraps all the way around the cylinder: one has $q(1+i \R)\subset \overline{U_C\cup U_B}$.

The product fibration $X:=\tau^*( P|_{U_B}) \times_{U_C} (Q|_{U_C})\to U_B$ is an LHS over $U_C$, with 2-form $(\Omega_P)_-\oplus \Omega_Q$. The pair $(J_M,J_Q)\in \EuJ(M,Q)$ induces a fully compatible almost complex structure $I$ in $X$. We consider the space $\mathrm{Tw}(J_M,J_Q)$ of fully compatible almost complex structures $I'$ in $X$ such that the support $\supp( I- I') \subset\overline{U_B}$  is contained in $U_B$. 

We then write $\EuScript{Z}(J_M,J_Q,I')$ for the moduli space of pairs $(u,v)$ where $u$ is $J_M$-holomorphic except over $U_B$, $v$ is $J_Q$-holomorphic except over $U_C$, and over $U_B$ the pair $(\tau^* u, v)$ is $I'$-holomorphic; and $\widetilde{\EuScript{Z}}(J_M,J_Q,I')$ for its enlargement to pairs with bubbles attached. We will call $\EuScript{Z}(J_M,J_Q,I')$  the moduli space of  \emph{$I'$-tweaked} pseudo-holomorphic pairs.

\subsection{Low-energy pseudo-holomorphic sections}
We consider pairs $(u,v)$ where $u$ is a $J_P$-holomorphic section of $\pi_P$, $v$ a $J_Q$-holomorphic section of $\pi_Q$, and $(u(x),v(x))\in \wh{V}_x$ when $x$ lies in $\wh{\Gamma}$. We first identify a subset of the $J_P$-holomorphic sections of $\pi_P$. The latter is, by construction, a trivial fibration, and a $J_P$-holomorphic section of $P$ is a map $u\colon B \to M$ satisfying
\[    u_s +  (\phi_t^*J_M)(u) u_t  =0 . \]
If we put $\tilde{u}(s,t)=\phi_t \circ u(s,t)$ then this equation becomes $ \tilde{u}_s + J_M(\tilde{u}_t  - X ( \tilde{u})) =0$, i.e.,
\begin{equation} 
\label{tilde u}  \tilde{u}_s + J_M \tilde{u}_t  -  \nabla K ( \tilde{u}) =0. 
\end{equation}
A class of solutions are those satisfying \begin{equation}\label{Morse trajectories}
 \tilde{u}_t =0,\quad  \tilde{u}_s =  \nabla K ( \tilde{u}) .  
\end{equation}
Let us call these \emph{gradient sections}. From equations (\ref{tilde u}, \ref{Omega P}), the action of a section $u$ equals
\begin{align} 
\notag \action(u)   =\int_B{u^*\Omega_M} 
 	&= \int_B{ \omega_M( u_s , u_t  - X (u) )\, ds\,  dt} \\
\notag 	&= \int_B{ \omega_M(\tilde{u}_s, \tilde{u}_t-X(\tilde{u}) ) \, ds\,  dt} \\
	& = \int_B{|\tilde{u}_s|^2 ds\,  dt } ,
\end{align}
so for a gradient section one has
\begin{equation}  \action(u)=   \int_B{ |( \nabla K) (\tilde{u}) |^2 ds\, dt}.  
\end{equation}
In particular, if $e(K) =\inf_u \action(u)$, where $u$ ranges over gradient sections, then $e(K)<\infty$ and $e(\lambda K)=\lambda^2 e(K)$ for $\lambda\in \R$.

We require that $u(x)\in V_x$ when $x\in \wh{\Gamma}$; for a gradient section to satisfy this condition we need $(\nabla K)\circ(\tilde{u})$ to be tangent to $V$. But from Lemma \ref{tangent}, when $J_M\in \EuJ(M;V)$, $\nabla K$ is \emph{everywhere} tangent to $V$.  

We can precisely identify the moduli space $\EuM^{\mathrm{grad}}(J_M,J_Q)$ of pairs $(u,v)\in \EuScript{Z}(J_M,J_Q)$, where $u$ is a gradient section and $v$ is \emph{horizontal}:  there is a canonical identification
\begin{equation} \label{gradient type moduli space}
\EuM^{\mathrm{grad}}(J_M,J_Q) \cong \bigcup_{x\in \crit(K|_V): \, \rho(x)\in \fix(\mu)} {W^u(K|_V;x)} , \end{equation} 
where $W^u(K|_V; x)$ is the unstable manifold of $x$ \emph{as a critical point of $K|_V$}.  This identification sends $y\in W^u(K|_V;x)$ to the unique pair $(u,v)\in \EuM^{\mathrm{grad}}(J)$ such that $u(p)=y$ and $v$ is the horizontal section corresponding to  $\rho(x) \in \fix(\mu)$. 

The unstable manifolds $W^u(K|_V;x)$ lie inside the fibres of $\rho$, and coincide with the unstable manifolds $W^u(m;x)$
of the Morse function $m$ on the fibre.
Thus 
\[ \EuM^{\mathrm{grad}}(J) \cong \fix(\mu)\times \bigcup_{x\in \crit(m)} {W^u(m;x)}  ,\]
and in particular, $ \EuM^{\mathrm{grad}}(J) $ is stratified by the Morse indices of critical points of $m$. 

The local dimension $d$ of $\EuM^{\mathrm{grad}}(J)$ is equal to $k-i(x)$, where $i(x)$ the Morse index of $x\in \crit(m)$.  The $d$-dimensional part $\EuM^{\mathrm{grad}}(J)$ is cut out transversely, and its Zariski tangent space $T_{(u,v)} \EuM^{\mathrm{grad}}(J)$ at a pair $(u,v)$ is contained in $T_{(u,v)} \EuM(J)$, the Zariski tangent space of $\EuM(J)$.
\begin{lem}
For all $(u,v)\in \EuM^{\mathrm{grad}}(J)$, one has $\ker D_{(u,v)}= T_{(u,v)} \EuM^{\mathrm{grad}}(J)$ and $\coker D_{(u,v)}=0$. Hence $\ind D_{(u,v)}=k-i(x)$.
\end{lem}
\begin{pf}
This follows from the Weitzenb\"ock formulae for Cauchy--Riemann operators in a Hermitian vector bundle.  Working in $P\to B$, with the symplectic form $\Omega_P+ \pi_P^*(ds\wedge dt)$, and considering the natural linearised symplectic connection in $\nabla^u$ in $u^*\Tv P$, we have 
\begin{align}
\label{Weitz1}
D_{u}^*D_{u} X  &= \frac{1}{2} \nabla^{u*} \nabla^{u} X + \frac{J}{2} F_\nabla(\partial_s,\partial_t) X , && X\in C^\infty_c(B; u^*\Tv P)\\
\label{Weitz2}
 D_u D_u^* \xi &= \frac{1}{2} \nabla^{u}( \nabla^u)^* \xi + \frac{J}{2} F_\nabla(\partial_s,\partial_t) \xi ,&& \xi \in \Omega^{0,1}_c(B;u^* \Tv P)
\end{align}
by a computation that can conveniently be done using the formulae of \cite[App. C]{McDuffSal}. The curvature term $F_\nabla(\partial_s,\partial_t)$, a section of the endomorphism bundle $\mathfrak{sp}(u^*\Tv P)$, was identified in (\ref{curvature2}); it vanishes identically when $F_\Omega$ does. Similar formulae apply to $D_v$ and its formal adjoint. Since both our LHFs are flat, the curvature terms in the Weitzenb\"ock formulae vanish. Now, $\ker D_{(u,v)}$ is the subspace of $\ker D_u\times \ker D_v$ on which the matching condition holds.  
Whilst (\ref{Weitz1}, \ref{Weitz2}) initially apply only to compactly supported sections, it follows by continuity that they are valid for smooth, exponentially decaying sections, in particular to elements of $\ker D_u$. Thus we see  $T_{(u,v)} \EuM^{\mathrm{grad}}(J) \subset \ker D_{(u,v)}$.  The reverse inclusion is automatic, hence equality holds.  Similarly, to show that $\coker D_{(u,v)}=0$, it suffices to show that $\ker D_{(u,v)}^*=0$. But $\ker D_{(u,v)}^*=\ker D_u^* \times \ker D_v^*$, which by the formulae is equal to $\ker \nabla_u^* \times \ker \nabla_v^*$. But $\ker \nabla_v^*=0$ because $K$ is Morse--Smale, and $\ker \nabla_u^*=0$ because $\mu$ is non-degenerate.
\end{pf}
\begin{lem}\label{grad1}
The number $\# \EuM^{\mathrm{grad}}(\bar{x}, x)_0$ of isolated gradient-type pairs asymptotic to $\bar{x}\in \fix(\mu)$ and to $x\in \wh{V}\cap \wh{V}'$ is equal to $1$ if $\rho(x)=\bar{x}$ and $x$ is a maximum for $m$, and zero otherwise.
\end{lem}
\begin{pf}
This is an immediate consequence of the structure of $\wh{V}\cap \wh{V}'$ 
(Proposition \ref{bijection}) together with the identification (\ref{gradient type moduli space}).
\end{pf}

Now we consider the gradient-type trajectories in the moduli spaces defining the homotopy $H$. The moduli spaces for the first homotopy, $H_1$, involve pairs $(u,v)$ where $v$ hits a cycle $Z_t$ in some fibre $Q_t$, where $t$ moves from $p'$ to $p$. For generically-chosen cycles, no point in these moduli spaces will be a gradient-type pair (for this amounts to a finite-dimensional intersection problem in $M$, of a codimension $k$ submanifold intersecting a $0$-manifold). The second homotopy, $H_2$, is more interesting. The contributions of gradient-type pairs to $H_2$ come from fibre products
$ \EuM^{\mathrm{grad}} \times_{\ev_p} Y  $ where $Y$ is the global angular chain bounding $Z$. One has 
\[  \EuM^{\mathrm{grad}}(\bar{x}, x) \times_{\ev_p} Y  \cong  W^u(x;m) \cap Y. \] 
By construction, the algebraic intersection number $\rho^{-1}(\bar{x}) \cap Y$ is $1$ for each $\bar{x}\in N$. So as to avoid issues concerning orientations of moduli spaces, we note that we can choose $Y$ so that  $\rho^{-1}(\bar{x}) \cap Y$ is transverse intersection in precisely one point whenever $\bar{x}\in \fix(\mu)$. We can then make the following conclusion.

\begin{lem}\label{grad2}
$  \# \EuM^{\mathrm{grad}} (\bar{x}, x) \times_{\ev_p} Y = 1$ when $\rho(x)=\bar{x}$ and $x$ is a minimum for $m$, while this number is zero in all other cases.
\end{lem}

\begin{lem}
Fix $(J_M,J_Q)\in \EuJ(M,Q)$. Then we can find a $C^\infty$-dense set of $I'\in \mathrm{Tw}(J_M,J_Q)$ so that $\EuScript{Z}(J_M,J_Q,I')\setminus \EuM^{\mathrm{grad}}(J_M,J_Q)$ is cut out transversely. 
\end{lem}
\begin{pf}
The LHFs $(P,\pi_P,\Omega_P)$ and $(Q,\pi_Q,\Omega_Q)$ have the simplifying feature that they are flat. Notice also that the set $U':=\overline{U_B\cup U_C}\subset \wh{S}$ has the property that any pair of horizontal sections $(u,v)$ over $U'$ (or rather, over its image in $\wh{S}_\Gamma$), satisfying the matching condition $\wh{\mathsf{V}}$, extends to an element of $\EuM^{\mathrm{grad}}(J_M,J_Q)$. By the argument of \cite[Lemma 2.26]{Seidel03},  any horizontal $(J_M,J_Q)$ can be tweaked by an $I'\in \mathrm{Tw}(J_M,J_Q)$ so that all non-horizontal pairs in $\EuScript{Z}(J_M,J_Q,I')$ (or indeed in $\widetilde{\EuScript{Z}}(J_M,J_Q,I')$) are regular.
\end{pf}

We now argue that it can be arranged that the contributions of the gradient-type pairs account for all the pseudo-holomorphic pairs of sufficiently low action. We need to make the notation more explicit, writing not just $J_P$ but $J_{P,K}$ for the almost complex structure on $P$ associated with the function $K$, and $\wh{\mathsf{V}}(K)$ for the matching condition associated with $K$.

\begin{lem}\label{Morse Bott}
Choose a function $K$ as in (\ref{def K}). There exists a constant $\epsilon>0$ such that the following holds. Replace $K$ by $\lambda K$ for some $\lambda \in (0, \epsilon]$. Choose $I' \in \mathrm{Tw}(J_{P,\lambda K}, J_Q)$ and suppose that $I'$ is a perturbation of size $\leq \epsilon$ in Floer's $C^\infty_{\varepsilon}$-norm.  Then, if  $(u,v)$ is any pair of finite-action, regular, $I'$-tweaked $(J_P,J_Q)$-holomorphic sections, subject to the matching condition $\wh{\mathsf{V}}(\lambda K)$, with $\action(u)+\action(v) \leq  \epsilon$, we have $(u,v)\in\EuM^{\mathrm{grad}}(J_M,J_{P,\lambda K})$.
\end{lem}

\begin{pf}
Suppose the lemma is not true. Let $J_{P,n}$ be the version of $J_P$ constructed using $\frac{1}{n}K$ in place of $K$.   Then there is a sequence of tweakings $I'_n$, perturbing $(J_{P,n},J_Q)$ by an amount which goes to zero in $C^\infty_\epsilon$ as $n\to \infty$, and of regular, $I'_n$-tweaked $(J_{P,n},J_Q)$-holomorphic pairs $(u_n,v_n)$, such that (i) $ \action(u_n)+\action(v_n) \to 0$, and (ii) neither $u_n$ nor $v_n$ is horizontal. It is better to replace $u_n$ by $\tilde{u}_n$, so that Floer's equation (\ref{tilde u}) holds, because $(\tilde{u}_n,v_n)$ is subject to a Lagrangian matching condition which is independent of $n$. Moreover, and satisfies a version of Floer's equation (\ref{tilde u}) whose coefficients converge in $C^\infty$ as $n\to \infty$. Gromov--Floer compactness implies that some subsequence of  $(\tilde{u}_n,v_n)$ converges to a pair $(\tilde{u}_\infty,v_\infty)$ such that $\action(u_\infty)=\action(v_\infty)=0$. Because the limiting energies are zero, no bubbling or trajectory-breaking is possible. Thus $u_\infty\colon B\to M$ must be a constant map, and $v_\infty$ must be horizontal.

Up to this point of the paper, we have worked exclusively with non-degenerate matched fibrations. Now we need briefly to invoke Morse--Bott methods so as to handle the case $\lambda=0$. It is well-known that the linear Fredholm theory underpinning Floer homology extends to the Morse--Bott case, provided one introduces exponentially weighted Sobolev spaces $L^{p}_{k,\delta}$ (see \cite{Don}, for instance). This more general theory extends without difficulty to our setting of matched fibrations. We take $p=2$ (so we have a Hilbert space) and take the weights $\delta$ very small, with signs chosen so that $L^{2}_{k,\delta}$-sections may increase at a mild exponential rate over the ends. Thus, if $D_0$ denotes the deformation operator at $(u_\infty,v_\infty)$, where $\lambda=0$, then $\ker D_0 \cong \R^k$ is the space of pairs which are zero over  $C$ and constant over $B$. 

We consider the 1-parametric moduli space $\EuScript{Z}^{\mathrm{par}}$ of triples $(\lambda,u,v)$, where $\lambda\in \R$ and $(u,v)$ is an index $0$ holomorphic pair for the function $\lambda K$. It is the zero set of a non-linear map $\Phi$, which takes the following shape:
\[ \Phi (\lambda , w) = 
\frac{1}{2} \left( d w+ \alpha_\lambda  \circ J \circ  \alpha_\lambda ^{-1} \circ dw \circ j
\right), \quad  \lambda \in \R. \]
Here we consider $(u,v)$ as a section $w$ of $E_\Gamma \to S_\Gamma$, while $\{\alpha_t\}$ is a 1-parameter family of automorphisms of $T E_\Gamma$ starting at $\alpha_0=\id$. Moreover, $\alpha_t$ acts as the identity on $\Th E_\Gamma$. The deformation operator $\EuScript{D}=D_{(0,w)}\Phi$ is given by
\[ \EuScript{D} (t, \dot{w}) =  D_0 \dot w + \frac{t}{2} ( \beta \circ J - J\circ \beta) \circ dw \circ j , \]
where $\beta= (d/dt)(\alpha_t)|_{t=0}$. But $\beta(\Th E_\Gamma)=0$, and so if $w=(u_\infty,v_\infty)$ is our pair of horizontal sections then $ ( \beta \circ J - J\circ \beta) \circ dw =0$.  Hence $\ker \EuScript{D} =\ker D_0 \oplus \R$. We have $(\im  \EuScript{D})^\perp = (\im D_0)^\perp\cong \R^k$.

A neighbourhood of $(0,u_\infty,v_\infty)$  in $\EuScript{Z}^{\mathrm{par}}$ can be described using a Kuranishi model. The neighbourhood is homeomorphic to a small neighbourhood of $0$ in $\kappa^{-1}(0)$, where $\kappa$ is the Kuranishi map
\[ \kappa \colon \ker \EuScript{D}\to \im(\EuScript{D})^\perp ,\quad (t, \dot{w})\mapsto \Pi\circ \Phi (t,\exp_{g_t} \dot{w}).  \]
Here $\Pi$ is orthogonal projection to $\im(\EuScript{D})^\perp$ and $g_t$ is a 1-parameter family of metrics. We have
\[ \kappa (t, \dot{w}) =   \Pi\circ \frac{t}{2} (\beta J-J\beta)\circ \dot{w}\circ j + O(|\dot{w}|^2+t^2). \]
Write $A$ for the linear map  $\ker D_0 \to \im(\EuScript{D})^\perp$,
$\dot{w} \mapsto  \Pi\circ \frac{1}{2} (\beta J-J\beta)\circ \dot{w}\circ j$. Using the inverse function theorem, we can find coordinates charts near the origins of $ \ker D_0$ and of $ \im(\EuScript{D})^\perp$ so that $\kappa \colon \R\times \ker D_0\to \im(\EuScript{D})^\perp $ takes the form
\[ \kappa(t,x) \mapsto t A(x). \] 
We argue that $A$ must be a linear isomorphism. For it maps between spaces of the same dimension, and if it had a kernel, the fibres $\EuScript{Z}^{\mathrm{par}}_\lambda$ of the projection $\EuScript{Z}^{\mathrm{par}}\to \R$ would have points of positive local dimension for small $\lambda\neq 0$. Moreover, for $n\gg 0$, $(n^{-1},u_n,v_n)$ would be such a point. This contradicts regularity of $(u_n,v_n)$. 

Thus $A$ is an isomorphism, and there exist coordinates in which $\kappa$ is given by the map
\[   \R\times \R^k  \to \R^k,  \quad (\lambda, x)\mapsto \lambda x,  \]
and $\pr_2 \colon \R\times\R^k \to \R$ corresponds to the projection $\EuScript{Z}^{\mathrm{par}}\to \R$. The strand $\R \times \{0\} \subset \kappa^{-1}(0)$ parametrises gradient sections. For large $n$, $(u_n,v_n)$ lies in this strand, which is what we want.
\end{pf}

\section{Proof of the Gysin sequence: algebraic aspects}\label{alg}

\subsection{$\R$-graded homological algebra}
The algebraic mechanism we use to prove Theorem \ref{Quasiiso} is closely related to that used by Seidel to establish the exactness of the sequence describing the effect of Dehn twists on Floer homology. However, because the the symplectic action functional is not exact and the action spectrum not necessarily discrete, two new ingredients are needed: a completion (equivalent to linearity over the Novikov ring $\Lambda_R$), and the observation that the low-order terms are defined over the base ring $R$.

If $R$ is a commutative unital ring, an \emph{$\R$-graded $R$-module} is a $R$-module $V$ with a given direct sum decomposition $V= \bigoplus_{r\in \R}{V_r}$. Its \emph{support} is $\supp(V):= \{r\in \R: V_r \neq 0\}$. 
Given an interval $I\subset \R$, we say that $V$ has \emph{gap} $I$ if $r,s \in \supp(V)$ implies $|r-s| \notin I$. A homomorphism $f\colon V\to V'$ between $\R$-graded modules has \emph{order} $I$ if $f(V_r)\subset \bigoplus_{s\in I}{V'_{r+s}}$ for all $r$. 

The following lemma \cite[Lemma 2.31]{Seidel03} is Seidel's variation on a well-known principle.
\begin{lem}
Suppose that $(D,\delta)$ is a finitely-supported $\R$-graded $R$-module with gap $[\epsilon,2\epsilon)$ for some $\epsilon>0$, equipped with a differential of order $[0,\infty)$; thus we may write $\delta=\delta_{\mathrm{low}}+\delta_{\mathrm{high}}$, where $\delta_{\mathrm{low}}$ is a differential of order $[0,\epsilon)$, and $\delta_{\mathrm{high}}$ a map of order $[2\epsilon,\infty)$. Then, if $\delta_{\mathrm{low}}$ is acyclic,  so is $\delta$; that is, $H(D,\delta_{\mathrm{low}})=0$ implies $H(D,\delta)=0$.
\end{lem} 
One can weaken the finite-support condition by demanding instead that $\supp(D)$ is bounded above, but the condition cannot be dropped altogether: let $D= R[x]\oplus R[y]$. Let $D_r$ be zero except when $r\in \Z_{\geq 0}$, in which case it is spanned by $x^r$ and $y^r$. Let $\delta(f(x),g(y))=(0, (1-y)f(y))$. Then (taking $0<\epsilon<1$) $\delta_{\mathrm{low}}$ is acyclic, but $\delta$ is not: the homology $H(D,\delta)$ is one-dimensional, spanned by $(0,1)$. 

The \emph{completion} $\hat{V}=\varprojlim_{r}{V/V_{\geq r}}$ of the $\R$-graded module $V$ is the submodule of $\prod_{r}{V_r}$ consisting of functions $v \colon \R \to V$,   with $v(r)\in V_r$, such that $\supp(v)$ has the `Novikov' property that $\supp(V) \cap (-\infty,c]$ is finite for any $c\in \R$. Completion is functorial: a homomorphism $f$ of degree $[0,\infty)$ extends `by continuity' to a homomorphism $\hat{f}$ between completions. We can still speak of a map $F \colon \hat{V} \to \hat{V'}$ having order $I$: we mean that, for all $r$, we have
\[ F(\hat{V}_r) \subset  \big( \bigoplus_{s\in I}{V'_{r+s} }\big)\, \hat{} ,\]
where $\hat{V}_r\subset \hat{V}$ is the image of $V_r$ under $V\to \hat{V}$, and the hat on the right-hand side denotes completion inside $\hat{V}'$. 

What remains of the lemma when one drops the finiteness (or boundedness from above) hypothesis on $\supp(D)$, but retains the gap assumption, is that $H(D,\delta_{\mathrm{low}})=0$ implies $H(\hat{D},\hat{\delta})=0$. (In the example with $D= k[x]\oplus k[y]$, one has $\hat{D}=k[[x]]\oplus k[[y]]$ and $(0,1)=\hat{\delta}(0,\sum_{n\geq 0}{y^n})$.) This assertion readily follows from Seidel's lemma and the following observation.
\begin{lem}
Let $(D,\delta)$ be a differential $R$-module, $\{ F^q D \}$ a decreasing filtration. Suppose that $H(D/F^q D)=0$ for all $q$. Then $H(\varprojlim_q{D/F^q D})=0$.
\end{lem}
The proof is a simple exercise. (An unnecessarily fancy proof is to observe that both$ \{D/F^q D\}$ and $\{H(D/F^q D)\}$ form inverse systems satisfying the Mittag--Leffler condition, so $\varprojlim^1 {H(D/F^q D)}= \varprojlim {H(D/F^q D)} =0$, and hence $H(\varprojlim{D/F^q D} )=0$, cf.  \cite{Weibel}.)

If one weakens the gap assumption, precautions are needed to safeguard even the conclusion that the completed $\R$-graded module is acyclic.

For $A$ an $\R$-graded module, let $\underline{A}$ denote the direct sum of shifts $\bigoplus_{r\in \R}{A[r]}$. If $A$ has a differential $d$ then $\underline{A}$ then has an induced differential $\underline{d}$ which respects the direct sum and restricts to $A[r]$ as $d$.

\begin{lem}
Suppose that $(A,d)$ is a finitely-supported $\R$-graded module with gap $[\epsilon,2\epsilon)$, equipped with a differential of order $[0,\epsilon)$, for some $\epsilon>0$.
Let $D=\underline{A}$. Let $\delta$ be a differential on $\hat{D}$ such that
\begin{enumerate}
\item[(i)] 
$\delta$ is invariant under the $\R$-action on $\hat{D}$ by shifts and continuous with respect to the completion.
\item[(ii)]
$\delta$ maps $A\subset \hat{D}$ into $\widehat{\bigoplus}_{r\geq 0}{A[-r]}$.
\item[(iii)]  
$\delta=\delta_{\mathrm{low}}+\delta_{\mathrm{high}}$, where $\delta_{\mathrm{low}}$ is equal to $\widehat{\underline{d}}$, the differential induced by $d$, and $\delta_{\mathrm{high}}$ is a homomorphism of order $[2\epsilon,\infty)$. 
\item[(iv)] $H(A,d)=0$. 
\end{enumerate}
Then $H(\hat{D},\delta)=0$.
\end{lem}
\begin{rmk}
We can think of $\underline{A}$ as the tensor product $A \otimes_R R[\R]$, so that elements of  $A[-r]$ can be written as $a \otimes t^r$, $a\in A$. When $A$ has bounded support, the completion $\underline{A}\, \hat{}$ of $\underline{A}$ is naturally identified with $A\otimes_R \Lambda_R$. Then (i) just means that $\Lambda_R$-linear, and $\widehat{\underline{d}}$ is the $\Lambda_R$-linear extension of $d$. 
\end{rmk}

\begin{pf}
Assumption (iv) implies that $H(\underline{A},\underline{d})=0$. It then follows from (iii) that $H(D,\delta_{\mathrm{low}})=0$.
Set
\[ F^p D = \widehat{\bigoplus}_{r\geq p\epsilon}{A[-r]} \] 
Thus $F^p D \supset F^{p+1}D$ and, by (i) and (ii), $\delta(F^p F)\subset F^p D$. This decreasing filtration of $D$ induces a singly-graded spectral sequence $(E^*_k, d_k)_{k\geq 0}$ with $E^*_1 = H(\mathrm{gr}_*(D))$ (where $\mathrm{gr}_*$ denotes the associated graded module $\bigoplus_p{F^p D/F^{p+1}D}$). It follows from (iii) that $E^1_*=0$. Indeed, $A$ is the direct sum of `clumps' $C^i$, preserved by $d$, whose support has width $<\epsilon$, and $F^p D/F^{p+1}D$ has corresponding clumps $\bigoplus_{r\in [p\epsilon,(p+1)\epsilon)}{C^i[-r]}\subset \bigoplus_{r\in [p\epsilon,(p+1)\epsilon)} A=F^pD /F^{p+1}D$. By (iii), the differential on the clumps is induced by the differential $\delta_{\mathrm{low}}$ on $C_i$, which is acyclic. A spectral sequence argument as in \cite[Lemma 2.31]{Seidel03} then shows that $F^p D/F^{p+1}D$ is acyclic.

The spectral sequence $E^*_k$ might not converge. However, if we look at the quotient $D/ F^q D$, with its bounded-below, exhaustive filtration induced by the filtration of $D$, we obtain a spectral sequence converging to $H(D/F^q D)$. The $E^*_1$ page is still zero. Thus $H(D/F^q D)=0$. By the previous lemma, $H(\varprojlim_q {D/F^qD})=0$, i.e., $H(\hat{D})=0$.
\end{pf}
Let us give property (ii) a name:  given $R$-modules $A$ and $B$ and a $\Lambda_R$-linear map $f\colon A\otimes \Lambda_R\to A\otimes \Lambda_R$, let us say $f$ is \emph{positive} if $f(a)=\sum{b_i t^{r_i}}$ for elements $b_i\in B$ and $r_i\geq 0$.

The following lemma is a generalisation of a rotated version of \cite[Lemma 2.32]{Seidel03}.
\begin{lem}[Double mapping cone lemma.]\label{double}
Fix  $\epsilon>0$. Suppose that $A$, $A'$ and $A''$ are finitely supported $\R$-graded $R$-modules satisfying
\begin{enumerate}
\item[{\sc gap}:] $A'$ and $A$ have gap $(0,3\epsilon)$ while $A''$ has  gap $(0,2\epsilon)$. If $r\in \supp(A')$ and $s\in \supp(A)$ then $|s-r|\notin (0, 4\epsilon)$.
\end{enumerate}
Let  $C=A\otimes_R \Lambda_R$, $C'=A'\otimes_R \Lambda_R$ and $C''=A''\otimes_R \Lambda_R$, and suppose that these three $\R$-modules are equipped with $\Lambda_R$-linear differentials $d_C$, $d_{C'}$ and 
$d_{C''}$. Suppose that
\[  \begin{CD} C' @>{b}>> C@>{c}>> C''   \end{CD} \]
is a sequence of differential maps, and $ h \colon  C' \to C''$ 
a null-homotopy of $c\circ b$. Assume  
\begin{enumerate}
\item[{\sc pos}:] The maps $d_C$, $d_{C'}$, $d_{C''}$, $b$, $c$ and $h$ are all positive.
\item[{\sc ord 1}:] The differentials $d_C$, $d_{C'}$ and $d_{C''}$ each have order $[2\epsilon, \infty)$. The map $b$ has order $[2\epsilon,\infty)$; $c$ and $h$ have order $[0,\infty)$. 
\item [{\sc ord 2}:]
One can write $c=c_{\mathrm{low}}+c_{\mathrm{high}}$ and $h=h_{\mathrm{low}}+h_{\mathrm{high}}$, where the linear maps $c_{\mathrm{low}}$ and $h_{\mathrm{low}}$ have order $[0,\epsilon)$ while $c_{\mathrm{high}}$ and $h_{\mathrm{high}}$ have order  $[ 2\epsilon,\infty)$.
\end{enumerate}
Further assume
\begin{enumerate}
\item [{\sc low 1}:] the maps $c_{\mathrm{low}}$ and $h_{\mathrm{low}}$ are induced by maps $c_0\colon A\to A''$ and $h_0\colon A'\to A''$, respectively; and
\item [{\sc low 2}:] the map $(h_0,c_0) \colon A \oplus A'  \to A''$ is a linear isomorphism.
\end{enumerate}
The induced differential map 
\[  (h, c) \colon \cone(b)  \to C'',  \] 
is certainly positive and has order $[0,\infty)$ when we consider the mapping cone as the completion of the direct sum $C'\oplus C$ of $\R$-graded modules.  The conclusion of the lemma is that $(h,c)$ is a quasi-isomorphism.
\end{lem}
\begin{pf}
We shall invoke the last lemma to prove that the `double mapping cone' $\cone(h,c)$ is acyclic. As an $\R$-module, $\cone(h,c) = C ' \oplus C \oplus C''= (A'\oplus A\oplus A'')\otimes \Lambda_R$, so requirement (i) of the previous lemma holds. Its differential $ \delta$ has the block form
\[ \delta =  \left[ \begin{array}{ccc}
	 d_{C'}	& 	0 	&	0 \\
	 b		&	-d_C &	0 \\
	- h		& 	c	&	d_{C''}		 
\end{array} \right] , \]
and is positive by {\sc pos}, so (ii) holds. We use {\sc ord 2} to decompose the differential as  $\delta= \delta_{\mathrm{low}}  + \delta_{\mathrm{high}}$, where
\[ \delta_{\mathrm{low}}= \left[ \begin{array}{ccc}
	0 	 	& 	0 	&	0 \\
	0		&	0 	&	0 \\
	-h_{\mathrm{low}}		& 	c_{\mathrm{low}}	&	0		 
\end{array} \right]  , \quad 
\delta_{\mathrm{high}} =  \left[ \begin{array}{ccc}
	d_{C'}	& 	0 	&	0 \\
	b		&	-d_C &	0 \\
	- h_{\mathrm{high}}	& 	c_{\mathrm{high}}	&	d_{C''}		 
\end{array}  \right] . \]
By {\sc ord 1} and {\sc ord 2}, $\delta_{\mathrm{low}}$ has order $[0,\epsilon)$, and squares to zero, while $\delta_{\mathrm{high}}$ has order $[2\epsilon,\infty)$.  By {\sc low 1}, $\delta_{\mathrm{low}}$ is induced by a differential $d$ on $A' \oplus A\oplus A''$.
The {\sc gap} conditions imply that $A' \oplus A \oplus A''$ has gap $[\epsilon, 2\epsilon)$. This deals with (iii).  Finally, (iv) holds because of assumption {\sc low 2}, and the proof is complete.
\end{pf}

\subsection{Completing the proof}
All the ingredients in our proof of Theorem \ref{Gysin} are now to close to hand.   
Recall from the introduction that the Gysin sequence is derived from the exact sequence of a mapping cone: indeed, it is an immediate consequence of Theorem \ref{Quasiiso}. The definition of $\rho_*$ is algebraic. To prove Theorem \ref{Quasiiso} we shall invoke the double mapping cone lemma (\ref{double}), taking the following as our input data.  
\begin{itemize}
\item 
$ A  = A' = R \fix(\phi_h^{-1}\circ \mu)$, with $\R$-gradings concentrated in degree zero.
\item
$A'' =  R (\wh{V} \cap \wh{V}' )$, with $\R$-gradings still to be specified. 
\item
The differential $d_C=d_{C'}$ on $C=C'=A\otimes_R \Lambda_R=\CF(\mu)$ is the Floer-theoretic differential, for some regular almost complex structure. Likewise, $d_{C''}$ is defined as a Floer-theoretic differential on $\CF(\wh{V},\wh{V}')$. The maps $c$ and $h$ were constructed in Section \ref{geom}.
\item
From Proposition \ref{bijection}, we have a canonical bijection between $\wh{V}\cap\wh{V}'$ and $\fix(\mu)\amalg\fix(\mu)$. We assign $\R$-degrees to elements of $\wh{V}\cap\wh{V}'$ so that the map $(h_0,c_0)$ induced by this bijection preserves degree.
\item
We choose $\epsilon>0$ as follows.
First, as a consequence of Gromov--Floer compactness, there is some $\epsilon_0>0$ so that $d_C$, $d_{C''}$ and $b$ all have order $[2\epsilon_0,\infty)$. There is $\epsilon_1$ so that $A''$ has gap $(0,2\epsilon_1)$.  By Lemma \ref{Morse Bott}, there is an $\epsilon_3>0$ so that, for any pair $(u,v)$ contributing to $c$ or $h$ which has action $\leq 2\epsilon_3$, $u$ is a gradient section and $v$ is horizontal. Finally, there is $\epsilon_4$ so that any gradient section $v$ has action $\action(v)\leq \epsilon_4/2$. We take $\epsilon < \min(\epsilon_0,\epsilon_1,\epsilon_2,\epsilon_3,\epsilon_4)$.
\end{itemize}
Now let us check the conditions of the double mapping cone lemma (\ref{double}).
The assumptions {\sc gap} and {\sc ord 1} hold by choice of $\epsilon$. 
The differentials $d_C=d_{C'}$ and $d_{C''}$ and the map $b$ are positive by construction.  Because $\epsilon<\epsilon_3$,  Lemma \ref{Morse Bott} tells us that $c_0$ and $h_0$ precisely capture the low (meaning $<\epsilon$) action contributions to $c$ and $h$, i.e.  {\sc low 1} holds. It follows that $h$ and $c$ are positive, so {\sc  pos} holds.\footnote{This explains why we chose to concentrate $A$ and $A'$ in degree 0  rather than grading them by representatives the action functional.} Moreover, {\sc ord 2} holds because $\epsilon<\epsilon_4$. The crucial condition {\sc low 2} comes from Lemmas \ref{bijection}, \ref{grad1} and \ref{grad2}.

This completes the proof of Theorem \ref{Gysin}. The first item in Addendum \ref{add1} is, at this point, a triviality: it is clear that the maps $e$ and $\rho^!$ intertwine the action by quantum cap product by $QH^*(N)$ (in the case of $e$, this is just the fact that the cap product on fixed point Floer homology \emph{is} an action of the small quantum cohomology algebra, and that this algebra is super-commutative). Hence so does the connecting map $\rho_*$. The remaining item of  this addendum, concerning orientations, will be addressed in the next section.

As to Addendum \ref{add2}, suppose that $\mu=\id$ and that  $\HF(\wh{V},\wh{V})\cong H_*(\wh{V};\Lambda_R)$. We perturb $\mu$ to a Hamiltonian automorphism $\phi_H$, generated by a $C^2$-small Hamiltonian $H$. Then $\Fix(\mu)=\crit(H)$, and the Floer complex $\CF(\phi_H)$ is canonically identified with the Morse complex for $H$. Similarly, $\CF(\wh{V}, (\phi_K\times \phi_H)\wh{V})$ is identified with the Morse complex for $H\times K$ restricted to $V$. We have established that the low-action contributions to $C\rho^!$ and $e$ are  Morse-theoretic; they compute the transfer map $\rho^!$ on Morse homology, and the cap product by $e(V)$ on Morse homology. This gives the result.

\section{Refinements}\label{refinements}
\subsection{On the borderline}
We now analyse the breakdown of the Gysin sequence in the borderline case where {\sc mon} holds but {\sc mas} just fails because $m^{\min}_{\wh{V}}=k+1$. We saw at the outset (Example \ref{sphere ex}) that Theorem \ref{Gysin} can then also fail.

Consider the moduli space $\widetilde{\EuScript{N}}$ of parametrized $J$-holomorphic discs $\delta \colon (D,\partial D) \to (M_-\times N,\wh{V})$ of Maslov index $k+1$.  

Because index $k+1$ discs have the smallest positive action among holomorphic discs attached to $\wh{V}$, $\EuScript{N}:=\widetilde{\EuScript{N}}/\aut(D)$ is compact. For generic $J$, it is also smooth of dimension $(k+1)+(2n-k) = 2n+1$.  According to Kwon--Oh \cite{KwonOh}, such discs must also be somewhere injective, which means that `transversality of evaluation' holds:  if we fix smooth, closed submanifolds $A$ and $B$ of $\wh{V}$ (or more generally smooth singular cycles), of respective codimensions $\alpha$ and $\beta$, the evaluation maps $\ev_1$ and  $\ev_{-1} \colon \widetilde{\EuScript{N}} \to \wh{V}$ are, for generic $J$, transverse to $A$ and $B$ respectively. Thus $\ev_1^{-1}(A)\cap \ev_{-1}^{-1}(B) $ is  smooth of dimension $2n+1-\alpha-\beta$. Its quotient by the subgroup $\R\subset \aut(D)=PSU(1,1)$, consisting of those $\alpha\in \aut(D)$ with $\alpha(1)=1$ and $\alpha(-1)=-1$, is compact. The action is free, since translation-invariant holomorphic discs are trivial, and so $\ev_1^{-1}(A)\cap \ev_{-1}^{-1}(B) /\R$ is compact, and smooth of dimension $2n-\alpha-\beta$.

In particular, if we take $A=\{x\}$ and $B=\rho^{-1}Z$, where $Z$ represents the Poincar\'e dual to $e(V)$, then $\dim \ev_1^{-1}(A)\cap \ev_{-1}^{-1}(B) /\R = 2n -(2n-k)-(k+1)=-1$, so the intersection is empty.

Now consider a global angular chain $Y$ bounding $\rho^{-1}(Z)$ (as in Subsection \ref{global angular}). The space
$(\ev_1^{-1}(\{x\})\cap \ev_{-1}^{-1}(Y)) /\R$ is a compact $0$-manifold. 
\begin{defn}
Define
$ \nu_Y= \# ( (\ev_1^{-1}(x)\cap \ev_{-1}^{-1}(Y)) /\R) \in \Z/2$.
\end{defn}
Note that $\nu_Y$ is independent of $x$ and of (regular) $J$.  For example, $\nu_Y=0$ if, for some $J$, there is a point $x\in \wh{V}$ which does \emph{not} lie on a $J$-holomorphic disc of index $k+1$. We do not address the dependence of $\nu_Y$ on $Y$. 
\begin{thm}\label{borderline}
Suppose that {\sc mon} holds and $m^{\min}_{\wh{V}}=k+1$. The conclusions of Theorem \ref{Gysin} hold when one replaces the map $e$ by $e+ t^{(k+1)\lambda} \nu_Y\,\id$, where $\lambda$ is the monotonicity constant of $\wh{V}$. 
\end{thm}

\begin{pf}
The construction of the map $e$ obviously goes through. Less obviously,  so does that of $\rho^!$: limits of sequences of sections of index 0 or 1 cannot bubble off discs, since these have index $\geq k+1 \geq 2$.  We will attempt to show that $\rho^! \circ e$ is nullhomotopic and see that we encounter an obstruction. 

The first chain homotopy, $H_1$, can still be defined as a map, because the relevant sections have index $k$.  When we come to prove that $H_1$ is a homotopy, a potential problem arises: in the relevant moduli spaces $\EuScript{Z}$ of index $k+1$ sections, a disc $\delta \colon (D,\partial D) \to (M_-\times N,\wh{V})$ of Maslov index $k+1$ may bubble off. Indeed, $\EuScript{Z}$ has extra ends parametrising limiting objects which consist of
\begin{itemize}
\item a pair $(u,v)$ of pseudo-holomorphic sections of index $0$; and
\item a bubble $\delta \colon (D,\partial D) \to (M_-\times N,\wh{V})$ such that $\delta(1)=v(p)$ and $\delta(-1) \in \im \rho^{-1}Z$. 
\end{itemize}
However, for generic almost complex structures, 
\[ \ev_1^{-1}(v(p))\cap \ev^{-1}_{-1}(\rho^{-1}Z)=\emptyset,\] 
so the problem is illusory and $H_1$ is a chain homotopy after all.

The second chain homotopy, $H_2$, is constructed via index $k$ sections (again, no bubbles). When try to prove that it is a homotopy by analysing a 1-dimensional moduli space we must again consider index $k+1$ sections, so  bubbling is possible. The moduli space now has ends corresponding to configurations consisting of:
\begin{itemize}
\item A pair $(u,v)$ of pseudo-holomorphic sections of index $0$; and
\item a bubble $\delta \colon (D,\partial D) \to (M_-\times N,\wh{V})$ such that $\delta(1)=v(p)$ and $\delta(-1) \in \im Y$. 
\end{itemize}
These configurations contribute the following boundary components to the 1-dimensional  moduli space:
\begin{equation}\label{boundary}
 \bigcup_{(x_-,x_+)} {\bigcup_{(u,v)\in \EuScript{Z}(x_-,x_+)^0}{ \EuScript{N}(v(p),Y) }. } 
 \end{equation} 
This follows from the gluing theorem for holomorphic discs proved in Theorem \cite[4.1.2]{BirCor}.  When $\nu_Y=0$, there is for each $x_-\in \fix(\nu)$ and each $x_+\in \wh{V}\cap \wh{V}'$, an even number of extra boundary components in the part of the moduli space consisting of sections asymptotic to $(x_-,x_+)$. Hence $h:=H_1+H_2$ is still a nullhomotopy.  In general, the appearance of the boundary (\ref{boundary}) implies that $\partial h +  h \partial =  \rho^! \circ e + \nu_Y t^{ (k+1) \lambda} \, \rho^!$ (the Novikov weight $\lambda(k+1)$ is the area of the minimal Maslov-index discs). That is, $h$ is a nullhomotopy of $\rho^!+  t^{\lambda(k+1)} \nu_Y \id$.

The analysis of low-action sections goes through unchanged. We may assume that $(k+1)\lambda \gg \epsilon$, which means that the new term  $\nu_Y t^{ (k+1) \lambda} \, \id$ counts as a `high-action' term. The algebraic part of the proof then goes over unchanged.
\end{pf}

\subsection{Compatibility with quantum cap product}
The quantum cohomology $QH^*(N)=H^*(N;\Lambda_{\Z/2})$ acts on $\HF_*(\mu)$ by quantum cap product. This makes it a module over the algebra $QH^*(N)$ with its quantum product. Since the map $e$ is itself defined as the quantum cap product by the Euler class of $\rho\colon V\to N$, and since $QH^*(N)$ is (super) commutative, $e$ intertwines the action of $QH^*(N)$. 

The Lagrangian Floer homology $\HF(\wh{V},\wh{V}')$ is a module over $QH^*(M_-\times N)$, and hence also over $QH^*(N)$ via the pullback $\pr_2^*\colon QH^*(N)\to QH^*(M_-\times N)$ (which is an algebra homomorphism). It follows easily from the definitions that $\rho^!$ intertwines the actions of $QH^*(N)$.

\subsection{Orientations}\label{ori}
To define Floer homology over $\Lambda_\Z$, one needs to specify a system of \emph{coherent orientations}. This can be reduced to topology: to specify coherent orientations for the Lagrangian Floer homology  $\HF_*(L_0,L_1)$ of $L_0,L_1\subset P $, it suffices to give \emph{relative spin-structures} on the two Lagrangians, that is, to give a stable vector bundle $\xi\to P$, together with stable spin-structures on $\xi|_{L_i} \oplus TL_i$ for $i \in \{0, 1\}$ \cite[Chapter 6]{FOOO}, \cite{Seidel06}. 

The fixed point groups $\HF_*(\mu)$ are defined over $\Lambda_\Z$ in a canonical way, as explained in \cite[II, Section 3.7]{PerutzI}. 

In the setting of Lagrangian matching conditions, one wants the cobordism-maps to be well-defined over $\Lambda_{\Z}$. This turns out to be rather straightforward: one takes a vector bundle $\eta\to S_\Gamma$ over the surface $S_\Gamma$, and considers its restrictions $ \eta|_{\Gamma_+}$ and $ \eta|_{\Gamma_- }$ as bundles over $\Gamma$.  Take $\xi\to \Gamma$ to be their direct sum $\xi= \eta|_{\Gamma_+}\oplus \eta|_{\Gamma_- }$. Giving a spin-structure for $TQ\oplus \xi|_Q$, where $Q$ is the Lagrangian boundary condition. This choice defines a lift of the cobordism-maps over $\Lambda_\Z$(we refer to \cite{WW08} for a full discussion). Note that one can always work with the canonically-defined groups  $\HF_*(\mu)$ (again, cf. \cite[II, Section 3.7]{PerutzI}). 

When $V \subset M$ is an $S^k$-fibred coisotropic submanifold, the natural identification $\wh{V}\to V$ (by projection to $M$) induces an isomorphism $T\wh{V} \cong TV$ which we shall consider as an equality.  If we put $\xi = p_1^* TM \to M\times N$, we will have
$  \xi |_{\wh{V}} \oplus T\wh{V} =   TM|_V  \oplus  T V$. Thus there is a canonical short exact sequence
\[ 0\to TV \oplus TV \to \xi |_{\wh{V}} \oplus T\wh{V}  \to N_{V/M}\to 0  \]
(the sequence splits as soon as one chooses a Riemannian metric on $M$ and so identifies $N_{V/M}$ with $TV^\perp$). Now, a direct sum $U\oplus U$ of two copies of the same $\mathrm{O}(r)$-vector bundle always has a canonical spin-structure because its structure group reduces to the diagonal subgroup $\mathrm{O}(r) \subset \SO(2r)$, and the inclusion homomorphism $\mathrm{O}(r)\to \SO(2r)$ lifts to $Spin(2r)$ since it kills $\pi_1$. Applying this last principle to $U=T V$, we see that to give a spin-structure on $ \xi |_{ \wh{V}} \oplus T\wh{V} $ it suffices to give a spin-structure $\sigma$ in $N_{V/M}$. 

Now, $\sigma$ also induces a spin-structure $(\id\times \mu)^* \sigma$ on the normal bundle to $\wh{V}'= (\id\times \mu)\wh{V}$, and hence gives rise to stable spin structures on both $\xi|_{\wh{V}}\oplus T\wh{V}$ and $\xi|_{\wh{V}'}\oplus T\wh{V}'$.  Thus, \emph{given $\sigma$, the Floer homology $\HF(\wh{V},\wh{V}')$ is canonically defined with $\Lambda_R$-coefficients for $R$ an arbitrary commutative ring.}

\begin{lem}
When $V$ is $S^1$-fibred, a co-orientation for $V$ induces a relative spin structure for $\wh{V}$. When $V$ is the unit-sphere bundle $S(E)$ in a Euclidean vector bundle $E\to N$,  a spin-structure in $E$ induces a relative spin structure for $\wh{V}$.
\end{lem}
\begin{pf}
In the $S^1$-fibred case, a co-orientation trivialises $N_{V/M}$, and so gives a canonical spin structure on $N_{V/M}$. When $V=S(E)$, we have $\Tv V \oplus \underline{\R} \cong \rho^*E$, and hence $w_2(N_{V/M})=w_2(\Tv V) = \rho^*w_2(E)$. By the classical Gysin sequence, this class vanishes if and only if $\mathrm{rank}(E)=2$ or $w_2(E)=0$. In the latter case, a spin-structure in $E$ induces a stable spin-structure on $\Tv V$. But $N_{V/M}\cong (\Tv V)^*$ via the symplectic form, so we get a stable spin-structure in $N_{V/M}$.
\end{pf}

Our conclusion (which is the second item of Addendum \ref{add1}) is as follows.

\begin{prop}
Suppose that the hypotheses of Theorem \ref{Gysin} hold. Assume \emph{either} that  $V$ is a co-orientable hypersurface in $M$, \emph{or} that $V=S(E)$ with $w_1(E)=0=w_2(E)$. Then $\HF(\wh{V},(\id\times \mu)\wh{V})$ can be defined over $\Lambda_R$ for any commutative unital ring $R$ in such a way that the conclusions of Theorem \ref{Gysin} hold.
\end{prop}
(The bundle $\eta\to S_\Gamma$ required for the definitions of $C\rho$ and $h$ can be specified by putting $\eta|_B= \Tv Q$ and $\eta|_C= 0$.

\section{Connected sums of three-manifolds}\label{3manifolds}
This final section explains the (conjectural) connection between the symplectic Gysin sequence and gauge theory on 3- and 4-manifolds.

\subsection{Connected sums and indefinite singularities}
Let $Y_1$ and $Y_2$ be closed, oriented smooth 3-manifolds. When $f_1\colon Y_1 \to S^1$ and $f_2 \colon Y_2\to S^1$ are \emph{harmonic Morse functions}, that is,  circle-valued Morse functions with only indefinite critical points, the connected sum $Y_1\, \# \,Y_2$ inherits a harmonic Morse function with $(\# \crit(f_1) + \#\crit (f_2) + 2)$ critical points. Indeed, let $I \subset S^1$ be a closed interval containing only regular points of $f_1$ and $f_2$.  Trivialise the two fibrations over $I$, identifying them with $I \times \Sigma_1 \to I$ 
and $I \times \Sigma_2 \to I$. There is an elementary cobordism $W$ from $\Sigma_1 \amalg \Sigma_2$ to $\Sigma_1\#\Sigma_2$, carrying a Morse function with a single critical point $c_2$, of index $2$.  Likewise, there is  an elementary cobordism $W'$ from $\Sigma_1\#\Sigma_2$
to $\Sigma_1 \amalg \Sigma_2$, carrying a Morse function with a single critical point $c_1$ of index $1$. Glue $W$ to $W'$ along $\Sigma_1 \# \Sigma_2$ to obtain a 3-manifold $V$ with a Morse function $V \to I$. Define a closed 3-manifold  by gluing $V$ into $ f_1^{-1}(S^1\setminus \interior(I)) \amalg   f_2^{-1}(S^1\setminus \interior(I)) $ in the obvious way. Denote this new manifold by $Y_1 \#_{I} Y_2 $. It comes with a Morse function $f_1 \#_{I} f_2\colon  Y_1 \#_{I} Y_2 \to S^1$, with two `new' critical points, both over $\interior{I}$. It is easy to check that  $Y_1 \#_{I} Y_2 $ is diffeomorphic to the connected sum $Y_1 \, \# \,Y_2$;  there is a disconnecting 2-sphere $S\subset Z$ lying over the interval $I$, the union of the ascending disc of $c_2$ and the descending disc of $c_1$. 

There is a natural cobordism $X$ from $Y_1 \amalg Y_2 $ to $Y_1 \#_{[a,b]} Y_2$, together with a function $F\colon X\to S^1 \times [0,1]$ such that $F|_{(Y_1 \amalg Y_2) } = (f_1 \amalg f_2,0 ) $ and $F|_{ Y_1 \#_{[a,b]} Y_2} = (f_1 \#_{[a,b]} f_2,1)$. To build it, begin with $W$. This carries a Morse function $q$ with just one critical point. We normalise $q$ so that $q(\Sigma_1 \amalg \Sigma_2)=\{0\}$, $q(\Sigma_1 \# \Sigma_2 ) = \{1\}$, and $q(W)=[0,1]$. Let $H = \{z: |z| \leq 1, \imag(z) \geq 0 \} \subset \C$ be a closed half-disc, and let $m \colon H\to [0,1]$ be the modulus function, $z\mapsto |z|$. Let $X_0 = m^* W$; so $X_0$ carries a natural map $F_0\colon  X_0\to H$. 

Take the trivial cobordism $(Y_1 \amalg Y_2) \times [0,1]$, and let $F_1 = (f_1 \amalg f_2, \id)\colon (Y_1 \amalg Y_2) \times [0,1]  \to  S^1 \times  [0,1] $. Embed $H$ into $[a,b]\times [0,1]$  by an affine linear map so that its straight edge is mapped to $[a,b]\times \{1\}$. Thereby embed $H$ into $S^1\times[0,1]$. Trivialise $F_1$ over $H$ (extending the existing trivialisation over the straight edge). Now define $X$ by excising $F_1^{-1} (H)$ from $(Y_1 \amalg Y_2) \times [0,1]$ and gluing in $X_0$ in its place, in a way which should be clear. Thus
\[   X = X_0 \cup \left (  (Y_1 \amalg Y_2) \times [0,1]  \setminus F^{-1}(H) \right). \] 
The map $F\colon X\to S^1\times[0,1]$ is obtained by gluing $F_0$ and (the restriction of ) $F_1$.

The function $\pr_2 \circ F \colon X\to [0,1]$ has just one critical point, namely, the unique point in the arc $\crit(F_0)$ which lies over the 
image of $\ii \in H$ under its embedding into $S^1\times [0,1]$. This critical point $c$ is non-degenerate of index 1. Hence $X$ is the elementary cobordism from the disjoint union to the connected sum.

The critical manifold of the map $F\colon X \to S^1\times[0,1]$ is an `indefinite fold', mapping injectively to the base; the local model is $(t;x_1,x_2,x_3)\mapsto (t; x_1^2+x_2^2-x_3^2)$. This makes it a simple example of a \emph{broken fibration} in the sense of \cite{PerutzI}. 

We shall concentrate on the case where $f_1$ and $f_2$ both have empty critical set, hence are fibrations. In this case, the set of critical values of $F$ is an arc in the cylinder: it resembles the arc $\Gamma$ from Figure \ref{gysincyl}.

\subsection{Symmetric products}
This section describes a class of examples of the symplectic Gysin sequence.
In these examples, $M$ is $\sym^n(\Sigma)$, the $n$th symmetric product of a closed Riemann surface $\Sigma$, for some $n\geq 1$. The symplectic form $\omega_\lambda$ on $M$ is a K\"ahler form representing one of the classes $\eta_\Sigma + \lambda \theta_\Sigma$ with $\lambda>0$, where $\eta_\Sigma$ is Poincar\'e dual to the class of $\{z\}+\sym^{n-1}(\Sigma)$ (for any $z\in \Sigma$) and $\theta_\Sigma$ is the pullback by the Abel--Jacobi map of the principal polarisation of the Picard torus $\mathrm{Pic}^n(\Sigma)$. 
By a standard formula, one has 
\[ c_1(TM) =  \left ( n+ \frac{\chi(\Sigma)}{2} \right) \eta_\Sigma - \theta_\Sigma.   \]
Since $\theta_\Sigma$ integrates trivially on spheres, one finds that 
$M$ is monotone (negatively monotone) provided that $n+\chi(\Sigma)/2 > 0 $ ($<0$).

There is a rich supply of automorphisms of $\sym^n(\Sigma)$ coming from self-diffeomorphisms of $\Sigma$. Indeed, if we fix an area-form $\alpha$ on $\Sigma$, the group of area-preserving diffeomorphisms acts on the symplectic mapping class groups of the symmetric products $\sym^n(\Sigma)$; there is a natural homomorphism
\begin{equation}\label{kappa}  \kappa_{n,\lambda} \colon (\aut/\ham) (\Sigma ,\alpha) \to (\aut/\ham)(\sym^n (\Sigma) ,\omega_\lambda) . \end{equation} 
This homomorphism was first constructed by gauge-theoretic means by D. Salamon \cite{Salamon99}, but it can also be described in elementary terms \cite[II: Section 1]{PerutzI}.

When $\mu=\kappa_{\lambda,n}(\phi)$ for some area-preserving diffeomorphism $\phi  \in \aut(\Sigma)$, the fixed-point Floer homology $\HF(\mu)$ was conjectured by D. Salamon to be isomorphic to a version of Seiberg--Witten monopole Floer homology for $Y := \torus_\phi$ \cite{Salamon99}.  Salamon's conjecture can be rephrased, using the calculations of \cite{Perutz08}, as follows. 
\begin{conj}\label{SalConj}
Take a Riemann surface $\Sigma$ with a positive area form $\alpha$, normalised to have area 1, inducing a closed 2-form $\alpha_\phi$ on $ \torus_\phi$. Then there is an isomorphism of $\Lambda_{\Z/2}$-modules
\begin{equation} \label{HF vs HM}
 \HF_*(\kappa_{k,\lambda}(\phi) )\cong   
\bigoplus_{\langle c_1(\mathfrak{s} ),[\Sigma] \rangle = 2(n+1-g)}{\mathrm{HM}_*(T_\phi, \mathfrak{s}; c_\lambda )} \quad (\lambda>0),  
\end{equation} 
where on the right-hand side, a representative of 
\[c_\lambda:=-8\pi^2 (\lambda^{-1}+n) [\alpha_\phi]-32\pi^2c_1(\Tv \torus_\phi) \in H^2(\torus_\phi;\R)\] defines a perturbation of the Chern--Simons--Dirac functional. Moreover, there is an isomorphism (\ref{HF vs HM})  induced by a quasi-isomorphism of the underlying chain complexes, such that the quantum cap product by $\eta_{\Sigma}$ goes over to the `$U$-map' \cite[Chapter 1]{KronMrow} on $\mathrm{HM}_*$, up to chain homotopy.
\end{conj}
The conventions here are those of Kronheimer--Mrowka \cite{KronMrow} (note that the absence of decorations on the symbol $\mathrm{HM}$ reflects the absence of reducibles)

Returning to symplectic geometry of $M=\sym^n(\Sigma)$: it was shown in \cite{PerutzI} that one can allocate to any embedded circle $\gamma \subset \Sigma$ an $S^1$-fibred coisotropic hypersurface $ V_\gamma  \subset M$, up to Hamiltonian isotopy. Its reduced space is symplectomorphic to $N:=\sym^{n-1}(\Sigma_\gamma)$, where $\Sigma_\gamma$ is the surface obtained by surgery along $\gamma$.  The symplectic structure $\bar{\omega}_\lambda$ is a K\"ahler form representing  $\eta_{\Sigma_\gamma} + \lambda \theta_{\Sigma_\gamma}$. Thus there is an $S^1$-bundle
\[ \rho\colon V_\gamma \to N ; \quad  \rho^*\bar{\omega}_\lambda  =  \omega_\lambda|_{V_\gamma}.  \]
The Euler class $e(V)$ is Poincar\'e dual to the difference of the divisors $ p+\sym^{n-2}(\Sigma_\gamma)$ and $q+\sym^{n-2}(\Sigma_\gamma)$, where $p$ and $q$ lie in the two respective discs glued to $\Sigma\setminus \gamma$ to form  $\Sigma_\gamma$.

Suppose that $\gamma$ is a \emph{separating} circle. Then
\[N=\coprod_{k=0}^{n-1}{\sym^k(\Sigma_1)\times \sym^{n-1-k} (\Sigma_2)}, \]
where $\Sigma_\gamma$ has connected components $\Sigma_1$ and $\Sigma_2$.
Write $\wh{V}^k$ for the component of $\wh{V}_\gamma$ that lies over $\sym^k(\Sigma_1)\times \sym^{n-1-k} (\Sigma_2)$. It was shown in \cite[II, Section  4]{PerutzI} that strong negativity for $\wh{V}^k$ in $\sym^n(\Sigma)_-
\times \sym^k(\Sigma_1)\times \sym^{n-1-k} (\Sigma_2)$ can be arranged, for all $k\in \{0,\dots, n-1\}$, when $n\leq \frac{1}{2}\min {(g(\Sigma_1),g(\Sigma_2))}$.  From Theorem \ref{iso} we obtain the following.
\begin{prop}\label{sep}
Assume $\gamma$ separating, and let $\Sigma_1$ and $\Sigma_2$ be the connected components of $\Sigma_\gamma$. Let $\mu=\kappa_{k,\lambda}(\phi_1)\times \kappa_{n-1-k,\lambda}(\phi_2)$.  When $2n \leq\min {(g(\Sigma_1),g(\Sigma_2))}$, there is a canonical isomorphism
\[  \HF_*(\wh{V},(\id \times \mu)\wh{V}) \cong H \cone \left( \eta_{\Sigma_1}\otimes 1 - 1\otimes \eta_{\Sigma_2} \right)   \]
where $ \eta_{\Sigma_1}\otimes 1 - 1\otimes \eta_{\Sigma_2}$ acts as a degree $(-2)$ endomorphism of $\CF(\mu)$.
\end{prop}

\begin{prop}Suppose that Conjecture \ref{SalConj} is true. Then there is an isomorphism 
\begin{equation} 
\HF(\wh{V},(\id\times \mu)\wh{V})  \cong  \bigoplus_{c_1(\mathfrak{t} ) = 2(n+1-g)}{\mathrm{HM}_*(Y_1 \# Y_2, \mathfrak{t}; \tilde{c}_\lambda)},  
\end{equation}
where $\tilde{c}_\lambda$ is the class on $Y_1 \# Y_2$ induced by the $c_\lambda$-classes on the $Y_i$.
\end{prop}
\begin{pf}
We invoke an unpublished theorem  of T. Mrowka and P. Ozsv\'ath in Seiberg--Witten theory. (There are closely related results due to S. K. Donaldson and M. Furuta in Yang--Mills theory \cite{Don}.) The theorem is that there is a natural isomorphism
\[   HM_*(Y_1 \# Y_2 ) \cong \cone(U\otimes 1 - 1\otimes U)  \]
where $U\otimes 1 - 1\otimes U$ acts as a degree $(-2)$ endomorphism of $HM_*(Y_1)\otimes_{\Lambda_{\Z/2}} HM_*(Y_2)$. (We work over $\Lambda_{\Z/2}$ for convenient comparison with symplectic Floer homology.) The cobordism $X$ induces a homomorphism $HM_*(X)\colon HM_*(Y_1)\otimes_{\Lambda_{\Z/2}} HM_*(Y_2)\to HM_*(Y_1 \# Y_2 )$, and the composite $HM_*(X)\circ (U\otimes 1 - 1\otimes U)$ is nullhomotopic because the basepoints in $Y_1$ and $Y_2$, used in the definitions of the $U$-maps, are homotopic in $X$. Any version of monopole Floer homology can be used, and in particular closed 2-form perturbations are allowed,  provided that they extend over $X$. 

In our case, we can take a near-symplectic form on $X$ whose zero set is $\crit(F)$ (cf. \cite{PerutzI}). Comparing the Mrowka--Ozsv\'ath isomorphism with Proposition \ref{sep} establishes the result.
\end{pf}
The last proposition is closely related to the conjecture from \cite{PerutzI} that the `Lagrangian matching invariants' for broken fibrations defined there are equal to the Seiberg--Witten invariants of the underlying 4-manifolds.  

\subsection{Non-separating circles}
A simpler situation occurs when $\gamma\subset \Sigma$ is a non-separating circle.  Since $n+ \chi(\Sigma)/2 = (n-1) + \chi(\Sigma_\gamma)/2$, $N$ is (anti-)monotone precisely when $M$ is. The product $M_- \times N$ is also (anti-)monotone, with monotonicity constant $n+\chi(\Sigma)/2$, when this constant is positive (negative), and its minimal Chern number is 
\[ c_1^{\min}  (M_-\times N )= |n+ \frac{\chi(\Sigma)}{2}|= 
 |n-1 + \frac{\chi(\Sigma_\gamma) }{2} |. \]
The minimal Maslov index is  \cite{PerutzI}
\[  m^{\min}_{\wh{V} } =  2 c_1^{\min}(M_-\times N) =  |2n+\chi(\Sigma)|=  |2(n-1) + \chi(\Sigma_\gamma)|  . \]  
The circle-bundle $\rho\colon V_\gamma \to N$ has $e(V)=0$,  though there is no preferred trivialisation. Thus we can make the following deduction from Theorem \ref{Gysin} and its strongly negative counterpart.
\begin{prop}
Assume $\gamma$ non-separating. Take any $\mu\in \aut(N,\bar{\omega}_\lambda)$, for example, $\mu=\kappa_{n-1,\lambda}(\phi)$.  Fix a trivialisation of $\rho$. When $n >  g(\Sigma)$ or $2n \leq g(\Sigma)-1$,  there is then a canonical isomorphism of relatively graded modules,
\[  \HF(\wh{V},( \id\times\mu) \wh{V})\cong   \HF(\mu) \oplus \HF(\mu)[1]. \] 
\end{prop}

This sequence also has a topological counterpart. Here one starts with the connected 3-manifold $Y=\torus_\phi$, and forms its \emph{internal} connected sum by excising a pair of disjointly embedded discs and gluing in a 1-handle $S^2\times [0,1]$ to obtain a new 3-manifold $Y'$ which has the diffeomorphism type of $Y \# (S^1\times S^2)$.  When $Y$ is fibred over $S^1$, by a map $f$,  $Y'$ maps to $S^1$ by a Morse function $f$ with precisely two critical points, of indices 1 and 2, just as before.  Moreover, there is an elementary cobordism from $X$ from $Y$ to $Y'$, and a map $F$ (broken fibration) from $X$ to an annulus, interpolating between $f$ and $f'$, with an arc of critical points forming a non-degenerate critical manifold. The author does not know an explicit source in the literature for the behaviour of monopole Floer homology under 1-handle attachment, but the behaviour of another conjecturally isomorphic theory, Heegaard Floer homology is well-understood, and tallies with our result.

\subsection{Flat connections}
Just as symplectic geometry of symmetric products mimics Seiberg--Witten theory, so symplectic geometry of moduli spaces of flat connections models instanton Floer homology (cf. Dostoglou--Salamon \cite{DostSal}), at least in cases where reducible flat connections can be avoided. 

Consider the twisted character variety $M_g$ of a punctured genus $g\geq 1$ surface:
\[ M_g = \{ (A_1,B_1;\dots ; A_g,B_g)\in \SU(2)^{2g}: \prod_{i=1}^g[A_i,B_i]=-I \}/\SU(2). \] 
Here $\SU(2)$ acts by conjugation. Among many possible references for the geometry of $M_g$, we mention \cite{Thaddeus}, because it discusses the locus $V$ to be introduced momentarily.\footnote{There are also interesting unpublished studies by M. Callahan \cite{Callahan} and P. Seidel \cite{Seidel97} of the role of $V$ in symplectic monodromy problems.} The space $M_g$ is a canonically smooth and symplectic manifold of dimension $6g-6$. Assuming $g\geq 2$, the subspace $ V:= \{    A_g = I \}  \subset M_g $ 
has codimension 3, and there is a natural projection $\rho\colon V\to M_{g-1}$. This map is an $S^3$-bundle with vanishing Euler class (the locus 
$\{ B_g=A_g=I\} \subset M_g$ defines a section $S$). It is easily verified to be the reduction map of an $S^3$-fibred coisotropic submanifold (to see this one can use either the Atiyah--Bott symplectic structure, thinking of the $A_i$ and $B_i$ as holonomies of flat connections, or Goldman's finite-dimensional reformulation of it).

Because $M_g$ and $M_{g-1}$ are simply connected, one has \[m^{\min}_{\wh{V}}=2c_1^{\min}(M_{g-}\times M_{g-1})=4 = 3+1,\] 
which puts us on the borderline for applicability of the Gysin sequence. The global angular chain (actually a cycle if we take the zero-cycle to represent the Euler class) may be taken to be the section $S$. 

Regarding elements of $M_g$ as representations of $\pi_1(\Sigma)$, where $\Sigma$ is a genus $g$ surface equipped with a standard basis of curves, 
one sees that there is a natural homomorphism $\theta\colon \pi_0 \diff^+(\Sigma)\to \aut(M_g,\omega)$, $\phi \mapsto (\phi^{-1})^*$.  According to Theorem \ref{borderline}, one has
\[  \HF_*(\wh{V}, \theta(\phi)\wh{V})\cong \cone(\nu_S t^a \, \id)  
\] 
for some $a>0$. Thus  $\HF_*(\wh{V}, \theta(\phi)\wh{V})$ is either always (i.e., for all $\phi$) zero, or else it is always  isomorphic, as a relatively graded module, to $\HF_*(\theta(\phi))\oplus \HF_*(\theta(\phi))[3]$. 

However, $\HF_*(\wh{V},\wh{V})$ is non-zero for any $g\geq 2$: by a theorem of Albers \cite[Corollary 2.11]{Albers}, a monotone Lagrangian $L\subset M^{2n}$ of minimal Maslov index $\geq 2$ can only have vanishing Floer self-homology if $[L] =0 \in H_n (M;\Z/2)$; but $[V]$ is non-zero in $H_*(M_g;\Z/2)$, hence $[\wh{V}]\neq 0 \in H_*(M_g \times M_{g-1})$. We make the following conclusion. 

\begin{thm}
A choice of section of $\rho\colon V\to M_{g-1}$ determines isomorphisms 
\[ \HF_*(\wh{V},\wh{V})\cong \HF_*(\theta(\phi))\oplus \HF_*(\theta(\phi))[3]\]
for all $\phi \in \diff^+(\Sigma)$.
\end{thm}

This appears to tally with the behaviour of $\SO(3)$-instanton Floer homology under internal connected sum, though we do not attempt to make the connection precise. It would be interesting to find a symplectic counterpart for the behaviour of external connected sums in instanton theory, as described in \cite{Fukaya,Don}.

\end{document}